\documentclass[OJMO,manuscript]{cedram}

\usepackage[all]{xypic}

\usepackage{bm}
\usepackage{booktabs, makecell}
\usepackage{natbib}
\usepackage{tikz}
\usepackage{mathtools}
\usetikzlibrary{positioning}
\usepackage{algorithm, algorithmicx, algpseudocode}
\algnewcommand{\IfThen}[2]{\State\algorithmicif\ #1\ \algorithmicthen\ #2}
\algnewcommand{\IfThenElse}[3]{\State\algorithmicif\ #1\ \algorithmicthen\ #2 \algorithmicelse\ #3}

\newcommand{\reals}{{\mathbb R}}

\newcommand{\idm}{\operatorname{I}}

\DeclareMathOperator*{\argmin}{arg\,min}

\newcommand{\state}{x}
\newcommand{\states}{{\bm x}}
\newcommand{\ctrl}{u}
\newcommand{\ctrls}{{\bm u}}

\newcommand{\auxstate}{y}
\newcommand{\auxstates}{{\bm \auxstate}}
\newcommand{\auxxstate}{z}
\newcommand{\auxxstates}{{\bm \auxxstate}}
\newcommand{\auxctrl}{v}
\newcommand{\auxctrls}{{\bm \auxctrl}}
\newcommand{\auxxctrl}{w}
\newcommand{\auxxctrls}{{\bm \auxxctrl}}

\newcommand{\costate}{\lambda}
\newcommand{\costates}{{\bm \costate}}

\newcommand{\diffstate}{\auxstate}

\newcommand{\diffctrl}{\auxctrl}
\newcommand{\diffctrls}{\bm \auxctrl}

\newcommand{\dualvar}{\mu}
\newcommand{\dualvars}{{\bm \mu}}

\newcommand{\horizon}{\tau}

\newcommand{\dimstate}{{n_\state}}
\newcommand{\dimctrl}{{n_\ctrl}}

\newcommand{\dyn}{f}

\newcommand{\contdyn}{\mathrm{f}}
\newcommand{\auxdyn}{\phi}

\newcommand{\concatdyn}{F}

\newcommand{\traj}{{\dyn^{[\horizon]}}}

\newcommand{\trajfunc}{{\traj}}

\newcommand{\augtrajfunc}{{g}}
\newcommand{\augtraj}{{g}}

\newcommand{\cost}{h}
\newcommand{\statecost}{c}

\newcommand{\obj}{\mathcal{J}}

\newcommand{\initstate}{\bar \state_0}

\newcommand{\costogo}{c}

\renewcommand{\H}{P}
\newcommand{\h}{p}
\newcommand{\G}{Q}
\newcommand{\g}{q}
\newcommand{\A}{A}
\newcommand{\B}{B}

\newcommand{\J}{J}
\renewcommand{\j}{j}
\newcommand{\jcst}{j^0}

\newcommand{\Q}{Q}
\newcommand{\q}{q}
\newcommand{\R}{R}

\newcommand{\diff}{\delta}
\newcommand{\lin}{\ell}
\newcommand{\qua}{q}
\newcommand{\model}{m}

\newcommand{\nxt}{{\operatorname{next}}}

\newcommand{\var}{\ctrls}
\newcommand{\fixedvar}{\var}
\newcommand{\currvar}{{\var^{(k)}}}

\newcommand{\dimvar}{{\horizon \dimctrl}}
\newcommand{\varaux}{\auxctrls}

\newcommand{\Oracle}{\operatorname{Oracle}}
\newcommand{\stepsize}{\gamma}
\newcommand{\stepsizescaled}{\bar \stepsize}
\newcommand{\reg}{\nu}

\newcommand{\dec}{\operatorname{dec}}
\newcommand{\inc}{\operatorname{inc}}
\newcommand{\prev}{\operatorname{prev}}

\newcommand{\Bell}{\operatorname{BP}}
\newcommand{\BellLQ}{\operatorname{LQBP}}
\newcommand{\BellL}{\operatorname{LBP}}
\newcommand{\Bellapprox}{\widehat{\Bell}}

\newcommand{\Roll}{\operatorname{Roll}}

\newcommand{\Forward}{\operatorname{Forward}}
\newcommand{\Backward}{\operatorname{Backward}}

\newcommand{\linesearch}{\operatorname{LineSearch}}
\newcommand{\DynProg}{\operatorname{DynProg}}

\newcommand{\infeasible}{\texttt{infeasible}}

\newcommand{\statexp}{z}
\newcommand{\yaw}{\theta}

\newcommand{\Pol}{\operatorname{Pol}}

\newcommand{\CheckSubPb}{\operatorname{CheckSubProblem}}
\newcommand{\valid}{\operatorname{valid}}
\newcommand{\True}{\operatorname{True}}
\newcommand{\False}{\operatorname{False}}

\newcommand{\coderef}{{\url{https://github.com/vroulet/ilqc}}}

\definecolor{ForestGreen}{RGB}{34,139,34}
\newcommand{\revised}[1]{#1}

\title[\hfill Iterative Linear Quadratic Optimization]{Iterative Linear Quadratic Optimization for Nonlinear Control: \\	Differentiable Programming Algorithmic Templates}

\author{\firstname{Vincent} \lastname{Roulet}}
\address{Google Brain, Seattle, USA\footnote{Work completed at the University of Washington before joining Google.}}
\email{vroulet@google.com}

\author{\firstname{Siddhartha} \lastname{Srinivasa}}
\address{Paul G. Allen School of Computer Science and Engineering, University of Washington, Seattle, USA}
\email{siddh@cs.washington.edu}

\author{\firstname{Maryam} \lastname{Fazel}}
\address{Department of Electrical and Computer Engineering, University of Washington, Seattle, USA}
\email{mfazel@uw.edu}

\author{\firstname{Zaid} \lastname{Harchaoui}}
\address{Department of Statistics University of Washington, Seattle, USA}
\email{zaid@uw.edu}
\thanks{This work was supported by NSF DMS-1839371, DMS-2134012, CCF-2019844, CIFAR-LMB, NSF TRIPODS II DMS-2023166 and faculty research awards.} 


\begin{abstract}
Iterative optimization algorithms depend on access to information about the
objective function. In a differentiable programming framework, this information,
such as gradients, can be automatically derived from the computational graph. We
explore how nonlinear control algorithms, often employing linear and/or
quadratic approximations, can be effectively cast within this framework.
Our approach illuminates shared components and differences between gradient
descent, Gauss-Newton, Newton, and differential dynamic programming methods in
the context of discrete time nonlinear control. Furthermore, we present
line-search strategies and regularized variants of these algorithms, along with
a comprehensive analysis of their computational complexities. We study the
performance of the aforementioned algorithms on various nonlinear control
benchmarks, including autonomous car racing simulations using a simplified car
model. All implementations are publicly available in a package coded in a
differentiable programming language.
\end{abstract}

\begin{document}
	\maketitle

	\section{Introduction}
	We consider nonlinear control problems in discrete time with finite horizon,
i.e., problems of the form
\begin{align}
	\min_{\substack{\state_0, \ldots, \state_\horizon \in \reals^\dimstate \\ \ctrl_0 \ldots, \ctrl_{\horizon-1} \in \reals^\dimctrl}}
	 & \sum_{t=0}^{\horizon-1} \cost_t(\state_t, \ctrl_t)  + \cost_\horizon(\state_\horizon) \label{eq:discrete_ctrl_pb}\\
	\mbox{subject to} \quad & \state_{t+1} = \dyn_t(\state_t, \ctrl_t), \quad 
  \mbox{for} \ t \in \{0, \ldots, \horizon-1\},  \quad \state_0 = \initstate, \nonumber
\end{align}
where at time $t$, $\state_t \in \reals^{\dimstate}$ is the state of the system,
$\ctrl_t \in \reals^\dimctrl$ is the control applied to the system, $\dyn_t:
\reals^\dimstate \times \reals^\dimctrl \rightarrow \reals^\dimstate$ is the
discrete dynamic, $\cost_t: \reals^\dimstate \rightarrow \reals$ is the cost on
the state and control variables  and $\initstate \in \reals^\dimstate$ is a
given fixed initial state. Problem~\eqref{eq:discrete_ctrl_pb} is entirely
determined by the initial state and the controls.

Problems of the form~\eqref{eq:discrete_ctrl_pb} have been tackled in various
ways, from direct approaches using nonlinear optimization
\citep{betts2010practical, wright1990solution, wright1991partitioned,
wright1993interior, de1988differential, dunn1989efficient, rao1998application}
to convex relaxations using semi-definite
optimization~\citep{boyd1997semidefinite}. Numerous packages exist for such
problems such as CasAdi~\citep{andersson2018cassadi},
Pyomo~\citep{bynum2021pyomo}, JumP~\citep{dunning2017jump},
IPOPT~\citep{wachter2006implementation}, or SNOPT~\citep{gill2005snopt},
Crocoddyl~\citep{jallet2023proxddp}, acados~\citep{Verschueren2021}. A popular
approach of the former category proceeds by computing at each iteration the
linear quadratic regulator associated with a linear quadratic approximation of
the problem around the current candidate
solutions~\citep{jacobson1970differential,li2007iterative, sideris2005efficient,
tassa2012synthesis}. The computed feedback policies are then applied either
along the linearized dynamics or along the original dynamics to output a new
candidate solution. Such canonical nonlinear control algorithms efficiently
incorporate second-order information into the optimization procedure by
exploiting the dynamical structure of the problem. This approach lends itself to
an integration in a differentiable programming framework to extend this paradigm
beyond first-order oracles. 

Differentiable programming consists of 
the implementation of functions in a programming language that enables access to
derivatives of these functions by automatic
differentiation~\citep{baur1983complexity,  rumelhart1985learning,
lecun1988theoretical, schmidhuber1990making, gilbert1992automatic, werbos1994roots,
griewank2008evaluating, baydin2018automatic, bolte2020mathematical, tensorflow2015-whitepaper,
paszke2017automatic}. Automatic differentiation itself has roots in the control
literature, and its use is pervasive in numerous
domains~\citep{griewank2008evaluating}, in particular deep
learning~\citep{zhang2021dive, goodfellow2016deep}. Canonical nonlinear control
algorithms incorporating second order information
can also be integrated in reinforcement learning pipelines
~\citep{recht2019tour, kakade2020information}, and may then benefit from a
differentiable programming viewpoint to isolate their underlying principles.
These algorithms have indeed generally be presented through linear algebraic
manipulations instantiated separately for each algorithm, which hinder a global
perspective~\citep{murray1984differential, de1988differential, li2007iterative,
sideris2005efficient, tassa2012synthesis}. 

The motivation of this work is to cast all such algorithms in a common
differentiable programming viewpoint to delineate the discrepancies between the
different algorithms and identify the common subroutines. We review the
implementation of (i) a Gauss-Newton method~\citep{sideris2005efficient}, a.k.a.
Iterative Linear Quadratic Regulator (ILQR), (ii) a Newton
method~\citep{de1988differential, liao1991convergence, dunn1989efficient}, (iii)
a differential dynamic programming approach based on linear approximations of
the dynamics and quadratic approximations of the costs, a.k.a. iterative Linear
Quadratic Regulator (iLQR)~\citep{tassa2012synthesis}, (iv) a differential
dynamic programming approach based on quadratic approximations of both dynamics
and costs, usually simply called DDP~\citep{jacobson1970differential}, and
consider regularized variants of the aforementioned algorithms with their
corresponding line searches. In turn, the differentiable programming viewpoint
informs efficient handling of memory by appropriate check-pointing. An
extended related work discussion is in
Appendix~\ref{app:related_work}.\\[-1em]

{\bf Outline.}
In Sec.~\ref{sec:lin_quad} we recall how linear quadratic control problems are
solved by dynamic programming and used as a
building block for nonlinear control algorithms. The implementation of classical
optimization oracles such as a gradient step, a Gauss-Newton step, or a Newton
step is presented in Sec.~\ref{sec:classical_optim}. Sec.~\ref{sec:ddp} details
the rationale and implementation of differential dynamic programming
approaches. Sec.~\ref{sec:comput_cplxity} presents the computational
complexities of each oracle in terms of space and time complexities in a
differentiable programming framework. All algorithms are tested on several
synthetic problems in Sec.~\ref{sec:exp}: swinging-up a fixed pendulum, and 
autonomous car racing with simple dynamics. Code is available at \coderef.

Appendix~\ref{app:notations},~\ref{app:related_work},~\ref{app:proofs},~\ref{app:line_searches}
detail notations, related work, proofs and line-search procedures respectively.
A summary of all algorithms with detailed pseudocode and computational schemes is given in
Appendix~\ref{app:summary}. Alternative implementations using check-pointing and 
different linear algebra solvers are presented in Appendix~\ref{app:cplxity} 
and~\ref{app:sparse_solvers} respectively.
Experimental setups and additional experiments are detailed in
Appendix~\ref{app:exp_details} and~\ref{app:exp_sup}.\\[-1em]

{\bf Notation.}
For a sequence of vectors $\state_1, \ldots, \state_{\horizon} \in
\reals^\dimstate$, we denote by semicolons their concatenation s.t. $\states =
(\state_1;\ldots;\state_\horizon) \in \reals^{\horizon \dimstate}$. For a
function $f:\reals^d\rightarrow \reals^n$, we denote by $\nabla f(x) \coloneqq
(\partial_{x_i} f_j(x))_{\substack{1\leq i\leq d, 1\leq  j\leq n}} \in \reals^{d
\times n}$ the transpose of the Jacobian of $f$ on
$x$. For a  function $f:\reals^d \times \reals^p \rightarrow  \reals^n$, we
denote  for $x\in \reals^d$, $y\in \reals^p$, $\nabla_x f(x, y)= (\partial_{x_i}
f_j(x, y) )_{\substack{1\leq i\leq d, 1\leq  j\leq n}} \in \reals^{d \times n}$
the partial transpose Jacobian of $f$ w.r.t. $x$ on $(x, y)$.  

For a multivariate function $f:\reals^d \rightarrow \reals^n$ composed of
coordinates $f_j:\reals^d \rightarrow \reals$ for $j \in \{1, \ldots, n\}$, we
denote its Hessian  $x\in \reals^d$ as a tensor $\nabla^2 f(x) \coloneqq (\nabla^2
f_1(x),\ldots,\nabla^2f_n(x)) \in \reals^{d\times d \times n}$. For a
multivariate function $f: \reals^d \times \reals^p \rightarrow \reals^n$
composed of coordinates $f_j:\reals^d \times \reals^p \rightarrow \reals$ for $j
\in \{1, \ldots, n\}$, we decompose its Hessian on $x \in \reals^d$, $y\in
\reals^p$ by defining, e.g.,    $\nabla_{xx}^2 f( x,  y) = (\nabla_{xx}^2 f_1(
x,  y), \ldots, \nabla_{xx}^2 f_n( x,  y)) \in \reals^{d \times d \times n}$.
The quantities $ \nabla_{yy}^2 f( x,  y) \in \reals^{p \times p \times n},
\nabla_{xy}^2 f( x,  y) \in \reals^{d \times p \times n},  \nabla_{yx}^2 f( x,
y) \in \reals^{p \times d \times n}$ are defined similarly. 

For a function $f:\reals^d \rightarrow \reals^ n$, and $x\in \reals^d$, we
define the finite difference expansion of $f$ around $x$, the linear expansion
of $f$ around $x$ and the quadratic expansion of $f$ around $x$ as,
respectively,
\begin{align}
\diff_f^x(y) \coloneqq f(x+y) - f(x), \qquad 
\lin_f^\state(y) \coloneqq \nabla f(x)^\top y, \label{eq:lin_approx}\qquad 
\qua_f^\state(y)  \coloneqq \nabla f(x)^\top y +\frac{1}{2} \nabla^2f(x)[y, y, \cdot].
\end{align}
The linear and quadratic approximations of $f$ around $x$  are then  $f(x+
\diffstate) \approx f(x) +\lin_f^\state(\diffstate)$ and $f(x+ \diffstate)
\approx f(x) + \qua_f^\state(\diffstate)$ respectively. Tensor notations, such
as $\nabla^2f(x)[y, y, \cdot]$, inspired from~\citep{nesterov2018lectures}, are
detailed in Appendix~\ref{app:notations}.

	\section{From Linear Quadratic Control Problem to Nonlinear Control Algorithm}\label{sec:lin_quad}
	Algorithms for nonlinear control problems revolve around
solving linear quadratic control problems  by dynamic programming. Therefore, we start
by recalling the rationale of dynamic programming and how discrete time control
problems with linear dynamics and quadratic costs can be solved by dynamic
programming.

\subsection{Dynamic Programming}
The idea of dynamic programming is to decompose dynamical problems such
as~\eqref{eq:discrete_ctrl_pb} into a sequence of nested subproblems defined  by
the \emph{cost-to-go} $c_t$, from $ \state_t$ at time $t \in \{0, \ldots, \horizon-1\}$: 
\begin{align}
	\costogo_t(\state_t) \coloneqq \min_{\substack{\ctrl_t,\ldots, \ctrl_{\horizon-1} \in \reals^\dimctrl\\\auxstate_t,\ldots, \auxstate_\horizon \in \reals^\dimstate}} \quad & \sum_{s=t}^{\horizon-1} \cost_{s}(\auxstate_{s}, \ctrl_{s}) + \cost_\horizon(\auxstate_\horizon)  \nonumber
	\\
	\mbox{subject to} \quad & \auxstate_{s+1} =  \dyn_{s}(\auxstate_{s}, \ctrl_{s}) \quad \mbox{for} \ s \in \{t, \ldots, \horizon-1\}, \nonumber \quad \auxstate_{t} =  \state_{t}. \nonumber
\end{align}
The cost-to-go from $\state_\horizon$ at time $\horizon$ is simply the last
cost, namely,
$
\costogo_\horizon(\state_\horizon) = \cost_\horizon(\state_\horizon),
$
and the original problem~\eqref{eq:discrete_ctrl_pb} amounts to compute
$\costogo_0(\bar \state_0)$. The cost-to-go functions define nested subproblems
that are linked for $t\in \{0, \ldots, \horizon-1\}$ by \emph{Bellman's
equation}~\citep{bellman1971introduction}
\begin{align}\label{eq:bellman}
	\costogo_{t}( \state_t)  & = \min_{\ctrl_t \in \reals^\dimctrl}  \cost_t(\state_t, \ctrl_t) + 
	\min_{\substack{\ctrl_{t+1},\ldots, \ctrl_{\horizon-1} \in \reals^\dimctrl
			\\\auxstate_{t+1},\ldots, \auxstate_\horizon \in \reals^\dimstate}}    
		\sum_{s=t+1}^{\horizon-1} \cost_{s}(\auxstate_{s}, \ctrl_{s}) + \cost_\horizon(\auxstate_\horizon) \nonumber\\
		& \hspace{105pt} \mbox{subject to} \hspace{20pt}   \auxstate_{s+1} =  \dyn_{s}(\auxstate_{s}, \ctrl_{s}) \ \mbox{for}\ s \in \{t+1, \ldots, \horizon-1\}, \ \auxstate_{t+1} =  \dyn_{t}(\state_{t}, \ctrl_{t}) \nonumber\\
	& = \min_{\ctrl_t \in \reals^\dimctrl}  \cost_t(\state_t, \ctrl_t) + \costogo_{t+1}(\dyn_t( \state_t, \ctrl_t)).
\end{align}
The optimal control at time $t$ from state $\state_t$ is given by
$\ctrl_t=\pi_t(\state_t)$, where $\pi_t$, called a \emph{policy}, is given by
\begin{equation}\nonumber
\pi_t( \state_t) \coloneqq \argmin_{\ctrl_t \in \reals^\dimctrl} \left\{ \cost_t(\state_t, \ctrl_t) + \costogo_{t+1}(\dyn_t( \state_t, \ctrl_t))\right\}.
\end{equation}
Define the procedure that back-propagates ($\Bell$) the cost-to-go functions as
\[
\Bell:
	\dyn_t, \cost_t, \costogo_{t+1} \rightarrow \left(\begin{array}{c}
		\costogo_{t}: \state\rightarrow	\min_{\ctrl \in \reals^\dimctrl}  \left\{ \cost_t(\state, \ctrl) + \costogo_{t+1}(\dyn_t( \state, \ctrl))\right\},\\
		\pi_t :\state \rightarrow	 \argmin_{\ctrl\in \reals^\dimctrl} \left\{ \cost_t(\state, \ctrl) + \costogo_{t+1}(\dyn_t( \state, \ctrl))\right\}
	\end{array}\right).
\] 
A dynamic programming approach, formally described in Algo.~\ref{algo:dyn_prog},
solves problems of the form~\eqref{eq:discrete_ctrl_pb} as follows.
\begin{enumerate}
	\item Compute recursively the cost-to-go functions $\costogo_t$ for
	$t=\horizon, \ldots, 0$ using  Bellman's equation~\eqref{eq:bellman}, i.e.,
	compute from $\costogo_\horizon= \cost_\horizon$, 
	\[
	\costogo_t, \pi_t = \Bell(\dyn_t, \cost_t, \costogo_{t+1}) \quad \mbox{for} \ t \in \{\horizon-1, \ldots, 0\},
	\]
	and record at each step the policies $\pi_t$.
	\item Unroll the optimal trajectory that starts from time 0 at $\bar
	\state_0$, follows the dynamics $\dyn_t$, and uses at each step the optimal
control given by the computed policies, that is, starting from $\state_0^* =
	\initstate$, compute
	\begin{equation}\label{eq:true_roll_out}
	\ctrl^*_t= \pi_t(\state_t^*), \qquad \state_{t+1}^*= \dyn_t(\state_t^*, \ctrl^*_t) \qquad \mbox{for $t=0, \ldots, \horizon-1$.}
	\end{equation}
\end{enumerate}
The resulting command $\ctrls^* = (\ctrl_0^*; \ldots; \ctrl_{\horizon-1}^*)$ and
trajectory $\states^* = (\state_1^*; \ldots; \state_\horizon^*)$ are then
optimal for problem~\eqref{eq:discrete_ctrl_pb}. 
In the following, the 
dynamic programming ($\DynProg$) procedure, detailed\footnote{For
ease of reference and comparisons, all procedures,
algorithms, and computational schemes are grouped in Appendix~\ref{app:summary}.} in
Algo.~\ref{algo:dyn_prog} in Appendix~\ref{app:summary}, is denoted
\begin{equation}\label{eq:dyn_prog_algo}
  \DynProg: (\dyn_t)_{t=0}^{\horizon-1},(\cost_t)_{t=0}^\horizon , \initstate, \Bell \rightarrow \ctrl_0^*, \ldots, \ctrl_{\horizon-1}^*.
\end{equation}

The bottleneck of the approach is the ability to solve Bellman's
equation~\eqref{eq:bellman}, i.e.,  having access to the procedure $\Bell$
defined above. 

\subsection{Linear Dynamic, Quadratic Cost}
For linear dynamics and quadratic costs, problem~\eqref{eq:discrete_ctrl_pb}
takes the form
\begin{align*}
		\min_{\substack{\state_0, \ldots, \state_\horizon \in \reals^\dimstate \\ \ctrl_0 \ldots, \ctrl_{\horizon-1} \in \reals^\dimctrl}}
	& \sum_{t=0}^{\horizon-1} \left( \frac{1}{2}\state_t^\top \H_t \state_t+  \frac{1}{2}\ctrl_t^\top \G_t \ctrl_t + \state_t^\top \R_t \ctrl_t+ \h_t^\top \state_t  +  \g_t^\top \ctrl_t\right) +\frac{1}{2}\state_\horizon^\top \H_\horizon \state_\horizon+ \h_\horizon^\top \state_\horizon \label{eq:lin_quad_ctrl_pb}\\
	\mbox{subject to} \quad & \state_{t+1} = A_t \state_t + B_t\ctrl_t,  \quad \mbox{for} \ t \in \{0, \ldots, \horizon-1\}, \quad \state_0 = \initstate. \nonumber
\end{align*}
Namely, we have $\cost_t(\state_t, \ctrl_t)  =  \frac{1}{2}\state_t^\top \H_t
\state_t+  \frac{1}{2}\ctrl_t^\top \G_t \ctrl_t + \state_t^\top \R_t \ctrl_t+
\h_t^\top \state_t  +  \g_t^\top \ctrl_t $ and $\dyn_t(\state_t, \ctrl_t)= A_t
\state_t + B_t\ctrl_t$. In that case, under appropriate conditions on the
quadratic functions,  Bellman's equation~\eqref{eq:bellman} can be solved
analytically through a linear quadratic back-propagation ($\BellLQ$) as recalled
in Lemma~\ref{lem:lin_quad_exact}.  Note that the operation~$\BellLQ$ defined
in~\eqref{eq:bellman_lq} amounts to computing the Schur complement of a block of
the Hessian  of  the quadratic $\state, \ctrl \rightarrow \qua_t(\state, \ctrl)
+ \costogo_{ t+1}(\lin_t(\state, \ctrl))$, namely, the block corresponding to
the Hessian w.r.t. the control variables (see, e.g.,~\citep[Appendix
A.5.5]{boyd2004convex}). The proofs of Lemma~\ref{lem:lin_quad_exact} and
Corollary~\ref{corr:lin_quad} are standard and are given in
Appendix~\ref{app:proofs}.

\begin{lemm}\label{lem:lin_quad_exact}
For linear functions $\lin_t$ and quadratic functions $\qua_t, \costogo_{t+1}$ s.t. $\qua_t(\state, \cdot) + \costogo_{t+1}(\lin_t(\state, \cdot))$ is
strongly convex for any $x$, the procedure
\begin{equation}
\BellLQ	: (\lin_t, \qua_t, \costogo_{t+1})  \rightarrow \left(\begin{array}{c}
    \costogo_t:\state\rightarrow\min_{\ctrl \in \reals^\dimctrl} \left\{ \qua_t(\state, \ctrl) + \costogo_{t+1}(\lin_t(\state, \ctrl))\right\}  \\
    \pi_t :\state\rightarrow \argmin_{\ctrl \in \reals^\dimctrl}  \left\{ \qua_t(\state, \ctrl) + \costogo_{t+1}(\lin_t(\state, \ctrl))\right\}
\end{array}\right),
\label{eq:bellman_lq}
\end{equation}
can be implemented analytically as detailed in Algo.~\ref{algo:BellLQ}.
\end{lemm}
If  problem~\eqref{eq:discrete_ctrl_pb} consists of linear dynamics and
quadratic costs that are strongly convex w.r.t. the control variable, the
procedure $\BellLQ$ can be applied iteratively in a dynamic programming approach
to give the solution of the problem, as formally stated in
Corollary~\ref{corr:lin_quad}. 
\begin{coro}\label{corr:lin_quad}
Consider problem~\eqref{eq:discrete_ctrl_pb}
such that for all $t\in \{0, \ldots, \horizon-1\}$, $\dyn_t$ is linear,
$\cost_t$ is convex quadratic with $\cost_t(\state, \cdot)$ strongly convex
for any $\state$, and $\cost_\horizon$ is convex quadratic. Then, the solution
of problem~\eqref{eq:discrete_ctrl_pb} is given by 
\[
\ctrls^* = \DynProg((\dyn_t)_{t=0}^{\horizon-1},(\cost_t)_{t=0}^\horizon , \initstate, \BellLQ),
\]
with $\DynProg$ and $\BellLQ$
implemented in Algo.~\ref{algo:dyn_prog} and Algo.~\ref{algo:BellLQ} respectively.
\end{coro}

\subsection{Nonlinear Control Algorithm Example}\label{ssec:nonlin_algo}
Nonlinear control algorithms based on nonlinear optimization use linear or
quadratic approximations of the dynamics and the costs at a current candidate
sequence of controllers to apply a dynamic programming procedure to the
resulting problem \citep{bellman1971introduction, dunn1989efficient,
sideris2005efficient, li2007iterative, tassa2012synthesis}. For example, the
Iterative Linear Quadratic Regulator (ILQR) algorithm uses linear approximations
of the dynamics and quadratic approximations of the
costs~\citep{li2007iterative}. Each iteration of the ILQR algorithm is composed
of the three steps below illustrated in Fig.~\ref{fig:ilqr}.\\

{\bf Iterative Linear Quadratic Regulator Iteration.}
\begin{enumerate}
 \item \underline{Forward pass:} Given a set of control variables $\ctrl_0,
 \ldots, \ctrl_{\horizon-1}$, compute the trajectory $
 \state_1,\ldots,\state_\horizon$ as $\state_{t+1} = \dyn_t(\state_t, \ctrl_t)$
 starting from $\state_0 = \initstate$, and the associated costs
 $\cost_t(\state_t, \ctrl_t), \cost_\horizon(\state_\horizon)$,  for $t \in \{0,
 \ldots, \horizon-1\}$. Record along the computations, i.e., for $t\in \{0,
 \ldots, \horizon-1\}$, the gradients of the dynamics and the gradients and
 Hessians of the costs.
 \item \underline{Backward pass:} Compute the optimal policies associated with
 the linear quadratic control problem
 	\begin{align}
 	\min_{\substack{\diffstate_0, \ldots \diffstate_\horizon \in \reals^\dimstate \\ \diffctrl_0,\ldots, \diffctrl_{\horizon-1} \in \reals^\dimctrl} } &
 	\sum_{t=0}^{\horizon-1}  
 	\left(\frac{1}{2}\diffstate_t^\top \H_t \diffstate_t
 	+  \frac{1}{2}\diffctrl_t^\top \G_t \diffctrl_t + \diffstate_t^\top \R_t \diffctrl_t
 	+ \h_t^\top \diffstate_t  +  \g_t^\top \diffctrl_t
 	\right)
 	+ \frac{1}{2}\diffstate_\horizon^\top \H_\horizon \diffstate_\horizon 
 	+ \h_\horizon^\top \diffstate_\horizon\nonumber\\
 	\mbox{subject to} \quad & \diffstate_{t+1} = A_t \diffstate_t + B_t\diffctrl_t, \quad \mbox{for} \ t \in \{0, \ldots, \horizon-1\}, \quad \diffstate_0 = 0, \nonumber\\
 	\mbox{where} \quad & \H_t = \nabla_{\state_t\state_t}^2 \cost_t(\state_t, \ctrl_t) \ \  
 	\G_t =   \nabla_{\ctrl_t\ctrl_t}^2 \cost_t(\state_t, \ctrl_t)\ \ 
 	\R_t = \nabla_{\state_t\ctrl_t}^2 \cost_t(\state_t, \ctrl_t)\ \nonumber \\
 	&\hspace{1.5pt} \h_t = \nabla_{\state_t} \cost_t(\state_t, \ctrl_t)\  \hspace{14pt}
 	\g_t = \nabla_{\ctrl_t} \cost_t(\state_t, \ctrl_t)\  \nonumber\\
 	&\hspace{-0.5pt} \A_t = \nabla_{\state_t} \dyn_t(\state_t, \ctrl_t)^\top \ \hspace{5pt}
 	\B_t = \nabla_{\ctrl_t} \dyn_t(\state_t, \ctrl_t)^\top. \nonumber
 \end{align}
The problem above can be written compactly as 
 			\begin{align}
	\min_{\substack{\diffstate_0, \ldots \diffstate_\horizon \in \reals^\dimstate \\ \diffctrl_0,\ldots, \diffctrl_{\horizon-1} \in \reals^\dimctrl} } &
	\sum_{t=0}^{\horizon-1} \qua_{\cost_t}^{\state_t, \ctrl_t}(\diffstate_t, \diffctrl_t) + \qua_{\cost_\horizon}^{\state_\horizon}(\diffstate_\horizon) \label{eq:LQR_gauss_newton} \\
	\mbox{\textup{subject to}} \quad & \diffstate_{t+1} = \lin_{\dyn_t}^{\state_t, \ctrl_t}(\diffstate_t, \diffctrl_t), \quad \mbox{for} \ t \in \{0, \ldots, \horizon-1\}, \quad  \diffstate_0 = 0, \hspace{65pt} \nonumber
\end{align}
where $ \qua_{\cost_\horizon}^{\state_\horizon}(\diffstate_\horizon)
=\frac{1}{2}\diffstate_\horizon^\top \H_\horizon \diffstate_\horizon +
\h_\horizon^\top \diffstate_\horizon$ and $ \qua_{\cost_t}^{\state_t,
\ctrl_t}(\diffstate_t, \diffctrl_t)= \frac{1}{2}\diffstate_t^\top \H_t
\diffstate_t +  \frac{1}{2}\diffctrl_t^\top \G_t \diffctrl_t + \diffstate_t^\top
\R_t \diffctrl_t + \h_t^\top \diffstate_t  +  \g_t^\top \diffctrl_t$ are the
quadratic expansions of the costs and $ \lin_{\dyn_t}^{\state_t,
\ctrl_t}(\diffstate_t, \diffctrl_t) =  A_t \diffstate_t + B_t\diffctrl_t$ is the
linear expansion of the dynamics, both expansions being defined around the
current sequence of controls and associated trajectory. The optimal policies
associated to this problem are obtained by computing recursively, starting from
$\costogo_\horizon = \qua_{\cost_\horizon}^{\state_\horizon}$, 
\[
\costogo_t, \pi_t = \BellLQ(\lin_{\dyn_t}^{\state_t, \ctrl_t}, \qua_{\cost_t}^{\state_t, \ctrl_t}, \costogo_{t+1}) \quad \mbox{for} \ t\in \{\horizon-1, \ldots, 0\},
\]
where $\BellLQ$ presented in Algo.~\ref{algo:BellLQ} outputs affine policies of
the form $\pi_t:\auxstate_t \rightarrow K_t\auxstate_t + k_t$.
\item \underline{Roll-out pass:} Define the set of candidate policies as
$\{\pi_t^\stepsize: \diffstate \rightarrow K_t \diffstate  + \stepsize k_t \
\mbox{ for $\stepsize \geq 0$} \}$. The next sequence of controllers is then
given as $\ctrl_t^\nxt = \ctrl_t + \diffctrl^\stepsize_t$ , where
$\diffctrl_t^\stepsize$ is given by rolling out the policies $\pi_t^\stepsize$
from $\diffstate^\stepsize_0=0$ along the linearized dynamics as 
\[
\diffctrl_t^\stepsize = \pi_t^\stepsize(\diffstate_t^\stepsize), \quad  \diffstate_{t+1} = \lin_{\dyn_t}^{\state_t, \ctrl_t}(\diffstate_t^\stepsize, \diffctrl_t^\stepsize), \quad \mbox{for} \ t \in \{0, \ldots, \tau\}
\]
for $\stepsize$ found by a line-search such that 
$
\sum_{t=0}^{\horizon-1} \left(\cost_t(\state_t + \diffstate_t^\stepsize, \ctrl_t+ \diffctrl_t^\stepsize) - \cost_t(\state_t, \ctrl_t)\right) + \cost_\horizon(\state_\horizon + \diffstate_\horizon^\stepsize) - \cost_{\horizon}(\state_\horizon)
\leq {\stepsize} \costogo_0(0),
$
with $\costogo_0(0)$ the solution of the linear quadratic control
problem~\eqref{eq:LQR_gauss_newton}.
\end{enumerate} 
The procedure is then repeated on the next sequence of control variables.
Ignoring the line-search phase (namely, taking $\stepsize=1$), each iteration
can be summarized as computing $\ctrls^{\nxt} = \ctrls + \diffctrls$ where 
\[
\diffctrls = \DynProg((\lin_{\dyn_t}^{\state_t, \ctrl_t})_{t=0}^{\horizon-1}, (\qua_{\cost_t}^{\state_t, \ctrl_t})_{t=0}^{\horizon} , \diffstate_0, \BellLQ)
\]
for $\diffstate_0 =0$, where $\DynProg$ is the dynamic programming procedure implemented in Algo.~\ref{algo:dyn_prog}. Note that for convex costs $\cost_t$ such
that $\cost_t(\state, \cdot)$ is strongly convex,  the
subproblems~\eqref{eq:LQR_gauss_newton} satisfy the assumptions of
Cor.~\ref{corr:lin_quad}.

The iterations of the following nonlinear control algorithms can always be
decomposed into the three passes described above for the ILQR algorithm. The
algorithms  vary by (i) what approximations of the dynamics and the costs are
computed in the forward pass, (ii) how the policies are computed in the backward
pass, (iii) how the policies  are rolled out.

	\section{Classical Optimization Oracle}\label{sec:classical_optim}
	Problem~\eqref{eq:discrete_ctrl_pb} is entirely determined by the choice of the
initial state and a sequence of control variables, such that the objective
in~\eqref{eq:discrete_ctrl_pb} can be written   in terms of the control
variables  $\ctrls = (\ctrl_0;\ldots;\ctrl_{\horizon-1})$ as
\begin{align*}
 \obj(\ctrls) \coloneqq &  \sum_{t=0}^{\horizon-1} \cost_t(\state_t, \ctrl_t) +\cost_\horizon(\state_\horizon) \\
& \mbox{s.t.} \quad \state_{t+1} = f_t(\state_t, \ctrl_t) \quad \mbox{for} \ t \in \{0, \ldots, \horizon-1\}, \quad \state_0 = \initstate.
\end{align*}
The objective can be decomposed into the costs and  the control of $\horizon$
steps of a sequence of dynamics defined as follows.
\begin{defi}\label{def:traj_func} We define the control of $\horizon$ discrete
	time dynamics $(\dyn_t:\reals^\dimstate \times \reals^\dimctrl \rightarrow
	\reals^\dimstate)_{t=0}^{\horizon-1}$ as the function $\traj: \reals^\dimstate
	\times \reals^{\horizon\dimctrl} \rightarrow \reals^{\horizon \dimstate}$,
	which, given an initial point $\state_0 \in \reals^{\dimstate}$ and a sequence
	of controls $\ctrls = (\ctrl_0;\ldots;\ctrl_{\horizon-1})\in
	\reals^{\horizon\dimctrl}$, outputs the corresponding trajectory $\state_1,
	\ldots, \state_\horizon$, i.e., 
	\begin{align}
		\label{eq:traj}
		\traj(\state_0, \ctrls) & \coloneqq ( \state_1;\ldots;\state_\horizon) \\
		\mbox{s.t.} \quad \state_{t+1} &= \dyn_t(\state_t, \ctrl_t) \quad \mbox{for} \ t \in \{0,\ldots, \horizon-1\}. \nonumber
	\end{align}
\end{defi}
Overall, problem~\eqref{eq:discrete_ctrl_pb} can  be written as the minimization
of a composition
\begin{align}\label{eq:composite_pb}
\min_{\ctrls \in \reals^{\horizon \dimctrl}} \ & \left\{\obj(\ctrls)= \cost\circ\augtraj(\ctrls)\right\}, \quad 
\mbox{where} \quad  \cost(\states, \ctrls) = \sum_{t=0}^{\horizon-1} \cost_t(\state_t, \ctrl_t) +\cost_\horizon(\state_\horizon), \quad \augtraj(\ctrls) = (\traj(\initstate, \ctrls), \ctrls),
\end{align}
for $\states = (\state_1;\ldots;\state_\horizon)$ and $\ctrls =
(\ctrl_0;\ldots;\ctrl_{\horizon-1})$.
The implementation of classical oracles for problem~\eqref{eq:composite_pb}
relies on the dynamical structure of the problem encapsulated in the control
$\traj$ of the discrete time  dynamics $(\dyn_t)_{t=0}^{\horizon-1}$.

\subsection{Formulation}\label{ssec:optim_oracles_formulation} Classical
optimization algorithms rely on the availability of oracles for the
objective. Here, we consider these oracles to compute the minimizer of an
approximation of the objective around the current point with an optional
regularization term. Formally, at a point $\ctrls\in \reals^\dimvar$, given a
regularization $\reg \geq 0$, for an objective of the form 
\[
\min_{\ctrls \in \reals^{\horizon\dimctrl}}\quad \cost\circ \augtrajfunc(\ctrls),
\]
as in~\eqref{eq:composite_pb}, we consider
\begin{enumerate}[nosep]
	\item[(i)] a \emph{gradient} oracle to use a linear expansion of the
	objective, and to output, for $\reg>0$,
	\begin{align}
 \argmin_{\diffctrls\in \reals^\dimvar} \: &
	\left\{	\lin_{\cost \circ \augtrajfunc}^{\ctrls}(\diffctrls) +\frac{\reg}{2}\|\diffctrls\|_2^2\right\}
	= -\reg^{-1} \nabla (\cost\circ \augtraj)(\ctrls), \label{eq:grad_step}
	\end{align}
	\item[(ii)] a  \emph{Gauss-Newton} oracle to use a linear quadratic expansion
	of the objective, and to output
	\begin{align}
	 \argmin_{\diffctrls \in \reals^\dimvar}\: & 
		  \left\{\qua_{\cost}^{\augtrajfunc(\ctrls)}(\lin_\augtrajfunc^\ctrls(\diffctrls))  + \frac{\reg}{2}\|\diffctrls\|_2^2\right\}
		  = -(\nabla \augtraj(\ctrls)\nabla^2\cost(\augtraj(\ctrls)) \nabla \augtraj(\ctrls)+ \reg\idm)^{-1} \nabla (\cost\circ \augtraj)(\ctrls), \label{eq:gn_step}
	\end{align}
	\item[(iii)] a  \emph{Newton} oracle to use a quadratic expansion of the
	objective, and to output
\end{enumerate}
\begin{align}
	 \argmin_{\diffctrls \in \reals^\dimvar} \: & 
	\left\{\qua_{\cost\circ \augtrajfunc}^\ctrls(\diffctrls)  + \frac{\reg}{2}\|\diffctrls \|_2^2 \right\}
	= -(\nabla^2 (\cost\circ \augtraj)(\ctrls) + \reg \idm)^{-1}\nabla (\cost\circ \augtraj)(\ctrls), \label{eq:newton_step}
\end{align}
where $\lin_f^x$, $\qua_f^x$ are the linear and quadratic expansions of $f$ around $x$ as defined in the notations in
Eq.~\eqref{eq:lin_approx}.

Gauss-Newton and Newton oracles are generally defined without a regularization,
i.e., for $\reg=0$. However, in practice, a regularization may be necessary to
ensure that Gauss-Newton and Newton oracles provide a descent direction.
Moreover, the reciprocal of the regularization, $1/\reg$, can play the role of a
stepsize as detailed in Appendix~\ref{app:line_searches}.
\revised{The regularization $\reg$ can then vary with the 
iterates similarly as in trust region methods~\citep[Chapter 4]{nocedal2006numerical}.}
Lemma~\ref{lem:lin_quad_oracle} presents how the computation of the above
oracles can be decomposed into the dynamical structure of the problem. The proof
is detailed in Appendix~\ref{app:proofs}.

\begin{lemm}\label{lem:lin_quad_oracle} Consider a nonlinear dynamical problem
	summarized as 
	\[
	\min_{\ctrls\in \reals^{\horizon\dimctrl}} \cost\circ\augtraj(\ctrls), \quad 
		\mbox{where} \quad  \cost(\states, \ctrls) = \sum_{t=0}^{\horizon-1} \cost_t(\state_t, \ctrl_t) +\cost_\horizon(\state_\horizon), \quad \augtraj(\ctrls) = (\traj(\initstate, \ctrls), \ctrls),
	\]
	with $\traj$ the control of $\horizon$  dynamics $(\dyn_t)_{t=0}^{\horizon-1}$
	as defined in Def.~\ref{def:traj_func}.

	Let $\ctrls = (\ctrl_0; \ldots;\ctrl_{\horizon-1})$ and $\trajfunc(\initstate,
	\ctrls) = (\state_1;\ldots;\state_\horizon)$. 	Gradient~\eqref{eq:grad_step},
	Gauss-Newton~\eqref{eq:gn_step} and Newton~\eqref{eq:newton_step} oracles for
	$\cost\circ \augtrajfunc$ amount to solving for  $\diffctrls^*=(\diffctrl_0^*;
	\ldots;\diffctrl_{\horizon-1}^* )$ linear quadratic control problems of the
	form
	\begin{align}
		\min_{\substack{\diffctrl_0,\ldots, \diffctrl_{\horizon-1} \in \reals^{\dimctrl}\\\diffstate_0,\ldots, \diffstate_\horizon \in \reals^{\dimstate}}} \quad &  
		\sum_{t=0}^{\horizon-1} \qua_t(\diffstate_t, \diffctrl_t) + \qua_\horizon(\diffstate_\horizon)\label{eq:lin_quad_oracle}\\
		\mbox{subject to} \quad & \diffstate_{t+1}=  \lin_{\dyn_t}^{\state_t, \ctrl_t} (\diffstate_t, \diffctrl_t)\quad \mbox{for} \  t \in \{0,\ldots,\horizon-1\}, \quad  \diffstate_0 = 0, \nonumber
	\end{align}
	where for
	\begin{enumerate}[nosep]
		\item[(i)]  the gradient oracle~\eqref{eq:grad_step},
		$\qua_\horizon(\diffstate_\horizon) =
		\lin_{\cost_\horizon}^{\state_\horizon}(\auxstate_\horizon) $ and,  for $0
		\leq t \leq  \horizon-1$, 
		\[
	\qua_t(\diffstate_t, \diffctrl_t) = \lin_{\cost_t}^{\state_t, \ctrl_t}(\diffstate_t, \diffctrl_t) + \frac{\reg}{2} \|\diffctrl_t\|_2^2,
		\] 
		\item[(ii)] the Gauss-Newton oracle~\eqref{eq:gn_step},
		$\qua_\horizon(\diffstate_\horizon) =
		\qua_{\cost_\horizon}^{\state_\horizon}(\auxstate_\horizon) $ and,  for $0
		\leq t\leq \horizon-1$, 
			\[
		\qua_t(\diffstate_t, \diffctrl_t)	= \qua_{\cost_t}^{\state_t, \ctrl_t}(\diffstate_t, \diffctrl_t) + \frac{\reg}{2} \|\diffctrl_t\|_2^2,
\]

		\item[(iii)] for the  Newton oracle~\eqref{eq:newton_step},
		$\qua_\horizon(\diffstate_\horizon) =
		\qua_{\cost_\horizon}^{\state_\horizon}(\auxstate_\horizon) $ and, defining 
		\begin{gather}\label{eq:adjoint_dyn}
			\costate_\horizon  =\nabla\cost_\horizon(\state_\horizon), \quad  
			\costate_{t}  = \nabla_{\state_t} \cost_t(\state_t, \ctrl_t) +   \nabla_{\state_t} \dyn_t(\state_t, \ctrl_t) \costate_{t+1} \quad \mbox{for} \ t\in\{\horizon-1,\ldots,1\},
		\end{gather}
		we have, for $0\leq t\leq \horizon-1$,
\[
			\qua_t(\diffstate_t, \diffctrl_t) = \qua_{\cost_t}^{\state_t, \ctrl_t} (\diffstate_t, \diffctrl_t) + \frac{1}{2} \nabla^2 \dyn_t(\state_t, \ctrl_t)[\cdot, \cdot, \lambda_{t+1}](\diffstate_t, \diffctrl_t)+ \frac{\reg}{2} \|\diffctrl_t\|_2^2,
\] 
where for $f:\reals^\dimstate \times \reals^\dimctrl \rightarrow
\reals^\dimstate$, $\state\in \reals^\dimstate$, $\ctrl \in \reals^\dimctrl$,
$\costate\in \reals^\dimstate$, we define
\begin{align}\label{eq:additional_quad_newton}
\nabla^2\dyn(\state, \ctrl)[\cdot, \cdot,  \costate]: (\diffstate, \diffctrl) \rightarrow & \nabla_{\state \state}^2\dyn(\state,\ctrl )[\diffstate, \diffstate, \costate ] + 2\nabla_{\state \ctrl}^2\dyn(\state,\ctrl )[\diffstate, \diffctrl, \costate]  {+} \nabla_{\ctrl \ctrl}^2\dyn(\state,\ctrl )[\diffctrl, \diffctrl, \costate ].
\end{align}
	\end{enumerate}
\end{lemm} 

From an optimization viewpoint, gradient, Gauss-Newton or Newton oracles are
considered as black-boxes. Second order methods such as Gauss-Newton or Newton
methods generally require solving a linear system at a cubic cost in the
dimension of the problem~\citep[Chapter 4]{nesterov2018lectures}. Here, the
dimension of the problem in the control variables is $\horizon \dimctrl$, with
$\dimctrl$, the dimension of the control variables, usually small (see the
numerical examples in Sec.~\ref{sec:exp}), but $\horizon$, the number of time
steps, potentially large if, e.g., the discretization time step used to
define~\eqref{eq:discrete_ctrl_pb} from a  continuous time control problem is
small while the original time length of the continuous time control problem is
large. A cubic cost w.r.t. the number of time steps $\horizon$ is then a priori
prohibitive. 

A closer look at the implementation of all the above
oracles~\eqref{eq:grad_step},~\eqref{eq:gn_step},~\eqref{eq:newton_step}, shows
that they all amount to solving linear quadratic control problems as presented
in Lemma~\ref{lem:lin_quad_oracle}. Hence, they can be solved by a dynamic
programming approach detailed in Sec.~\ref{ssec:optim_oracles_implementation} at
a cost linear w.r.t. the number of time steps $\horizon$. As a consequence, if
the dimensions $\dimctrl, \dimstate$ of the control and state variables are
negligible compared to the horizon $\horizon$, the computational complexities of
Gauss-Newton and Newton oracles, detailed in Sec.~\ref{sec:comput_cplxity} are
of the same order as the computational complexity of a gradient oracle. This
observation was done by~\citet{de1988differential, dunn1989efficient} for a
Newton step and~\citet{sideris2005efficient} for a Gauss-Newton step.
\citet{wright1990solution} also presented how sequential quadratic programming
methods can naturally be cast in a similar way. Lemma~\ref{lem:lin_quad_oracle}
casts all classical optimization oracles in the same formulation, including a
gradient oracle. 

The linear quadratic control problems can be solved by different procedures than
dynamic programming such as using Riccati-based or parallel implementations as
detailed in Appendix~\ref{app:sparse_solvers}~\citep{wright1991partitioned}. We
focus on their resolution by dynamic programming to cast all algorithms in a
common framework.

\subsection{Implementation}\label{ssec:optim_oracles_implementation}

Given Lemma~\ref{lem:lin_quad_oracle}, for $\traj(\initstate, \ctrls)$ the control of $\horizon$ dynamics
$(\dyn_t)_{t=0}^{\horizon-1}$ defined in Def.~\ref{def:traj_func},  classical optimization oracles for
objectives of the form
\begin{align}\nonumber
\obj(\ctrls)= \cost\circ\augtraj(\ctrls), \quad 
	\mbox{where} \quad  \cost(\states, \ctrls) = \sum_{t=0}^{\horizon-1} \cost_t(\state_t, \ctrl_t) +\cost_\horizon(\state_\horizon), \quad \augtraj(\ctrls) = (\traj(\initstate, \ctrls), \ctrls),
\end{align}
can be implemented by (i) instantiating the linear quadratic control
problem~\eqref{eq:lin_quad_oracle} with the chosen approximations, (ii) solving
the linear quadratic control problem~\eqref{eq:lin_quad_oracle} by dynamic
programming as detailed in Sec.~\ref{sec:lin_quad}. Precisely, their
implementation can be split into the following three phases.

\begin{enumerate}
\item \underline{Forward pass:}  All oracles start by gathering the information
necessary for the step in a forward pass that takes the generic form of
Algo.~\ref{algo:forward} and can be summarized as 
\[
\obj(\ctrls), (\model_{\dyn_t}^{\state_t, \ctrl_t})_{t=0}^{\horizon-1}, (\model_{\cost_t}^{\state_t, \ctrl_t})_{t=0}^{\horizon-1}, \model_{\cost_\horizon}^{\state_\horizon}  = \Forward(\ctrls,(\dyn_t)_{t=0}^{\horizon-1}, (\cost_t)_{t=0}^\horizon, \initstate, o_\dyn, o_\cost )  
\]
that compute the objective $\obj(\ctrls)$ associated to the given sequence of
controls $\ctrls$ and record  approximations $(\model_{\dyn_t}^{\state_t,
\ctrl_t})_{t=0}^{\horizon-1}, (\model_{\cost_t}^{\state_t,
\ctrl_t})_{t=0}^{\horizon-1}, \model_{\cost_\horizon}^{\state_\horizon} $ of the
dynamics and the costs up to the orders $o_\dyn$ and $o_\cost$, 
\revised{respectively as
\begin{equation}\label{eq:model_approxs}
  \model_{\dyn_t}^{\state_t, \ctrl_t} = \begin{cases} 
    \lin_{\dyn_t}^{\state_t, \ctrl_t} &\mbox{if} \ o_\dyn = 1 \\
    \qua_{\dyn_t}^{\state_t, \ctrl_t} & \mbox{if} \ o_\dyn = 2
  \end{cases}, \quad  	\model_{\cost_t}^{\state_t, \ctrl_t} = \begin{cases}
    \lin_{\cost_t}^{\state_t, \ctrl_t} &\mbox{if} \ o_\cost = 1 \\
    \qua_{\cost_t}^{\state_t, \ctrl_t} & \mbox{if} \ o_\cost = 2
  \end{cases}, \quad \model_{\cost_\horizon}^{\state_\horizon} = \begin{cases}
    \lin_{\cost_\horizon}^{\state_\horizon} &\mbox{if} \ o_\cost = 1 \\
    \qua_{\cost_\horizon}^{\state_\horizon} & \mbox{if} \ o_\cost = 2.
  \end{cases}
\end{equation}
}
The orders of approximation $o_f, o_h$ for each algorithm are summarized 
in Fig.~\ref{fig:taxonomy}.
\item \underline{Backward pass:} Once approximations of the dynamics have been
computed, a backward pass on the corresponding linear quadratic control
problem~\eqref{eq:lin_quad_oracle} can be done as in the linear quadratic case
presented in Sec.~\ref{sec:lin_quad}. The backward passes of the gradient oracle
in Algo.~\ref{algo:backward_grad}, the Gauss-Newton oracle in
Algo.~\ref{algo:backward_gn} and the Newton oracle in
Algo.~\ref{algo:backward_newton} take generally the form 
\[
(\pi_t)_{t=0}^{\horizon-1}, \costogo_0 = \Backward((\model_{\dyn_t}^{\state_t, \ctrl_t})_{t=0}^{\horizon-1}, (\model_{\cost_t}^{\state_t, \ctrl_t})_{t=0}^{\horizon-1}, \model_{\cost_\horizon}^{\state_\horizon}, \reg).
\]
Namely, they take as input  a regularization $\reg\geq 0$ and some
approximations of the dynamics and the costs $(\model_{\dyn_t}^{\state_t,
\ctrl_t})_{t=0}^{\horizon-1}, (\model_{\cost_t}^{\state_t,
\ctrl_t})_{t=0}^{\horizon-1}, \model_{\cost_\horizon}^{\state_\horizon}$
computed in a forward pass,  and return a set of policies and the final
cost-to-go corresponding to the subproblem~\eqref{eq:lin_quad_oracle}. 
\item \underline{Roll-out pass:} Given the output of a backward pass defined
above, the oracle is computed by rolling out the policies along the linear
trajectories defined in the subproblem~\eqref{eq:lin_quad_oracle}. Formally,
given a sequence of policies $(\pi_t)_{t=0}^{\horizon-1}$, the oracles are then
given as $\diffctrls = (\diffctrl_0;\ldots;\diffctrl_{\horizon-1})$ computed,
for $\diffstate_0 =0$, by Algo.~\ref{algo:roll_out} as
\[
\diffctrls = \Roll(\diffstate_0,( \pi_t)_{t=0}^{\horizon-1},  (\lin_{\dyn_t} ^{\state_t, \ctrl_t})_{t=0}^{\horizon-1}).
\]
Here the policies $(\pi_t)_{t=0}^{\horizon-1}$ are output by one of the backward
passes in Algo.~\ref{algo:backward_grad}, Algo.~\ref{algo:backward_gn} or
Algo.~\ref{algo:backward_newton}. For the Gauss-Newton and Newton oracles, an
additional procedure checks whether the subproblems are convex at each iteration
as explained in more detail in Appendix~\ref{app:summary}.
\end{enumerate} 
Gradient, Gauss-Newton, and Newton oracles are implemented by, respectively,
Algo.~\ref{algo:gd}, Algo.~\ref{algo:gn}, Algo.~\ref{algo:newton}. Additional
line-searches are presented in Appendix~\ref{app:line_searches}. The
computational schemes of a gradient, a Gauss-Newton and a Newton oracle are illustrated in Fig.~\ref{fig:grad}, Fig.~\ref{fig:ilqr} and Fig.~\ref{fig:newton} respectively.\\

{\bf Gradient back-propagation.}~\\
For a gradient  oracle~\eqref{eq:grad_step}, the procedure $\BellLQ$ normally
used to solve linear quadratic control problems  simplifies to a linear back-propagation,
$\BellL$, presented in Algo.~\ref{algo:bp} that implements
\begin{equation}\label{eq:bellman_l}
	\BellL	: (\lin_t^\dyn, \lin_t^\cost, \costogo_{t+1}, \reg)  \rightarrow \left(\begin{array}{c}
	\costogo_t:\state\rightarrow\min_{\ctrl \in \reals^\dimctrl} \left\{ \lin_t^\cost(\state, \ctrl) + \costogo_{t+1}(\lin_t^\dyn(\state, \ctrl)) + \frac{\reg}{2} \|\ctrl\|_2^2\right\} \\
	\pi_t :\state\rightarrow \argmin_{\ctrl \in \reals^\dimctrl}  \left\{ \lin_t^\cost(\state, \ctrl) + \costogo_{t+1}(\lin_t^\dyn(\state, \ctrl))+ \frac{\reg}{2} \|\ctrl\|_2^2\right\}
\end{array}\right), 
\end{equation}
for linear functions $\lin_t^\dyn, \lin_t^\cost, \costogo_{ t+1}$. Plugging into
the overall dynamic programming procedure, Algo.~\ref{algo:bp}, the
linearizations of the dynamics and the costs, we get that the gradient oracle,
Algo.~\ref{algo:backward_grad}, computes affine cost-to-go functions of the form
$	\costogo_t(\diffstate_t)  = \j_t^\top \diffstate_t + \jcst_t $ with
\begin{gather*}
	\j_\horizon  =\nabla\cost_\horizon(\state_\horizon), \quad  
	\j_{t}  = \nabla_{\state_t} \cost_t(\state_t, \ctrl_t) +   \nabla_{\state_t} \dyn_t(\state_t, \ctrl_t) \j_{t+1} \quad \mbox{for} \ t\in\{0,\ldots, \horizon-1\}.
\end{gather*}
Moreover, the policies are independent of the state variables, i.e.,
$\pi_t(\diffstate_t) = k_t$, with
\[
k_t = - \reg^{-1}(\nabla_{\ctrl_t} \cost_t(\state_t, \ctrl_t) + \nabla_{\ctrl_t}\dyn_t(\state_{t}, \ctrl_t) \j_{t+1}) = - \reg^{-1} \nabla_{\ctrl_t} (\cost\circ \augtraj)(\ctrls).
\]
The roll-out of these policies is independent of the dynamics and output
directly the gradient up to a factor $-\reg^{-1}$. Note that we naturally
retrieve the gradient back-propagation algorithm~\citep{griewank2008evaluating}.

	\section{Differential Dynamic Programming Oracle}\label{sec:ddp}
	The original differential dynamic programming algorithm was developed
by~\citet{jacobson1970differential} and revisited by,
e.g.,~\citet{mayne1975first, murray1984differential, liao1992advantages,
tassa2014control}. The reader can verify from the aforementioned citations that
our presentation matches the original formulation in, e.g., the quadratic case,
while offering a larger perspective on the method that incorporates, e.g.,
linear quadratic approximations. Such approaches have also been called \emph{direct
multiple shooting} by~\citet{bock1984multiple}.

\subsection{Rationale}
Denoting $\cost$ the total cost as in~\eqref{eq:composite_pb} and $\traj$ the
control in $\horizon$ dynamics $(\dyn_t)_{t=0}^{\horizon-1}$, Differential
Dynamic Programming (DDP) oracles consist in solving approximately
\begin{equation}\nonumber
	\min_{\diffctrls \in \reals^{\horizon \dimctrl}} \cost(\traj(\initstate, \ctrls + \diffctrls), \ctrls+\diffctrls),
\end{equation}
by means of a dynamic programming procedure and using the resulting policies to
update the current sequence of controllers. For a consistent presentation with
the classical optimization oracles presented in Sec.~\ref{sec:classical_optim},
we consider a regularized formulation of the DDP oracles, that is, 
\begin{equation}\label{eq:ddp_obj3}
	\min_{\diffctrls \in \reals^{\horizon \dimctrl}} \cost(\traj(\initstate, \ctrls + \diffctrls), \ctrls+\diffctrls) + \frac{\reg}{2} \|\diffctrls\|_2^2,
\end{equation}
for some regularization $\reg \geq 0$. 

The objective in problem~\eqref{eq:ddp_obj3} can be rewritten as 
\begin{equation}\label{eq:decomp_diff_obj}
\cost(\traj(\initstate, \ctrls + \diffctrls), \ctrls+\diffctrls) = \cost(\traj(\initstate, \ctrls)) +  \diff_\cost^{\traj(\initstate, \ctrls)}(\diff_\traj^{\initstate, \ctrls}(0,\diffctrls), \diffctrls),
\end{equation}
where for a function $f$, $\delta_f^{x}$ is the finite difference expression of
$f$ around $x$ as defined in the notations in Eq.~\eqref{eq:lin_approx}. In
particular, $\diff_\traj^{\initstate, \ctrls}(0,\diffctrls)$ is the trajectory
defined by the finite differences of the dynamics given as
\[
\diff_{\dyn_t}^{\state_t, \ctrl_t}(\diffstate_t, \diffctrl_t) = \dyn_t(\state_t + \diffstate_t, \ctrl_t + \diffctrl_t) - \dyn_t(\state_t, \ctrl_t).
\]
The dynamic programming approach is then applied on the above dynamics. Namely,
the goal is to solve
\begin{align}\label{eq:ddp_obj_detailed}
	\min_{\substack{\diffctrl_0,\ldots, \diffctrl_{\horizon-1} \in \reals^{\dimctrl}\\\diffstate_0,\ldots, \diffstate_\horizon \in \reals^{\dimstate}}}  & \sum_{t=0}^{\horizon-1} \diff_{\cost_t}^{\state_t, \ctrl_t}(\diffstate_t, \diffctrl_t)  + \frac{\reg}{2}\|\diffctrl_t\|_2^2 + \diff_{\cost_\horizon}^{\state_\horizon} (\diffstate_\horizon) \\
	\mbox{subject to} \quad & \diffstate_{t+1} = \diff_{\dyn_t}^{\state_t, \ctrl_t} (\diffstate_t, \diffctrl_t) \quad \mbox{for} \ t\in \{0, \ldots,\horizon-1\}, \quad \diffstate_0 = 0 \nonumber,
\end{align}
by dynamic programming. Denote then $\costogo_t^*$ the cost-to-go functions
associated to problem~\eqref{eq:ddp_obj_detailed} for $t\in \{0, \ldots
\horizon\}$. These cost-to-go functions satisfy the recursive equation
\begin{equation}\label{eq:bellman_ddp}
	\costogo_t^*(\diffstate_t) = \min_{\diffctrl_t \in \reals^\dimctrl}  \left\{\diff_{\cost_t}^{\state_t, \ctrl_t}(\diffstate_t, \diffctrl_t)  + \frac{\reg}{2} \|\diffctrl_t\|_2^2+ \costogo_{t+1}^*( \diff_{\dyn_t}^{\state_t, \ctrl_t} (\diffstate_t, \diffctrl_t) )\right\},
\end{equation}
starting from $\costogo_\horizon^* = \diff_{\cost_\horizon}^{\state_\horizon}$
and such that our objective is to compute $\costogo_0^*(0)$. Since the dynamics
$\diff_{\dyn_t}^{\state_t, \ctrl_t}$ are not linear and the costs
$\diff_{\cost_t}^{\state_t, \ctrl_t}$ are not quadratic, there is no analytical
solution for the subproblem~\eqref{eq:bellman_ddp}. To circumvent this issue,
the cost-to-go functions are approximated as 
$
\costogo_t^*(\diffstate_t) \approx \costogo_t(\diffstate_t),
$
where $  \costogo_t$ is computed from approximations of the dynamics and the
costs. The approximation is done around the nominal value of the
subproblem~\eqref{eq:ddp_obj_detailed} which is $\auxctrls=0$ and corresponds to
$\auxstates=0$ and no change of the original objective
in~\eqref{eq:decomp_diff_obj}.

Denoting  $\model_f$  an expansion of a function $f$ around the origin such that
$f(x) \approx f(0) + \model_f(x)$, the cost-to-go functions are computed  with
an approximate back-propagation $\Bellapprox$ of cost-to-go functions:
\begin{equation}\label{eq:bell_approx}
	\Bellapprox:		
		\diff_t^\dyn, \diff_t^\cost, \costogo_{t+1} \rightarrow  \left(\begin{array}{c}
			\costogo_t: \diffstate  \rightarrow  (\diff_t^\cost {+} \costogo_{t+1}{\circ} \diff_t^\dyn)(0, 0) + \min_{\diffctrl \in \reals^\dimctrl}  
			\left\{
			+ \model_{\diff_t^\cost}(\diffstate, \diffctrl)  
			+ \model_{\costogo_{t+1}\circ \diff_t^\dyn}(\diffstate, \diffctrl) 
			+ \frac{\reg}{2}\|\diffctrl\|_2^2 \right\},\\
			\pi_t:  \diffstate  \rightarrow \argmin_{\diffctrl \in \reals^\dimctrl}  
			\left\{\model_{\diff_t^\cost}(\diffstate, \diffctrl) 
			+  \model_{ \costogo_{t+1}\circ \diff_t^\dyn}(\diffstate, \diffctrl)
			+ \frac{\reg}{2}\|\diffctrl\|_2^2  \right\} 
		\end{array}\right),
\end{equation}
applied to the finite differences $\diff_{\dyn_t}^{\state_t, \ctrl_t}
\rightarrow\diff_t^\dyn$ and $\diff_{\cost_t}^{\state_t, \ctrl_t} \rightarrow
\diff_t^\cost$. A DDP oracle computes then a sequence of policies by iterating
in a backward pass, starting from $\costogo_\horizon=
\model_{\diff_{\cost_\horizon}^{\state_\horizon}}$,
\begin{equation}\label{eq:bell_approx_rec}
	\costogo_t, \pi_t = \Bellapprox(\diff_{\dyn_t}^{\state_t, \ctrl_t},\diff_{\cost_t}^{\state_t, \ctrl_t}, \costogo_{t+1} ) \quad \mbox{for} \ t\in \{\horizon-1, \ldots, 0\}.
\end{equation}
Given a set of policies, an approximate solution is given by rolling out the
policies along the dynamics defining problem~\eqref{eq:ddp_obj_detailed}, i.e.,
by computing $\diffctrl_0, \ldots, \diffctrl_{\horizon-1}$ as
\begin{equation}\label{eq:rollout_ddp0}
	\diffctrl_t= \pi_t( \diffstate_t), \qquad \diffstate_{t+1} =\diff_{\dyn_t}^{\state_t, \ctrl_t} (\diffstate_t, \diffctrl_t)=  \dyn_t(\state_t + \diffstate_t, \ctrl_t + \diffctrl_t ) - \dyn_t(\state_t, \ctrl_t) \qquad \mbox{for $t=0, \ldots, \horizon-1$}.
\end{equation}
The main difference with the classical optimization oracles lies a priori in the
computation of the policies in~\eqref{eq:bell_approx_rec} detailed below and in
the roll-out pass that uses the finite differences of the dynamics. The constant part of the cost-to-go functions is used for line-searches as detailed in
Appendix~\ref{app:line_searches}.

\subsection{Detailed Derivation of the Backward Passes}
~

{\bf Linear Approximation.} 
If we consider a linear approximation for the composition of the cost-to-go
function and the dynamics, we have
\[
\model_{\costogo_{t+1}\circ \diff_{\dyn_t}^{\state, \ctrl}} 
=\lin_{\costogo_{t+1}\circ \diff_{\dyn}^{\state, \ctrl}} 
=  \lin_{\costogo_{t+1}}^{\diff_{\dyn}^{\state, \ctrl}(0, 0)}\circ \lin_{\diff_{\dyn}^{\state, \ctrl}} 
= \lin_{\costogo_{ t+1}}\circ \lin_{\dyn}^{\state, \ctrl},
\]
where we denote simply $\lin_{f} = \lin_f^{0}$ the linear expansion of a
function $f$ around the origin.  

Plugging this model into~\eqref{eq:bell_approx} and using linear approximations
of the costs, the recursion~\eqref{eq:bell_approx_rec} amounts to computing,
starting from $\costogo_\horizon =
\lin_{\diff_{\cost_\horizon}^{\state_\horizon, \ctrl_\horizon}} =
\lin_{\cost_\horizon}^{\state_\horizon, \ctrl_\horizon}$,
\begin{align}\label{eq:rec_bell_approx_lin}\nonumber
\costogo_t(\diffstate)	&	= \diff_{\cost_t}^{\state_t, \ctrl_t}(0, 0) + \min_{\diffctrl\in \reals^{\dimctrl}}\lin_{\diff_{\cost_t}^{\state_t, \ctrl_t}}(\diffstate, \diffctrl) + \costogo_{t+1}(\diff_{\dyn_t}^{\state_t, \ctrl_t}(0, 0)) +  \lin_{\costogo_{t+1}}(\lin_{\dyn_t}^{\state_t,\ctrl_t}(\diffstate, \diffctrl)) + \frac{\reg}{2} \|\diffctrl\|_2^2,\\
&	=  \min_{\diffctrl\in \reals^{\dimctrl}}\lin_{\cost_t}^{\state_t, \ctrl_t}(\diffstate, \diffctrl) +  \costogo_{t+1}(\lin_{\dyn_t}^{\state_t,\ctrl_t}(\diffstate, \diffctrl)) + \frac{\reg}{2} \|\diffctrl\|_2^2,\nonumber
\end{align}
 where in the last line we used that the cost-to-go functions $\costogo_t$  are
 necessarily affine, s.t. $\costogo_{ t+1}(\diffstate) = \costogo_{t+1}(0) +
 \lin_{\costogo_{t+1}}(\diffstate)$.  We retrieve then the same recursion as the
 one used for a gradient oracle~\eqref{eq:bellman_l}, with the same policies. Since
 the computed policies are constant, they are not affected by the dynamics along
 which a roll-out phase is performed. In other words, the oracle returned by
 using linear approximations in a DDP approach is just a gradient oracle.

{\bf Linear Quadratic Approximation.}
If we consider a linear quadratic  approximation for the composition of the
cost-to-go function and the dynamics, we have
\[
\model_{\costogo_{t+1}\circ \diff_{\dyn^{\state, \ctrl}}} 
= \qua_{\costogo_{t+1}}^{\diff_{\dyn}^{\state, \ctrl}(0, 0)}\circ \lin_{\diff_{\dyn}^{\state, \ctrl}}
= \qua_{\costogo_{t+1}} \circ \lin_{\dyn}^{\state, \ctrl}, \]
 where we denote simply $\qua_{f} = \qua_f^{0}$ the quadratic expansion of a
 function $f$ around the origin.  
Plugging this model into~\eqref{eq:bell_approx} and using quadratic
approximations of the costs, the recursion~\eqref{eq:bell_approx_rec} amounts to
computing, starting from $\costogo_\horizon =
\qua_{\diff_{\cost_\horizon}^{\state_\horizon, \ctrl_\horizon}} =
\qua_{\cost_\horizon}^{\state_\horizon, \ctrl_\horizon}$,
\begin{align}
	\costogo_t(\diffstate) & =\diff_{\cost_t}^{\state_t, \ctrl_t}(0, 0)  +\min_{\diffctrl\in \reals^{\dimctrl}} \qua_{\diff_{\cost_t}^{\state_t, \ctrl_t}}(\diffstate, \diffctrl) +   \costogo_{t+1}(\diff_{\dyn_t}^{\state_t, \ctrl_t}(0, 0))  + \qua_{\costogo_{t+1}}^{\diff_{\dyn}^{\state, \ctrl}(0, 0)}\circ \lin_{\diff_{\dyn}^{\state, \ctrl}}^{(0,0)} (\diffstate, \diffctrl) + \frac{\reg}{2} \|\diffctrl\|_2^2\nonumber\\
	& = \min_{\diffctrl\in \reals^{\dimctrl}}  \qua_{\cost_t}^{\state_t, \ctrl_t}(\diffstate, \diffctrl) +   \costogo_{t+1}(0)  +  \qua_{\costogo_{t+1}} ( \lin_{\dyn_t}^{\state_t, \ctrl_t}(\diffstate, \diffctrl)) + \frac{\reg}{2} \|\diffctrl\|_2^2. \label{eq:rec_bell_approx_lq}
\end{align}
If the costs $\cost_t$ are convex for all $t$ and $\qua_{\cost_t}^{\state_t,
\ctrl_t}(\diffstate, \cdot) + \frac{\reg}{2} \|\cdot\|_2^2$ is strongly convex
for all $t$ and all $\diffstate$, then the cost-to-go functions $\costogo_t$ are
convex quadratics for all $t$, i.e., $\costogo_{t+1}(\diffstate) =
\costogo_{t+1}(0) + \qua_{\costogo_{t+1}}(\diffstate)$. In that case, the
recursion~\eqref{eq:rec_bell_approx_lq} simplifies as
\begin{equation}\label{eq:rec_bell_approx_lq_simplified}
	\costogo_t(\diffstate) =\min_{\diffctrl\in \reals^{\dimctrl}}\qua_{\cost_t}^{\state_t, \ctrl_t}(\diffstate, \diffctrl) +  \costogo_{t+1}(\lin_{\dyn_t}^{\state_t,\ctrl_t}(\diffstate, \diffctrl))+ \frac{\reg}{2} \|\diffctrl\|_2^2,
\end{equation} 
and the policies are given by the minimizer of
Eq.~\eqref{eq:rec_bell_approx_lq_simplified}. The
recursion~\eqref{eq:rec_bell_approx_lq_simplified} is then the same as the
recursion done when computing a Gauss-Newton oracle. Namely, the backward pass
in this case is the backward pass of a Gauss-Newton oracle.  Though the output
policies are the same, the output of the oracle will differ since the roll-out
phase does not follow the linearized trajectories in the DDP approach. The
computational scheme of a DDP approach with linear quadratic approximations
presented in Fig.~\ref{fig:ddplq} is then almost the same as the one of a
Gauss-Newton oracle presented in Fig.~\ref{fig:ilqr}, except that in the
roll-out phase the linear approximations of the dynamics are replaced by finite
differences of the dynamics. This DDP approach amounts to the iterative Linear 
Quadratic Regulator (iLQR) developed by~\citet{tassa2012synthesis}.\\

{\bf Quadratic Approximation.} If we consider a quadratic approximation for the composition of the cost-to-go
function and the dynamics, we get
\[
\model_{\costogo_{t+1}\circ \diff_{\dyn}^{\state, \ctrl}} 
= \qua_{\costogo_{t+1}\circ \diff_{\dyn}^{\state, \ctrl}} 
= \frac{1}{2} \nabla^ 2 \dyn(\state, \ctrl)[\cdot, \cdot, \nabla \costogo_{t+1}(0)] 
+ \qua_{\costogo_{t+1}} \circ \lin_{\dyn}^{\state, \ctrl},
\]
where $\nabla ^ 2 \dyn(\state, \ctrl)[\cdot, \cdot, \costate]$ is defined
in~\eqref{eq:additional_quad_newton}. Plugging this model
into~\eqref{eq:bell_approx} and using quadratic approximations of the costs, the
recursion~\eqref{eq:bell_approx_rec} amounts to, starting from
$\costogo_\horizon = \qua_{\diff_{\cost_\horizon}^{\state_\horizon,
\ctrl_\horizon}} = \qua_{\cost_\horizon}^{\state_\horizon, \ctrl_\horizon}$,
\begin{align}\label{eq:rec_bell_approx_quad}
	\costogo_t(\diffstate) & = \diff_{\cost_t}^{\state_t, \ctrl_t}(0, 0)  
	+\min_{\diffctrl\in \reals^{\dimctrl}} 
	 \qua_{\diff_{\cost_t}^{\state_t, \ctrl_t}}(\diffstate, \diffctrl) 
	+   \costogo_{t+1}(\diff_{\dyn_t}^{\state_t, \ctrl_t}(0, 0))  
	+  \qua_{\costogo_{t+1}\circ \diff_{\dyn_t}^{\state, \ctrl}} (\diffstate, \diffctrl) 
	+ \frac{\reg}{2} \|\diffctrl\|_2^2 \\
& 	=\min_{\diffctrl\in \reals^{\dimctrl}} 
	 \qua_{\cost_t}^{\state_t, \ctrl_t}(\diffstate, \diffctrl) 
	+   \costogo_{t+1}(0)  
	+  \qua_{\costogo_{t+1}} \circ \lin_{\dyn_t}^{\state_t, \ctrl_t}(\diffstate, \diffctrl)  
	+ \frac{1}{2} \nabla^ 2 \dyn_t(\state_t, \ctrl_t)[\cdot, \cdot, \nabla \costogo_{t+1}(0)](\diffstate, \diffctrl)
	+ \frac{\reg}{2} \|\diffctrl\|_2^2. \nonumber
\end{align}
Provided that the costs are convex and that $\qua_{\cost_t}^{\state_t,
\ctrl_t}(\diffstate, \cdot) + \frac{1}{2} \nabla^ 2 \dyn_t(\state_t,
\ctrl_t)[\cdot, \cdot, \nabla \costogo_{t+1}(0)](\diffstate, \cdot) +
\frac{\reg}{2} \|\cdot\|_2^2$ is strongly convex for all $t$ and all
$\diffstate$, the cost-to-go functions $\costogo_t$ are convex quadratics for
all $t$. In that case, the recursion~\eqref{eq:rec_bell_approx_quad} simplifies
as
\begin{equation}\label{eq:rec_bell_approx_quad_simplified}
	\costogo_t(\diffstate) =
	\min_{\diffctrl\in \reals^{\dimctrl}}\qua_{\cost_t}^{\state_t, \ctrl_t}(\diffstate, \diffctrl) 
	+  \costogo_{t+1}(\lin_{\dyn_t}^{\state_t,\ctrl_t}(\diffstate, \diffctrl))
	+ \frac{1}{2} \nabla^ 2 \dyn_t(\state_t, \ctrl_t)[\cdot, \cdot, \nabla \costogo_{t+1}(0)](\diffstate, \diffctrl) 
	+ \frac{\reg}{2} \|\diffctrl\|_2^2,
\end{equation} 
and the policies are given by the minimizer of
Eq.~\eqref{eq:rec_bell_approx_quad_simplified}. The overall backward pass is
detailed in Algo.~\ref{algo:backward_ddp}. 

Compared to the backward pass of the Newton oracle in
Algo.~\ref{algo:backward_newton}, we note that the additional cost derived from
the curvatures of the dynamics is not computed the same way. Namely, the Newton
oracle computes this additional cost by using back-propagated adjoint variables
in Eq.~\eqref{eq:adjoint_dyn}, while in the DDP approach the additional cost is
directly defined through the previously computed cost-to-go function.
Fig.~\ref{fig:ddpq} illustrates the computational scheme of the implementation
of DDP with quadratic approximations and can be compared to the computational
scheme of the Newton oracle in Fig.~\ref{fig:newton}.

Note that, while we used second order Taylor expansions for the compositions and
the costs, the approximate cost-to-go-functions $\costogo_t$ are \emph{not}
second order Taylor expansion of the true cost-to-go functions $\costogo_t^*$,
except for $\costogo_\horizon$. Indeed, $\costogo_t$ is computed as an
approximate solution of the Bellman equation. The true Taylor expansion of the
cost-to-go function requires the gradient and the Hessian of the cost and the
dynamic  in Eq.~\eqref{eq:rec_bell_approx_quad} computed at the minimizer of the
subproblem. Here, since we only use an approximation of the minimizer, we do not
have access to the true gradient and Hessian of the cost-to-go function.

\subsection{Implementation}
The implementation of the DDP oracles follows the same steps as the ones given
for classical optimization oracles as detailed below. The implementation of a
DDP oracle with linear quadratic approximations is given in
Algo.~\ref{algo:ddp_lq} and illustrated in Fig.~\ref{fig:ddplq}. The
implementation of a DDP oracle with quadratic approximations is given in
Algo.~\ref{algo:ddp_q} and illustrated in Fig.~\ref{fig:ddpq}.

\begin{enumerate}
	\item \underline{Forward pass:} As for the classical optimization methods, the
	oracles start by gathering the
	information necessary for the backward pass using
	Algo.~\ref{algo:forward} that computes
	\[
\obj(\ctrls), (\model_{\dyn_t}^{\state_t, \ctrl_t})_{t=0}^{\horizon-1}, (\model_{\cost_t}^{\state_t, \ctrl_t})_{t=0}^{\horizon-1}, \model_{\cost_\horizon}^{\state_\horizon}  = \Forward(\ctrls,(\dyn_t)_{t=0}^{\horizon-1}, (\cost_t)_{t=0}^\horizon, \initstate, o_\dyn, o_\cost ),  
	\] 
	where $o_\dyn$ and $o_\cost$ define the order of the approximations
  $(\model_{\dyn_t}^{\state_t,
  \ctrl_t})_{t=0}^{\horizon-1}, (\model_{\cost_t}^{\state_t,
  \ctrl_t})_{t=0}^{\horizon-1}, \model_{\cost_\horizon}^{\state_\horizon} $ of the
  dynamics and the costs up to the orders $o_\dyn$ and $o_\cost$ as 
  in~\eqref{eq:model_approxs}.
	\item \underline{Backward pass:} As for the classical optimization oracles,
	the backward pass can generally be written
	\[
	(\pi_t)_{t=0}^{\horizon-1}, \costogo_0 = \Backward((\model_{\dyn_t}^{\state_t, \ctrl_t})_{t=0}^{\horizon-1}, (\model_{\cost_t}^{\state_t, \ctrl_t})_{t=0}^{\horizon-1}, \model_{\cost_\horizon}^{\state_\horizon}, \reg),
	\]
	 If linear approximations are used, the backward pass is given in
	 Algo.~\ref{algo:backward_grad}, if linear quadratic approximations are used,
	 the backward pass is given in Algo.~\ref{algo:backward_gn} and if quadratic
	 approximations are used, the backward pass is given in
	 Algo.~\ref{algo:backward_ddp}.
	\item \underline{Roll-out pass:} The roll-out phase differs by using finite
	differences of the original dynamics of problem~\eqref{eq:ddp_obj_detailed}
	rather than the linearized dynamics. Formally, given a sequence of policies
	$(\pi_t)_{t=0}^{\horizon-1}$, the oracles are then given as $\diffctrls =
	(\diffctrl_0;\ldots;\diffctrl_{\horizon-1})$ computed, for $\diffstate_0 =0$,
	by Algo.~\ref{algo:roll_out} as
	\[
	\diffctrls = \Roll(\diffstate_0,( \pi_t)_{t=0}^{\horizon-1},  (\diff_{\dyn_t} ^{\state_t, \ctrl_t})_{t=0}^{\horizon-1}),
	\]
	where $\diff_{\dyn_t} ^{\state_t, \ctrl_t}(\auxstate_t, \auxctrl_t) =
	\dyn_t(\state_t + \diffstate_t, \ctrl_t + \diffctrl_t) - \dyn_t(\state_t,
	\ctrl_t)$.
	
\end{enumerate}

	\section{Computational Complexity}\label{sec:comput_cplxity}
	\begin{figure}
	\begin{center}
		\begin{tikzpicture}
			\node(dyn) {\small \bf Dyn approx.};
			\node(cost) [right=2.5em of dyn] {\small \bf Cost approx.};
			\node(forward) [above right=1ex and -2em of dyn, align=center]{\bf Forward pass\\
			\small{Algo.~\ref{algo:forward}}} ;
			\node(bellman) [right=57pt of forward, align=center] {\bf Backward pass\\ 
			\small{~}};
			\node(rollout) [right=40pt of bellman, align=center] {\bf Roll-out\\ \small{Algo.~\ref{algo:roll_out}}};
			\node(algo) [right=6.5em of rollout, align=center] {\bf Oracle\\
			\small{~}};
			
			\node(dyn1) [below =1em of dyn] {1st order};
			\node(dyn2) [below =9em of dyn1] {2nd order};
			
			\node(cost1) [right=4em of dyn1] {1st order};
			\node(cost2) [below =2em of cost1] {2nd order};
			\node(cost3) [right =4em of dyn2] {2nd order};
			
			\node(bell1) [right =4.5em of cost1, align=center] {$\Backward_{\operatorname{GD}}$ \\
				\small{Algo.~\ref{algo:backward_grad}}};
			\node(bell2) [right =4.4em of cost2,align=center] {$\Backward_{\operatorname{GN}}$\\ 
				\small{Algo.~\ref{algo:backward_gn}}};
			\node(bell3) [right =4.4em of cost3, align=center] {$\Backward_{\operatorname{NE}}$\\ 
				\small{Algo.~\ref{algo:backward_newton}}};
			\node(bell4) [below =1em of bell3, align=center] {$\Backward_{\operatorname{DDP}}$\\ 
				\small{Algo.~\ref{algo:backward_ddp}}};
			
			\node(roll1) [right=6em of bell1]{None};
			\node(roll2) [right=4em of bell2]{Linearized dyn.};
			\node(roll3) [below=1.5em of roll2]{Original dyn.};
			\node(roll4) [right=4.5em of bell3]{Linearized dyn.};
			\node(roll5) [right=4.5em of bell4]{Original dyn.};
			
			\node (algo1) [right =75pt of roll1, align=center] {
			{~}\\
			GD\\
			\small{Algo.~\ref{algo:gd}}};
			\node (algo2) [right =40pt of roll2, align=center] {
			{~}\\
			GN (ILQR)\\
			\small{Algo.~\ref{algo:gn}}};
			\node (algo3) [right=35pt of roll3, align=center] {
			{~}\\
			DDP-LQ (iLQR)\\
			\small{Algo.~\ref{algo:ddp_lq}}};
			\node (algo4) [right=43pt of roll4, align=center] {
			{~}\\
			NE\\
			\small{Algo.~\ref{algo:newton}}};
			\node (algo5) [right=35pt of roll5, align=center] {
			{~}\\
			DDP-Q (DDP)\\
			\small{Algo.~\ref{algo:ddp_q}}};
			
			\draw (dyn1) -- (cost1) -- (bell1) -- (roll1) -- (algo1);
			\draw (dyn1) -- (cost2) -- (bell2) -- (roll2) -- (algo2);
			\draw (dyn1) -- (cost2) -- (bell2) -- (roll3) -- (algo3);
			\draw (dyn2) -- (cost3) -- (bell3) -- (roll4) -- (algo4);
			\draw (dyn2) -- (cost3) -- (bell4) -- (roll5) -- (algo5);
		\end{tikzpicture}
	\end{center}
	\caption{Taxonomy of non-linear control oracles. GD stands for gradient Descent, GN for Gauss-Newton, NE for Newton, DDP-LQ and DDP-Q stand for DDP  with  linear quadratic  or quadratic approx. The iterations of the algorithms use a line-search procedure presented in Algo.~\ref{algo:line_search} as illustrated in Algo.~\ref{algo:gn_algo}.\label{fig:taxonomy}}
\end{figure}

In Figure~\ref{fig:taxonomy}, we present a summary of the different algorithms
presented in this manuscript. We added in parentheses the names usually given
for these methods. Additional line-search mechanisms are presented in
Appendix~\ref{app:line_searches}. The overall implementations are detailed in
Appendix~\ref{app:summary}. We consider then the computational complexities of
the algorithms in a differentiable programming framework. 

{\bf Formal Computational Complexity.}~
We present in~Table~\ref{tab:oracle_complexities} the computational complexities
of the algorithms following the implementations described in
Sec.~\ref{sec:classical_optim} and Sec.~\ref{sec:ddp} and detailed in
Appendix~\ref{app:summary}. We ignore the additional cost of the line-searches
which requires a theoretical analysis of the admissible stepsizes depending on
the smoothness properties of the dynamics and the costs. We consider for
simplicity that the cost of evaluating a function $f:\reals^d \rightarrow
\reals^n$ is of the order of $O(n d)$, as it is the case if $f$ is linear. For
the computational complexities of the core operation of the backward pass, i.e,
$\BellLQ$ in Algo.~\ref{algo:BellLQ} or $\BellL$ in Algo.~\ref{algo:bp}, we
simply give the leading computational complexities, which, in the case of
$\BellLQ$, are the matrix multiplications and inversions.
  
The time complexities differ depending on whether linear or quadratic
approximations of the costs are used. In the latter case, matrices of size
$\dimctrl \times \dimctrl$ need to be inverted and matrices of size
$\dimstate\times \dimstate$ need to be multiplied. However, all oracles have a
linear time complexity with respect to the horizon $\horizon$.

We note that the space complexities of the gradient descent and the Gauss-Newton
method or the DDP approach with linear quadratic approximations are essentially
the same. On the other hand, the space complexity of the Newton oracle is a
priori larger.

{\bf Computational Complexity in a Differentiable Programming Framework.}~ The
decomposition of each oracle between forward, backward and roll-out passes has
the advantage to clarify the discrepancies between each approach. However, a careful implementation of these
oracles only requires storing in memory the function and the inputs given at
each time-step. Namely, the forward pass can simply keep in memory  $\cost_t,
\dyn_t, \state_t, \ctrl_t$ for $t \in \{0, \ldots, \horizon\}$. The backward
pass computes then, on the fly, the information necessary to compute the
policies. This amounts to a simple system of check-pointing, a strategy used in
differentiable programming to circumvent the memory cost of the reverse-mode of
automatic differentiation~\citep{griewank2008evaluating}.

Such an approach is detailed in Appendix~\ref{app:cplxity}. In summary, by
considering an implementation that simply stores in memory the inputs and the
programs that implement the functions, a Newton oracle and an oracle based on a
DDP approach with quadratic approximation have the same time and space
complexities as their linear quadratic counterparts up to constant factors. This
remark was done by~\citet{nganga2021accelerating} for implementing a DDP
algorithm with quadratic approximations.

\begin{table}
	\begin{center}
		\bgroup
		\def\arraystretch{1.5}

					\begin{tabularx}{\linewidth}{p{105pt}|p{280pt}}
			\multicolumn{1}{c}{}	  & \multicolumn{1}{l}{Time complexities of the forward pass in Algo.~\ref{algo:forward}}\\
			\toprule
			Function eval. \newline{\footnotesize ($o_\dyn=o_\cost =0$)}
			& $\horizon\Big(
			\underbrace{\dimstate^2 {+} \dimstate\dimctrl}_{\dyn_t} 
			{+} \underbrace{\dimstate {+} \dimctrl}_{\cost_t}
			\Big){=}O(\horizon(\dimstate^2 {+} \dimstate\dimctrl))$
			\\
			Lin. (GD)  \newline{\footnotesize ($o_\dyn=o_\cost =1$)}
			& $\horizon\Big(
			\underbrace{\dimstate^2 {+} \dimstate\dimctrl}_{\dyn_t, \nabla \dyn_t }
			{+} \underbrace{\dimstate {+} \dimctrl}_{\cost_t, \nabla \cost_t}
			\Big){=} O(\horizon(\dimstate^2 {+} \dimstate\dimctrl))$
			\\
			Lin.-quad. (GN/DDP-LQ)   \newline{\footnotesize ($o_\dyn=1, o_\cost =2$)}
			& $\horizon\Big(
			\underbrace{\dimstate^2 {+} \dimstate\dimctrl}_{\dyn_t, \nabla \dyn_t} 
			{+} \underbrace{\dimstate {+} \dimctrl}_{\cost_t, \nabla \cost_t} {+} 
			\underbrace{\dimstate^2 {+} \dimctrl^2 {+} \dimstate\dimctrl}_{\nabla^2 \cost_t}
			\Big) {=} O(\horizon(\dimstate {+} \dimctrl)^2)$
			\\
			Quad. (NE/DDP-Q)   \newline{\footnotesize ($o_\dyn=o_\cost =2$)}
			& $\horizon\Big(
			\underbrace{\dimstate^2 {+} \dimstate\dimctrl}_{\dyn_t, \nabla \dyn_t}  
			{+} \underbrace{(\dimstate^2 {+} \dimctrl^2 {+} \dimstate\dimctrl)\dimstate}_{\nabla^2 \dyn_t} 
			{+} \underbrace{\dimstate {+} \dimctrl}_{\cost_t, \nabla \cost_t} {+} 
			\underbrace{\dimstate^2 {+} \dimctrl^2 {+} \dimstate\dimctrl}_{\nabla^2 \cost_t}
			\Big){=} O(\horizon\dimstate(\dimstate {+} \dimctrl )^2)$
			\\
			\bottomrule
		\end{tabularx}
	
	\vspace{1ex}
	
					\begin{tabularx}{\linewidth}{p{105pt}|p{280pt}}
	\multicolumn{1}{c}{}	  & \multicolumn{1}{l}{Space complexities of the forward pass in Algo.~\ref{algo:forward}}\\
	\toprule
	Function eval. \newline{\footnotesize ($o_\dyn=o_\cost =0$)}
	& $0$
	\\
	Lin. (GD)  \newline{\footnotesize ($o_\dyn=o_\cost =1$)}
	& $\horizon\Big(
	\underbrace{\dimstate^2 {+} \dimstate\dimctrl}_{\nabla \dyn_t }
	{+} \underbrace{\dimstate {+} \dimctrl}_{\nabla \cost_t}
	\Big){=} O(\horizon(\dimstate^2 {+} \dimstate\dimctrl))$
	\\
	Lin.-quad. (GN/DDP-LQ)   \newline{\footnotesize ($o_\dyn=1, o_\cost =2$)}
	& $\horizon\Big(
	\underbrace{\dimstate^2 {+} \dimstate\dimctrl}_{\nabla \dyn_t} 
	{+} \underbrace{\dimstate {+} \dimctrl}_{\nabla \cost_t} {+} 
	\underbrace{\dimstate^2 {+} \dimctrl^2 {+} \dimstate\dimctrl}_{\nabla^2 \cost_t}
	\Big) {=} O(\horizon(\dimstate{+} \dimctrl)^2)$
	\\
	Quad. (NE/DDP-Q)  \newline{\footnotesize ($o_\dyn=o_\cost =2$)}
	& $\horizon\Big(
	\underbrace{\dimstate^2 {+} \dimstate\dimctrl}_{\nabla \dyn_t}  
	{+} \underbrace{(\dimstate^2 {+} \dimctrl^2 {+} \dimstate\dimctrl)\dimstate}_{\nabla^2 \dyn_t} 
	{+} \underbrace{\dimstate {+} \dimctrl}_{ \nabla \cost_t} {+} 
	\underbrace{\dimstate^2 {+} \dimctrl^2 {+} \dimstate\dimctrl}_{\nabla^2 \cost_t}
	\Big){=} O(\horizon\dimstate(\dimstate {+} \dimctrl )^2)$
	\\
	\bottomrule
\end{tabularx}

\vspace{1ex}
					\begin{tabularx}{\linewidth}{p{60pt}|p{280pt}}
			\multicolumn{1}{c}{} & \multicolumn{1}{l}{Time complexities of the backward passes in Algo.~\ref{algo:backward_grad}, \ref{algo:backward_gn}, \ref{algo:backward_newton}, \ref{algo:backward_ddp} and the roll-out in Algo.~\ref{algo:roll_out}}\\
			\toprule
			GD 
			& $\horizon\Big(
			\underbrace{\dimstate^2 {+} \dimstate\dimctrl}_{\textrm{Roll}} 
			{+} \underbrace{\dimstate^2 {+} \dimstate\dimctrl
			 }_{\BellL}
			\Big){=}O(\horizon(\dimstate^2 {+} \dimstate\dimctrl))$ 
			\\
			GN/DDP-LQ
			& $\horizon\Big(
			\underbrace{\dimstate^2 {+} \dimstate\dimctrl}_{\textrm{Roll}} 
			{+} \underbrace{\dimstate^3 {+} \dimctrl^3 {+} \dimctrl^2\dimstate}_{\BellLQ}
			\Big){=} O(\horizon(\dimstate {+} \dimctrl)^3)$
			\\
			NE/DDP-Q 
			& $\horizon\Big(
			\underbrace{\dimstate^2 {+} \dimstate\dimctrl}_{\textrm{Roll}} 
			{+} \underbrace{\dimstate^3 {+} \dimctrl^3 {+} \dimctrl^2\dimstate}_{\BellLQ} 
			{+} \underbrace{(\dimstate^2 {+} \dimctrl^2 {+} \dimstate\dimctrl)\dimstate}_{\nabla\dyn_t^2[\cdot, \cdot, \lambda]}\Big) {=} O(\horizon(\dimstate {+} \dimctrl)^3)$
			\\
			\bottomrule
		\end{tabularx}

		\egroup
		\caption{Space and time complexities of the oracles of Sec.~\ref{sec:classical_optim} and~\ref{sec:ddp}. 
			Acronyms are given in Fig.~\ref{fig:taxonomy}.
			\label{tab:oracle_complexities}}
    \vspace*{-20pt}
	\end{center}
\end{table}

	\section{Experiments}\label{sec:exp}
  The control environments considered are thoroughly described in
Appendix~\ref{app:exp_details}. The code is publicly available at~{\small
\coderef}. Additional experiments are presented in Appendix~\ref{app:exp_sup}, a
comparison of all algorithms is presented in Fig.~\ref{fig:conv_all_iteration}.

\revised{ All the following plots are in log-scale where on the vertical axis we
plot $\log\left(\obj(\currvar)/\obj(\var^{(0)})\right)$ with $\obj$ the
objective, $\currvar$  the set of controls at iteration $k$. The acronyms (GD,
GN, NE, DDP-LQ, DDP-Q) correspond to the taxonomy of algorithms presented in
Fig.~\ref{fig:taxonomy}. Finally, the algorithms are stopped if no valid
stepsizes have been found by line-search beyond machine precision $\varepsilon$,
or if the relative difference in terms of costs is smaller than machine
precision, where $\varepsilon \approx 10^{-16}$ as we ran these experiments in
double precision. Hence, if a curve stops, this means that the linesearch did
not find a valid stepsize beyond this point.}

\subsection{Linear Quadratic Approximation}
\begin{figure}
	\begin{center}
    \includegraphics[width=0.9\linewidth]{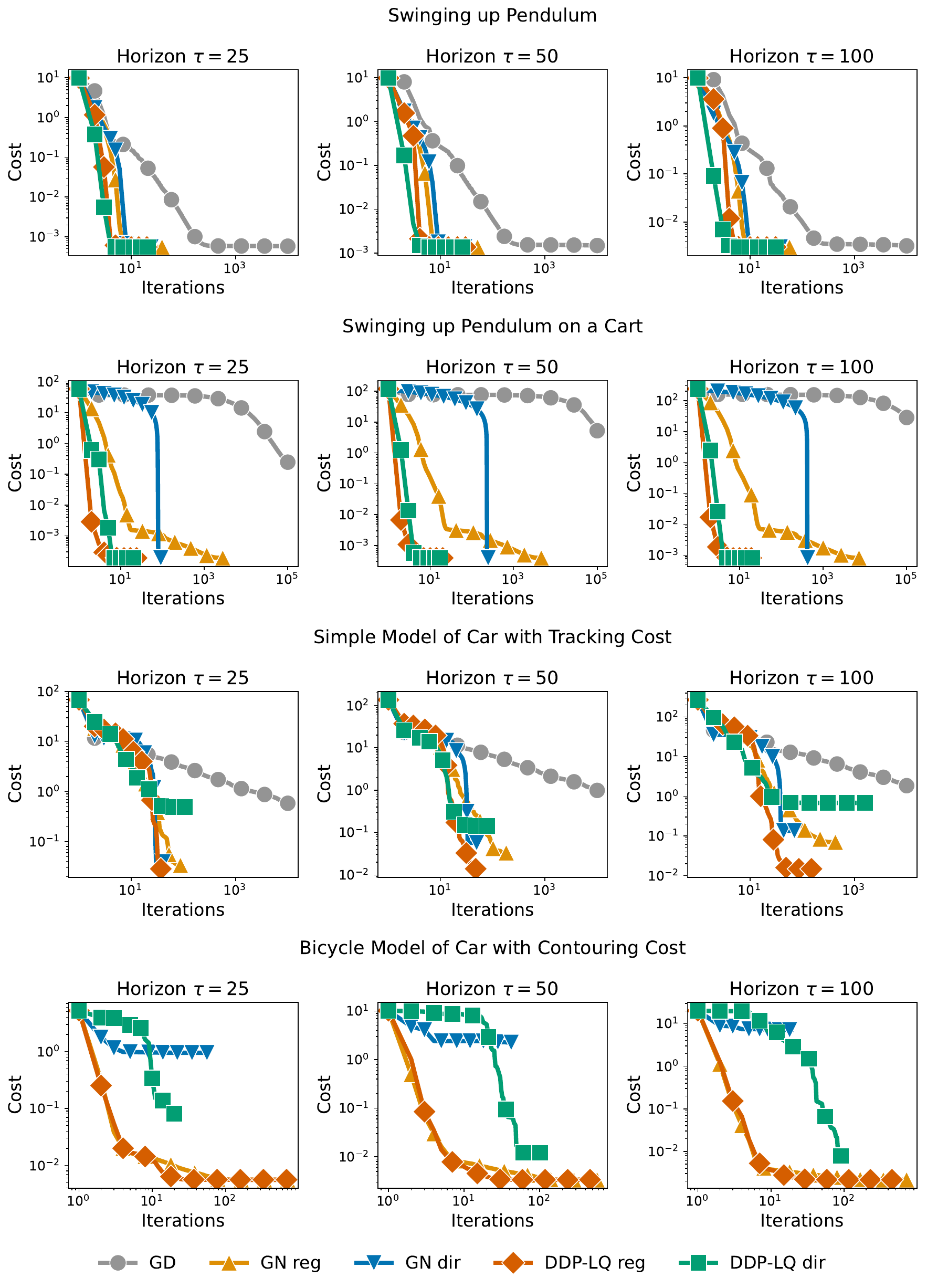}
    \caption{
    \revised{
		Cost along iterations on various control problems detailed in
    Appendix~\ref{app:exp_details} with algorithms using linear (GD) or
    linear-quadratic approximations (GN, DDP-LQ, see Fig.~\ref{fig:taxonomy} for
    taxonomy details) and directional (dir \eqref{eq:linesearch_armijo}) or
    regularized (reg \eqref{eq:linesearch_reg}) steps.
    \label{fig:conv_linquad}}
    }
	\end{center}
\end{figure}

\revised{
We compare first a gradient descent and nonlinear
control algorithms with linear quadratic approximations, i.e., GN or DDP-LQ with
directional or regularized steps. We make the following observations.
\begin{enumerate}
  \item Cost along iterations (Fig.~\ref{fig:conv_linquad})
  \begin{itemize}
    \item GN and DDP-LQ always outperform GD by several order of
    magnitudes.
    \item DDP-LQ generally performs better or on par with 
    GN, for the same steps (directional or regularized).
    \item For GN, regularized steps generally provide a more steady convergence
    than directional steps. The later do not find a valid stepsize in the real car
    example, which involves highly nonlinear dynamics (see
    Appendix~\ref{app:exp_details}). However, once in a quadratically convergent
    phase the directional steps can provide faster convergence than the regularized
    steps  (see e.g. GN dir vs GN reg on the pendulum on a cart example).
    \item For DDP-LQ, similar observations can be done. Only on the simple pendulum problem, 
    directional steps slightly outperform regularized steps
  \end{itemize}
  \item Cost along computational time (Fig.~\ref{fig:conv_linquad_time})
  \begin{itemize}
    \item In terms of time, regularized steps may require fewer evaluations
    during the line-search as they incorporate previous stepsizes and may provide
    faster convergence in time.
  \end{itemize}
  \item Gradient norm along iterations (Fig.~\ref{fig:norm_grad_linquad})
  \begin{itemize}
    \item Algorithms based on linear quadratic approximations generally display
    a late quadratic convergence in terms of gradient norm. One exception is the
    realistic model of the car, where only the regularized steps versions are
    able to make substantial progress along the iterations and such progress is
    only linear (in log-log plot scale).
  \end{itemize}
\end{enumerate}
}

\subsection{Quadratic Approximation}
\begin{figure}
	\begin{center}
    \includegraphics[width=0.9\linewidth]{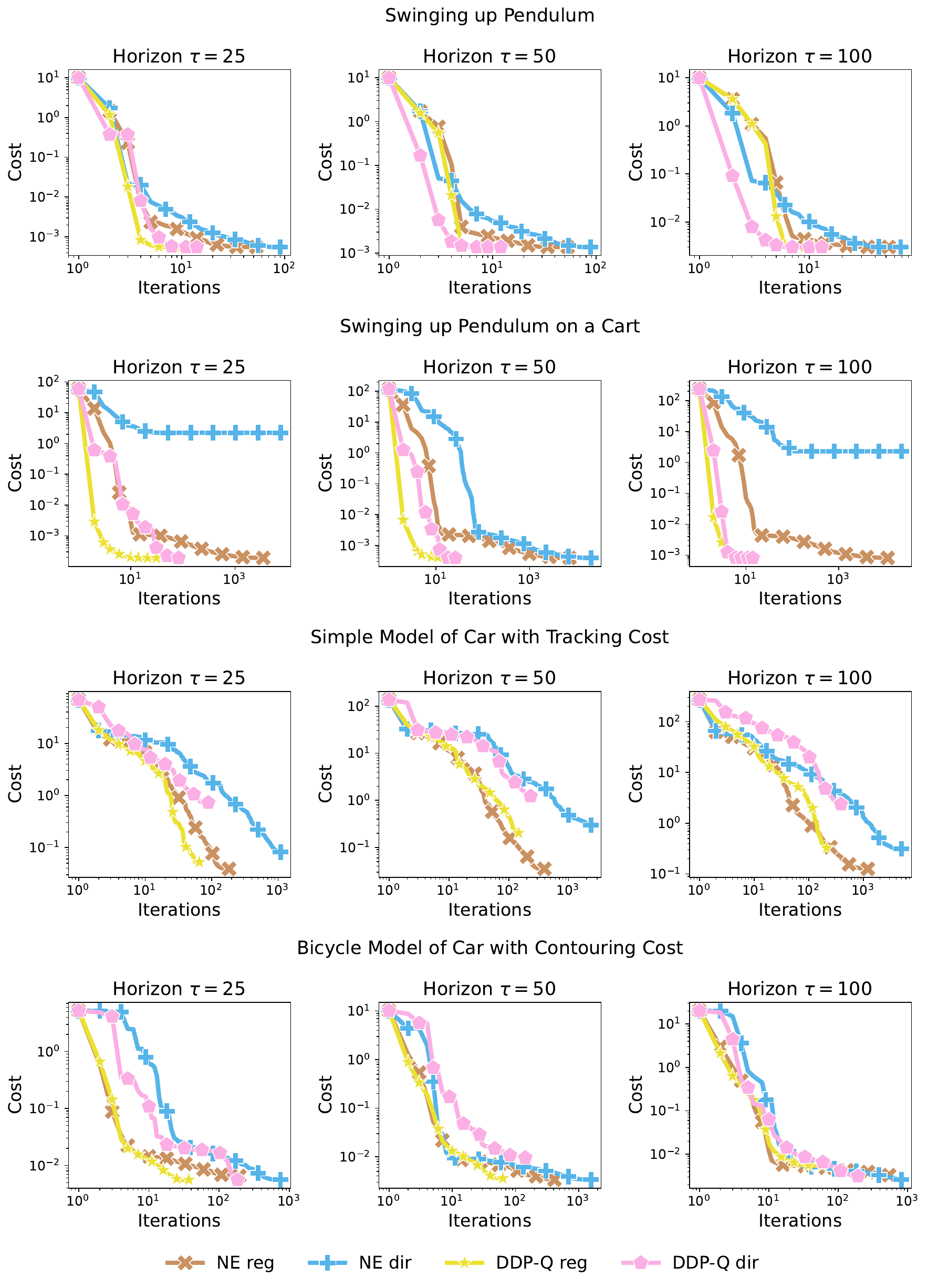}
    \caption{
    \revised{
		Cost along iterations on various control problems detailed in
		Appendix~\ref{app:exp_details} with algorithms using quadratic
		approximations (NE, DDP-Q, see Fig.\ref{fig:taxonomy} for taxonomy details)
		and directional (dir \eqref{eq:linesearch_armijo}) or regularized (reg
		\eqref{eq:linesearch_reg}) steps. \label{fig:conv_quad}}
    }
	\end{center}
\end{figure}

We compare now nonlinear control algorithms with
quadratic approximations, i.e., NE or DDP-Q.
\revised{
\begin{enumerate}
  \item Cost along iterations (Fig.~\ref{fig:conv_quad})
  \begin{itemize}
    \item As for the linear-quadratic approximations, the DDP
    approach (here DDP-Q) generally outperforms or performs on par with its
    Newton (NE) counterpart.
    \item For Newton, the regularized steps generally outperform 
    the directional steps.
    \item For DDP-Q, the regularized steps outperform the directional steps on
    all examples but the first pendulum examples.
  \end{itemize}
  \item Cost along computational time (Fig.~\ref{fig:conv_quad_time})
  \begin{itemize}
    \item In terms of time, all algorithms appear to generally perform on par.
  \end{itemize}
  \item Gradient norm along iterations (Fig.~\ref{fig:norm_grad_quad})
  \begin{itemize}
    \item We generally observe quadratic convergence in gradient norm 
    for all algorithms in late stage of training. However, quadratic convergence 
    of Newton generally appears later than DDP.
  \end{itemize}
\end{enumerate}
}

\clearpage

  \clearpage

	{\bf Acknowledgments.} 
  This work was supported by NSF DMS-1839371, DMS-2134012, CCF-2019844, CIFARLMB,
  NSF TRIPODS II DMS-2023166 and faculty research awards. The authors deeply thank
  Alexander Liniger for his help on implementing the bicycle model of a car. The
  authors also thank Dmitriy Drusvyatskiy, Krishna Pillutla and John Thickstun for
  fruitful discussions on the paper and the code. 

	\clearpage
	\appendix

  {\bf \Huge\selectfont\sffamily\bfseries Appendix}\\

  The appendix is organized as follows.
  \begin{enumerate}
    \item Appendix~\ref{app:notations} presents tensor notations used to describe
    algorithms with second-order information on the dynamics.
    \item Appendix~\ref{app:related_work} expands the discussion of related work.
    \item Appendix~\ref{app:proofs} details the proofs of the results claimed in the main text.
    \item Appendix~\ref{app:line_searches} presents line-search mechanisms incorporated
    in the algorithms to ensure their efficiency.
    \item Appendix~\ref{app:summary} details all pseudocode algorithms with associated computational graphs in a differentiable programming framework.
    \item Appendix~\ref{app:cplxity} gives additional complexities when implementing the
    oracles with checkpointing.
    \item Appendix~\ref{app:sparse_solvers} presents alternative ways to compute oracles using the structure of the subproblems.
    \item Appendix~\ref{app:exp_details} details the control environments on which the algorithms are tested.
    \item Appendix~\ref{app:exp_sup} presents additional experiments: the convergence of the algorithms in time, and the stepsize selected along the iterations by a line-search procedure.
  \end{enumerate}

  \section{Tensor Notation}\label{app:notations}
  A tensor $\mathcal{A} = (a_{i,j,k})_{1\leq i\leq d, 1\leq j\leq p, 1\leq k\leq
n} \in \reals^{d \times p \times n}$ is represented as a list of matrices
$\mathcal{A} = (A_1,\ldots, A_n)$ where $A_k = (a_{i,j,k})_{1 \leq i \leq d, 1
\leq j \leq p} \in \reals^{d\times p}$ for $ k\in \{1,\ldots n\}$. Given
$\mathcal{A}  \in \reals^{d\times p\times n}$ and  $P \in \reals^{d \times d'},
Q \in \reals^{p \times p'}, R \in \reals^{n \times n'}$, we denote
\[
\mathcal{A}[P, Q, R] \coloneqq
\left(\sum_{k=1}^{n} R_{k,1}P^\top A_k Q,  \ldots,  \sum_{k=1}^{n} R_{k,n'}P^\top A_k Q \right) \in \reals^{d'\times p'\times n'}.
\]
For $\mathcal{A}_0 \in \reals^{d_0\times p_0\times n_0}$, $P \in \reals^{d_0
\times d_1}, Q \in \reals^{p_0 \times p_1}, R \in \reals^{n_0 \times n_1}$
denote $\mathcal{A}_1= \mathcal{A}_0[P, Q, R] \in \reals^{d_1\times p_1\times
n_1}$. Then, for $S \in \reals^{d_1\times d_2}, T \in \reals^{p_1 \times p_2}, U
\in \reals^{ n_1\times n_2}$, we have 
$
\mathcal{A}_1[S, T, U] = \mathcal{A}_0[PS, QT, RU] \in \reals^{d_2\times p_2\times n_2}.
$
If $P, Q$ or $ R$ are identity matrices, we use the symbol ``$\: \cdot\: $'' in
place of the identity matrix. For example, we denote $\mathcal{A}[P, Q, \idm_n]
= \mathcal{A}[P,Q, \cdot] = \left(P^\top A_1 Q,  \ldots,  P^\top A_n Q \right)$.
If $P, Q$ or $R$ are vectors we consider the flattened object. In particular, for
$x\in \reals^d, y\in \reals^p$, we denote
$
\mathcal{A}[x, y, \cdot] =  \left(\begin{matrix}
	x^\top A_1 y,\ldots, x^\top A_ny
\end{matrix}
\right)^\top\in \reals^n,
$
rather than having $\mathcal{A}[x, y, \cdot] \in \reals^{1 \times 1\times n}$.
Similarly, for $z\in \reals^n$, we denote
$
\mathcal{A}[\cdot, \cdot, z] = \sum_{k=1}^n z_kA_k \in \reals^{d\times p}.
$ 
Such notations follow the ones used by \citet[Chapter 5]{nesterov2018lectures}
to study third-order derivatives.

  \section{Related Work}\label{app:related_work}
  Nonlinear control problems of the form~\eqref{eq:discrete_ctrl_pb} stem from the
discretization of generic optimal control problems in continuous time of the
form
\begin{align}
  \min_{\state(\cdot), \ctrl(\cdot)} \quad
  & \int_0^T h(x(t), u(T)) + h_T(x(T)) \label{eq:cont_ctrl_pb}\\
  \mbox{subject to} \quad  & \dot x(t) = f(x(t), u(t)), \quad x(0) = \bar x_0, \nonumber
\end{align}
Continuous optimal control problems of the form~\eqref{eq:cont_ctrl_pb} can be
tackled in various ways~\citep{diehl2006fast}. One can approach the problem from
a \emph{dynamic programming} perspective to derive the Hamilton-Jacobi-Bellman
equation, a partial differential equation in state
space~\citep{lions1982generalized}. Alternatively, one can derive necessary
optimality conditions for~\eqref{eq:cont_ctrl_pb} to derive a boundary value
problem. Such a method is referred to as an \emph{indirect method} and amounts
to a ``optimize then discretize'' approach~\citep{farshidian2017efficient}.
Finally, problem~\eqref{eq:cont_ctrl_pb} can be tackled by \emph{direct methods}
that consider finite dimensional approximations of the original infinite
dimensional problem~\eqref{eq:cont_ctrl_pb}. Direct methods amount to a
``discretize then optimize'' approach~\citep{diehl2006fast}, they can further be
split into different approaches. First, one may consider a finite representation
of the continuous control $u(t)$ as piecewise constant functions whose values
$q_1, \ldots, q_\tau$ at each piece define the finite number of degrees of
freedom. The problem still involves an ODE in the state variable, $\dot x(t) =
f(x(t), u_{q_{1:\tau}}(t))$, albeit a simpler one. Tackling the problem with
such a partial discretization is referred to as a \emph{single shooting}
method~\citep{diehl2006fast, bock1984multiple}. \emph{Collocation
methods}~\citep{von1993numerical} consider discretizing both the states and
controls, leading to a formulation like~\eqref{eq:discrete_ctrl_pb}, that can
benefit from advanced numerical integration methods. Finally, \emph{multiple
shooting}~\citep{diehl2006fast,bock1984multiple} combines both approaches. The
system is split in multiple windows and for each window a single shooting method
is used. We focus solely on the resulting discrete time nonlinear control
problems~\eqref{eq:discrete_ctrl_pb} and refer the interested reader to,
e.g.,~\citet{diehl2006fast} for an overview of the approaches mentioned above.

One of the first approaches for nonlinear discrete time control
problems~\eqref{eq:discrete_ctrl_pb} appear to be the Differential Dynamic
Programming (DDP) methods developed by~\citet{jacobson1970differential} and
further explored by~\citet{mayne1975first,murray1984differential,
liao1991convergence}. \citet{bock1984multiple} referred to such approaches as
\emph{direct multiple shooting}. Numerous variants of DDP have been developed to
account for constraints or noise in the dynamics~\citep{li2007iterative,
tassa2007receding, tassa2014control, giftthaler2018family}.

An implementation of a Newton method for nonlinear control problems of the
form~\eqref{eq:discrete_ctrl_pb} was developed after the DDP approach
by~\citet{de1988differential, dunn1989efficient}. A parallel implementation of a
Newton step and sequential quadratic programming methods were developed
by~\citet{wright1990solution, wright1991partitioned}, which led to efficient
implementations of interior point methods for linear quadratic control problems
under constraints by using the block band diagonal structure of the system of
KKT equations solved at each step~\citep{wright1991structured}. A detailed
comparison of the DDP approach and the Newton method was conducted
by~\citet{liao1992advantages}, who observed that the original DDP approach
generally outperforms its Newton counterpart. We extend this analysis by
comparing  regularized variants of the algorithms. Finally, the storage of
second order information for DDP and Newton can be alleviated with a careful
implementation in a differentiable programming framework as done in our
implementation and noted earlier by~\citet{nganga2021accelerating}. 

Simpler approaches consisting in taking linear approximations of the dynamics
and quadratic approximations of the costs were implemented as part of public
software~\citep{todorov2012mujoco}. Two variants have been presented. The
Iterative Linear Quadratic Regulator (ILQR) algorithm as originally formulated
by~\citet{li2007iterative} amounts naturally to a Gauss-Newton
method~\citep{sideris2005efficient}. A variant that mixes linear quadratic
approximations of the problem with a DDP approach, named iterative Linear
Quadratic Regulator (iLQR) was further analyzed empirically
by~\citet{tassa2012synthesis}. Here, we detail the line-searches for both
approaches and present their regularized variants. We provide detailed
computational complexities of all aforementioned algorithms that illustrate the
trade-offs between the approaches.

Nonlinear model predictive control methods generally use sparse linear algebra
solvers at each iteration~\citep{diehl2009efficient} using solvers like
IPOPT~\citep{wachter2006implementation} or SNOPT~\citep{gill2005snopt}. For
offline control problems like~\eqref{eq:discrete_ctrl_pb}, such sparse linear
algebra solvers can also be used to compute the Gauss-Newton or Newton oracles
seen as the solutions of a linear quadratic problem with underlying sparse band
diagonal structure as first observed
by~\citet{wright1990solution,wright1991partitioned}. These sparse linear algebra
solvers are an alternative to the dynamic programming procedures, presented in
this manuscript, that can be seen as solving Riccatti equations in discrete-time
with finite horizon. On the other hand, these sparse linear algebra solvers
cannot be used as a black-box to implement DDP methods since they output
directly the control variables solutions of the subproblem and do not a priori
give access to the policies. They can nevertheless be adapted to record
policies~\citep{Verschueren2021, jallet2023proxddp}. In this manuscript, we cast
both classical optimization oracles and DDP approaches in a common
differentiable programming framework to highlight their common ground and
discrepancies, which would not be possible from a purely algebraic viewpoint. We
aim at comparing these approaches purely in terms of iterations to understand
differences in behavior, and leave out the optimization of these implementations
in specific frameworks, using e.g. sparse linear algebra solvers to implement
each classical optimization oracle. This viewpoint was generalized to handle
nonlinear inequalities in model predictive control~\citep{diehl2009efficient} or
even generic graphs of computations~\citep{srinivasan2015graphical}. Alternative
methods cast as sequential quadratic programming
techniques~\citep{messerer2021survey, frasch2015parallel,
verschueren2016exploiting, houska2013quadratically} are also worth mentioning. 

For our experiments, we adapted the bicycle model of a miniature car developed
by~\citet{liniger2015optimization} in Python. We provide an implementation in
Python,  available at {\small\coderef}  for further exploration of the
algorithms. Similar implementations have been implemented in the trajax
library~\citep{jax2018github}. Numerous other packages exist to implement
nonlinear control algorithms such as CasAdi~\citep{andersson2018cassadi},
Pyomo~\citep{bynum2021pyomo}, JumP~\citep{dunning2017jump},
acados~\citep{Verschueren2021} that can take advantage of off-the-shelf interior
point solvers such as IPOPT~\citep{wachter2006implementation}, or
SNOPT~\citep{gill2005snopt}. Recently, \citet{bambade2022prox} developed a new
solver for quadratic programs with linear constraints using augmented
Lagrangian. This solver in turn led to new efficient nonlinear control
algorithms such as prox-DDP~\citep{jallet2023proxddp}.

  \section{Proofs}\label{app:proofs}
  
This section gathers proofs of propositions given in the main text.
\subsection{Linear Quadratic Control}
\begin{lemm}
For linear functions $\lin_t$ and quadratic functions $\qua_t, \costogo_{t+1}$ s.t. $\qua_t(\state, \cdot) + \costogo_{t+1}(\lin_t(\state, \cdot))$ is
strongly convex for any $x$, the procedure
\begin{equation}
\BellLQ	: (\lin_t, \qua_t, \costogo_{t+1})  \rightarrow \left(\begin{array}{c}
    \costogo_t:\state\rightarrow\min_{\ctrl \in \reals^\dimctrl} \left\{ \qua_t(\state, \ctrl) + \costogo_{t+1}(\lin_t(\state, \ctrl))\right\}  \\
    \pi_t :\state\rightarrow \argmin_{\ctrl \in \reals^\dimctrl}  \left\{ \qua_t(\state, \ctrl) + \costogo_{t+1}(\lin_t(\state, \ctrl))\right\}
\end{array}\right), \nonumber
\end{equation}
can be implemented analytically as detailed in Algo.~\ref{algo:BellLQ}.
\end{lemm}
\begin{proof} Consider $\lin_t, \qua_t, \costogo_{t+1}$ to be parameterized as
  $\lin_t(\state, \ctrl) = \A \state + \B\ctrl$, $\qua_t(\state, \ctrl) =
  \frac{1}{2}\state^\top \H \state+  \frac{1}{2}\ctrl^\top \G \ctrl +
  \state^\top \R \ctrl+ \h^\top \state  +  \g^\top \ctrl$, $
  \costogo_{t+1}(\state) = \frac{1}{2}\state^\top \J_{t+1}\state + \j_{t+1}^\top
  \state + \jcst_{t+1}.$ The cost-to-go function at time $t$ is
  \begin{align*} 
  \costogo_t(\state) &  = \frac{1}{2}\state^\top \H \state +
  \h^\top \state  + \jcst_{t+1}\\
    & \quad + \min_{\ctrl \in \reals^\dimctrl}\left\{  \frac{1}{2}(\A \state + \B
  \ctrl)^\top \J_{t+1}(\A \state+ \B \ctrl) + \j_{t+1}^\top (\A \state + \B
  \ctrl) + \frac{1}{2} \ctrl^\top \G \ctrl + \state^\top \R \ctrl +  \g^\top
  \ctrl\right\}.
  \end{align*} 
  Since $\cost(\state, \cdot) +
  \costogo_{t+1}(\lin(\state, \cdot))$ is strongly convex, we have that $\G +
  \B^\top \J_{t+1}\B \succ 0$. Therefore, the policy at time $t$ is
  \[
  \pi_t(\state) =   -(\G + \B^\top\J_{t+1}\B)^{-1} [(\R^\top + \B^\top
  \J_{t+1} \A) \state + \g + \B^\top \j_{t+1}].
  \]
  Using that $\min_{u \in \reals^\dimctrl} u^\top M u/2 + m^\top x = -m^\top M^{-1}m/2$ where, here,
  $M=\G + \B^\top \J_{t+1}\B$, $m=(\R^\top + \B^\top \J_{t+1} \A) \state + \g +
  \B^\top \j_{t+1}$, we get that the cost-to-go function at time $t$ is given by
  \begin{align*} \costogo_t(\state) & = \frac{1}{2}\state^\top \left(\H +
  \A^\top \J_{t+1}\A - (\R + \A^\top \J_{t+1}\B) (\G + \B^\top \J_{t+1}\B)^{-1}
  (\R^\top + \B^\top \J_{t+1} \A)\right) \state \\
    & \quad + \left(\h + \A^\top \j_{t+1} - (\R + \A^\top \J_{t+1}\B) (\G +
    \B^\top \J_{t+1}\B)^{-1} (\g+ \B^\top \j_{t+1})\right)^\top \state \\
    & \quad - \frac{1}{2}(\g + \B^\top \j_{t+1})^\top (\G +
  \B^\top\J_{t+1}\B)^{-1} ( \g + \B^\top \j_{t+1}) + \jcst_{t+1}. \end{align*}
\end{proof}

\begin{coro}
  Consider problem~\eqref{eq:discrete_ctrl_pb} such that for all $t\in \{0,
  \ldots, \horizon-1\}$, $\dyn_t$ is linear, $\cost_t$ is convex quadratic with
  $\cost_t(\state, \cdot)$ strongly convex for any $\state$, and
  $\cost_\horizon$ is convex quadratic. Then, the solution of
  problem~\eqref{eq:discrete_ctrl_pb} is given by 
  \[
  \ctrls^* = \DynProg((\dyn_t)_{t=0}^{\horizon-1},(\cost_t)_{t=0}^\horizon , \initstate, \BellLQ),
  \]
  with $\DynProg$ as defined in~\eqref{eq:dyn_prog_algo} and $\BellLQ$
  implemented in Algo.~\ref{algo:BellLQ}
\end{coro}
\begin{proof} Note that at time $t \in \{0, \ldots, \horizon-1\}$ for a given
  $\state \in \reals^\dimstate$, if $\costogo_{t+1}$ is convex, then
  $\costogo_{t+1}(f_t(\state, \cdot))$ is convex as the composition of a
  convex function and a linear function and $\costogo_{t+1}(f_t(\state,
  \cdot)) + \cost_t(\state, \cdot)$ is then strongly convex as the sum of a
  convex and a strongly convex function. Moreover, $\state, \ctrl \rightarrow
  \costogo_{t+1}(f_t(\state, \ctrl)) + \cost_t(\state, \ctrl)$ is jointly
  convex since  $\state, \ctrl \rightarrow \costogo_{t+1}(f_t(\state, \ctrl))$
  is the composition of a convex function with a linear function and $\cost_t$
  is convex by assumption. Therefore, $\costogo_t: \state \rightarrow
  \min_{\ctrl \in \reals^\dimctrl}  \costogo_{t+1}(f_t(\state, \ctrl)) +
  \cost_t(\state, \ctrl)$ is convex as the partial infimum of jointly convex
  function.
	
  In summary, at time $t \in \{0, \ldots, \horizon-1\}$, if $\costogo_{t+1}$
is convex, then (i)  $\costogo_{t+1}(f_t(\state, \cdot)) + \cost_t(\state,
\cdot)$ is strongly convex, and (ii) $\costogo_{t}$ is convex. This ensures
that the assumptions of Lemma~\ref{lem:lin_quad_exact} are satisfied at each
iteration  of Algo.~\ref{algo:dyn_prog} (line~\ref{line:dyn_prog_backward})
since $\costogo_\horizon=\cost_\horizon$ is convex. \end{proof}

\subsection{Oracle Decomposition}
Before presenting the proof of Lemma~\ref{lem:lin_quad_oracle}, we present below a compact formulation of the first and second order information of $\traj$ with respect to the first and second order information of the dynamics $(\dyn_t)_{t=0}^{\horizon-1}$. The decomposition done in this lemma is reused for the proof of Lemma~\ref{lem:lin_quad_oracle}.
\begin{lemm}\label{lem:grad_hess_detailed}
	Consider the control  $\traj$ of $\horizon$ dynamics $(\dyn_t)_{t=0}^{\horizon-1}$ as defined in Def.~\ref{def:traj_func} and an initial point $\state_0 \in \reals^\dimstate$.
	For $\states = (\state_1;\ldots; \state_\horizon)$ and $\ctrls=(\ctrl_0; \ldots;\ctrl_{\horizon-1})$, define
	\[
	\concatdyn(\states, \ctrls) = (\dyn_0(\state_0, \ctrl_0); \ldots;\dyn_{\horizon-1}(\state_{\horizon-1}, \ctrl_{\horizon-1})).
	\]
	The gradient of  the control  $\traj$ of the dynamics $(\dyn_t)_{t=0}^{\horizon-1}$  on $\ctrls\in \reals^{\horizon\dimctrl}$ can be written
	\[
	\nabla_\ctrls \traj (\state_0, \ctrls) = \nabla_{\ctrls} \concatdyn(\states, \ctrls) (\idm - \nabla_\states \concatdyn(\states, \ctrls))^{-1}. 
	\]
	The Hessian of  the control  $\traj$ of the dynamics $(\dyn_t)_{t=0}^{\horizon-1}$  on $\ctrls\in  \reals^{\horizon\dimctrl}$ can be written
	\begin{align*}
		\nabla^2_{\ctrls\ctrls} \traj(\state_0,\ctrls)  & {=} 
		\nabla^2_{\states\states}\concatdyn(\states, \ctrls)[N , N,  M ] 
		{+} \nabla^2_{\ctrls\ctrls} \concatdyn(\states, \ctrls)[\cdot, \cdot, M ]   
		{+} \nabla^2_{\states\ctrls}\concatdyn(\states, \ctrls)[N,  \cdot, M ]
		{+} \nabla^2_{\ctrls\states}\concatdyn(\states, \ctrls)[\cdot, N,  M ],
	\end{align*}
	where $M =  (\idm - \nabla_\states \concatdyn(\states, \ctrls))^{-1}$ and $N=  \nabla_\ctrls \traj (\state_0, \ctrls) ^\top $.
\end{lemm}
\begin{proof}
	Denote simply, for $\ctrls \in \reals^{\horizon \dimctrl}$,  $\auxdyn(\ctrls) = \traj(\state_0, \ctrls)$ with $\state_0$ a fixed initial state. By definition, the function $\auxdyn$ can be decomposed, for $\ctrls\in \reals^{\horizon\dimctrl}$, as $\auxdyn(\ctrls) = (\auxdyn_1(\ctrls); \ldots;\auxdyn_\horizon(\ctrls))$, such that
	\begin{equation}\label{eq:detailed_traj_func}
		\auxdyn_{t+1}(\ctrls)  = \dyn_t(\auxdyn_{t}(\ctrls), E_t^\top \ctrls) \quad \mbox{for} \ t \in \{0, \ldots, \horizon-1\},
	\end{equation}
	with $\auxdyn_{0}(\ctrls) = \state_0$ and for $t \in \{0, \ldots, \horizon-1\}$,  $E_t = e_t \otimes \idm_{\dimctrl}$ is such that $E_t^\top \ctrls = \ctrl_t$, with $e_t$ the $t+1$\textsuperscript{th} canonical vector in $\reals^{\horizon}$, $\otimes$ the Kronecker product and $\idm_{\dimctrl} \in \reals^{\dimctrl\times\dimctrl}$ the identity matrix. By taking the derivative of \eqref{eq:detailed_traj_func}, we get,
	denoting $\state_t = \auxdyn_{t}(\ctrls)$ for $t\in \{0, \ldots, \horizon\}$  and
	using that $E_t^\top \ctrls= \ctrl_t$, 
	\begin{equation}\nonumber
		\nabla\auxdyn_{t+1}(\ctrls)  = \nabla \auxdyn_{t}(\ctrls) \nabla_{\state_t}\dyn_t(\state_t, \ctrl_t) + E_t \nabla_{\ctrl_t} \dyn_t(\state_t,\ctrl_t) \quad \mbox{for} \  t\in \{0, \ldots, \horizon-1\}.
	\end{equation}
	So, for $\auxctrls= (\auxctrl_0;\ldots;\auxctrl_{\horizon-1}) \in \reals^{\horizon\dimctrl}$, denoting $\nabla \auxdyn(\ctrls)^\top \auxctrls = (\auxstate_1;\ldots;\auxstate_\horizon)$ s.t. $\nabla \auxdyn_{t}(\ctrls)^\top \auxctrls = \auxstate_t$ for $t\in \{1, \ldots, \horizon\}$, 
	we have, with $\auxstate_0=0$,
	\begin{equation}\label{eq:detailed_grad_traj_func_simp}
		\auxstate_{t+1} = \nabla_{\state_t}\dyn_t(\state_t, \ctrl_t)^\top \auxstate_t + \nabla_{\ctrl_t} \dyn_t(\state_t,\ctrl_t)^\top \auxctrl_t \quad \mbox{for} \  t\in \{0, \ldots, \horizon-1\}.
	\end{equation}
	Denoting $\auxstates = (\auxstate_1;\ldots;\auxstate_{\horizon})$, we have then 
	\[
	(\idm-A)\auxstates = B\auxctrls, \quad \mbox{i.e.},  \quad \nabla \auxdyn(\ctrls)^\top\auxctrls = (\idm-A)^{-1} B\auxctrls,
	\]
	where $A = \sum_{t=1}^{\horizon-1} e_te_{t+1}^\top \otimes A_t$ with $A_t = \nabla_{\state_t} \dyn_t(\state_t, \ctrl_t) ^\top $ for $t \in \{1, \ldots, \horizon-1\}$ and $B = \sum_{t=1}^{\horizon} e_te_t^\top \otimes B_{t-1}$ with $B_t = \nabla_{\ctrl_t} \dyn_t(\state_t, \ctrl_t)^\top$ for $t\in \{0, \ldots, \horizon-1\}$, i.e.
	\begin{align*}
		A = \left(	\begin{matrix}
			0 & A_1 & 0 &\ldots & 0 \\
			\vdots & \ddots & \ddots & \ddots & \vdots \\
			& & & \ddots & 0 \\
			\vdots & & & \ddots & A_{\horizon-1} \\
			0 & \ldots & & \ldots & 0 
		\end{matrix}\right), \quad B = \left(\begin{matrix}
			B_0 & 0 & \ldots & 0 \\
			0 & \ddots & \ddots & \vdots \\
			\vdots & \ddots & \ddots & 0 \\
			0 & \ldots & 0 & B_{\horizon-1}
		\end{matrix}\right).
	\end{align*}
	By definition of  $\concatdyn$ in the claim,  one easily check that $A = \nabla_\states \concatdyn(\states, \ctrls)^\top$ and $B = \nabla_\ctrls \concatdyn(\states, \ctrls)^\top$. Therefore, we get 
	\[
	\nabla_{\ctrls} \traj(\state_0, \ctrls)	= \nabla \auxdyn(\ctrls) = \nabla_\ctrls \concatdyn(\states, \ctrls) (\idm- \nabla_\states \concatdyn(\states, \ctrls))^{-1}.
	\]	
	For the Hessian, note that for $g:\reals^d \rightarrow \reals^p$, $f:\reals^p \rightarrow \reals$, $x\in \reals^d$,  we have  
	$
	\nabla^2 (f\circ g)(x) = \nabla g(x)\nabla^2f(x)\nabla g(x)^\top + \nabla^2 g(x) [\cdot, \cdot, \nabla f(x)] \in \reals^{d\times d}.
	$
	If $f :\reals^p \rightarrow \reals^n$, we have 
	$
	\nabla^2 (f\circ g)(x) = \nabla^2 f(x)[\nabla g(x)^\top, \nabla g(x)^\top, \cdot]+ \nabla^2 g(x) [\cdot, \cdot, \nabla f(x)]\in \reals^{d\times d \times n}.
	$
	Applying this on $\dyn_t\circ g_t$ for $g_t(\ctrls) = (\auxdyn_{t}(\ctrls), E_t^\top\ctrls)$, we get from Eq.~\eqref{eq:detailed_traj_func}, using that $\nabla g_t(\ctrls)= (\nabla \auxdyn_{t}(\ctrls), E_t)$,
	\begin{align} \nonumber
		\nabla^2\auxdyn_{t+1}(\ctrls) & = 
		\nabla^2\auxdyn_{t}(\ctrls)[\cdot, \cdot, \nabla_{\state_t} \dyn_t(\state_t, \ctrl_t)]
		\\ 
		& 	\quad + \nabla^2_{\state_t\state_t} \dyn_t(\state_t, \ctrl_t)[\nabla \auxdyn_{t}(\ctrls)^\top, \nabla \auxdyn_{t}(\ctrls)^\top, \cdot]+ \nabla^2_{\ctrl_t\ctrl_t}\dyn_t(\state_t, \ctrl_t)[E_t^\top , E_t^\top, \cdot]  \nonumber\\
		&\quad  +	\nabla^2_{\state_t\ctrl_t}\dyn_t(\state_t, \ctrl_t)[\nabla \auxdyn_{t}(\ctrls)^\top, E_t^\top , \cdot] 
		+	\nabla^2_{\ctrl_t\state_t}\dyn_t(\state_t, \ctrl_t)[E_t^\top , \nabla \auxdyn_{t}(\ctrls)^\top, \cdot],  \nonumber
	\end{align}
	for $t\in \{0, \ldots, \horizon-1\}$, with $\nabla^2\auxdyn_0(\ctrls) = 0$. Therefore, for $\auxctrls= (\auxctrl_0;\ldots;\auxctrl_{\horizon-1}), \auxxctrls = (\auxxctrl_0;\ldots;\auxxctrl_{\horizon-1}) \in \reals^{\horizon\dimctrl}$, $\dualvars = (\dualvar_1;\ldots;\dualvar_{\horizon}) \in \reals^{\horizon \dimstate}$, we get
	\begin{align}\nonumber
		\nabla^2\auxdyn(\ctrls)[\auxctrls, \auxxctrls, \dualvars] &= \sum_{t=0}^{\horizon-1}\nabla^2\auxdyn_{t+1}(\ctrls)[\auxctrls, \auxxctrls, \dualvar_{t+1}]\\
		& = \sum_{t=0}^{\horizon-1} \Big( \nabla^2_{\state_t\state_t}\dyn_t(\state_t, \ctrl_t)[\auxstate_{t}, \auxxstate_{t}, \costate_{t+1}] 
		+ \nabla^2_{\ctrl_t\ctrl_t}\dyn_t(\state_t, \ctrl_t)[\auxctrl_t, \auxxctrl_t, \costate_{t+1} ]\label{eq:detailed_hessian_traj}\\
		& \hspace{35pt} + \nabla^2_{\state_t\ctrl_t}\dyn_t(\state_t, \ctrl_t)[\auxstate_t,\auxxctrl_t, \costate_{t+1} ] + \nabla^2_{\ctrl_t\state_t}\dyn_t(\state_t, \ctrl_t)[\auxctrl_t,\auxxstate_t, \costate_{t+1} ]
		\Big),  \nonumber
	\end{align}
	where $\auxstates = (\auxstate_1;\ldots;\auxstate_\horizon) = \nabla \auxdyn(\ctrls)^\top\auxctrls$,  $\auxxstates = (\auxxstate_1;\ldots;\auxxstate_\horizon) = \nabla \auxdyn(\ctrls)^\top\auxxctrls$, with $\auxstate_0 = \auxxstate_0=0$ and  $\costates = (\costate_1;\ldots;\costate_\horizon) \in \reals^{\horizon\dimstate}$ is defined by
	\begin{align*}
		\costate_{t} & =  \nabla_{\state_t} \dyn_{t}(\state_{t}, \ctrl_{t}) \costate_{t+1} + \dualvar_t\qquad \mbox{for} \ t\in\{1,\ldots, \horizon-1\}, \quad \costate_\horizon  =\dualvar_\horizon.
	\end{align*}
	On the other hand, denoting $\concatdyn_t(\states, \ctrls) = \dyn_t(\state_t, \ctrl_t)$ for $t\in \{0, \ldots, \horizon-1\}$, the Hessian of $\concatdyn$ with respect to the variables $\ctrls$ can be decomposed as
	\[
	\nabla^2_{\ctrls\ctrls}\concatdyn(\states,\ctrls)[\auxctrls, \auxxctrls, \costates] = \sum_{t=0}^{\horizon-1} \nabla^2_{\ctrls\ctrls} \concatdyn_t(\states, \ctrls)[\auxctrls, \auxxctrls, \costate_{t+1}] = \sum_{t=0}^{\horizon-1} \nabla^2_{\ctrl_t\ctrl_t} \dyn_t(\state_t, \ctrl_t)[\auxctrl_t, \auxxctrl_t, \costate_{t+1}].
	\]
	The Hessian of $\concatdyn$ with respect to the variable $\states$ can be decomposed as 
	\begin{align*}
		\nabla^2_{\states\states}\concatdyn(\states, \ctrls) [\auxstates, \auxxstates, \costates] 
    & = \sum_{t=0}^{\horizon-1} \nabla^2_{\states\states} \concatdyn_t(\states,\ctrls)[\auxstates, \auxxstates, \costate_{t+1}] 
    = \sum_{t=1}^{\horizon-1} \nabla^2_{\state_t\state_t} \dyn_t(\state_t,\ctrl_t)[\auxstate_t, \auxxstate_t, \costate_{t+1}].
	\end{align*}
  Finally, the second cross-derivatives of $\concatdyn$ w.r.t. $\states$ and $\ctrls$ can be decomposed as 
  \begin{align*}
    \nabla^2_{\states\ctrls}\concatdyn(\states, \ctrls)[\auxstates, \auxxctrls, \costates]
    & = \sum_{t=0}^{\horizon-1} \nabla^2_{\states\ctrls} \concatdyn_t(\states, \ctrls)[\auxstates, \auxxctrls, \costate_{t+1}]
    = \sum_{t=1}^{\horizon-1} \nabla^2_{\state_t\ctrl_t} \dyn_t(\state_t, \ctrl_t)[\auxstate_t, \auxxctrl_t, \costate_{t+1}].
  \end{align*}  
  From Eq.~\eqref{eq:detailed_hessian_traj}, we then get 
	\begin{align*}
		\nabla^2\auxdyn(\ctrls)[\auxctrls, \auxxctrls, \dualvars]  & 
		{=} \nabla^2_{\states\states} \concatdyn(\states, \ctrls)[\auxstates, \auxxstates, \costates] 
		{+}  \nabla^2_{\ctrls\ctrls}\concatdyn(\states,\ctrls)[\auxctrls, \auxxctrls, \costates]  
		{+} \nabla^2_{\states\ctrls}\concatdyn(\states,\ctrls)[\auxstates, \auxxctrls, \costates] 
		{+} \nabla^2_{\ctrls\states} \concatdyn(\states, \ctrls)[\auxctrls, \auxxstates, \costates].
	\end{align*}
	Finally, by noting that  $\auxstates = (\nabla_\ctrls \concatdyn(\states, \ctrls)(\idm -\nabla_\states \concatdyn(\states, \ctrls))^{-1} )^\top \auxctrls$, $\auxxstates=(\nabla_\ctrls \concatdyn(\states, \ctrls)(\idm -\nabla_\states \concatdyn(\states, \ctrls))^{-1} )^\top \auxxctrls$, and $\costates = (\idm -\nabla_\states \concatdyn(\states, \ctrls))^{-1} \dualvars$, the claim is shown.
\end{proof}

\begin{lemm} Consider a nonlinear dynamical problem summarized as 
	\[
	\min_{\ctrls\in \reals^{\horizon\dimctrl}} \cost\circ\augtraj(\ctrls), \quad 
		\mbox{where} \quad  \cost(\states, \ctrls) = \sum_{t=0}^{\horizon-1} \cost_t(\state_t, \ctrl_t) +\cost_\horizon(\state_\horizon), \quad \augtraj(\ctrls) = (\traj(\initstate, \ctrls), \ctrls),
	\]
	with $\traj$ the control of $\horizon$  dynamics $(\dyn_t)_{t=0}^{\horizon-1}$
	as defined in Def.~\ref{def:traj_func}.

	Let $\ctrls = (\ctrl_0; \ldots;\ctrl_{\horizon-1})$ and $\trajfunc(\initstate,
	\ctrls) = (\state_1;\ldots;\state_\horizon)$. 	Gradient~\eqref{eq:grad_step},
	Gauss-Newton~\eqref{eq:gn_step} and Newton~\eqref{eq:newton_step} oracles for
	$\cost\circ \augtrajfunc$ amount to solving for  $\diffctrls^*=(\diffctrl_0^*;
	\ldots;\diffctrl_{\horizon-1}^* )$ linear quadratic control problems of the
	form
	\begin{align}
		\min_{\substack{\diffctrl_0,\ldots, \diffctrl_{\horizon-1} \in \reals^{\dimctrl}\\\diffstate_0,\ldots, \diffstate_\horizon \in \reals^{\dimstate}}} \quad &  
		\sum_{t=0}^{\horizon-1} \qua_t(\diffstate_t, \diffctrl_t) + \qua_\horizon(\diffstate_\horizon)\nonumber\\
		\mbox{subject to} \quad & \diffstate_{t+1}=  \lin_{\dyn_t}^{\state_t, \ctrl_t} (\diffstate_t, \diffctrl_t)\quad \mbox{for} \  t \in \{0,\ldots,\horizon-1\}, \quad  \diffstate_0 = 0, \nonumber
	\end{align}
	where for
	\begin{enumerate}[nosep]
		\item[(i)]  the gradient oracle~\eqref{eq:grad_step},
		$\qua_\horizon(\diffstate_\horizon) =
		\lin_{\cost_\horizon}^{\state_\horizon}(\auxstate_\horizon) $ and,  for $0
		\leq t \leq  \horizon-1$, 
		\[
	\qua_t(\diffstate_t, \diffctrl_t) = \lin_{\cost_t}^{\state_t, \ctrl_t}(\diffstate_t, \diffctrl_t) + \frac{\reg}{2} \|\diffctrl_t\|_2^2,
		\] 
		\item[(ii)] the Gauss-Newton oracle~\eqref{eq:gn_step},
		$\qua_\horizon(\diffstate_\horizon) =
		\qua_{\cost_\horizon}^{\state_\horizon}(\auxstate_\horizon) $ and,  for $0
		\leq t\leq \horizon-1$, 
			\[
		\qua_t(\diffstate_t, \diffctrl_t)	= \qua_{\cost_t}^{\state_t, \ctrl_t}(\diffstate_t, \diffctrl_t) + \frac{\reg}{2} \|\diffctrl_t\|_2^2,
\]

		\item[(iii)] for the  Newton oracle~\eqref{eq:newton_step},
		$\qua_\horizon(\diffstate_\horizon) =
		\qua_{\cost_\horizon}^{\state_\horizon}(\auxstate_\horizon) $ and, defining 
		\begin{gather}\nonumber
			\costate_\horizon  =\nabla\cost_\horizon(\state_\horizon), \quad  
			\costate_{t}  = \nabla_{\state_t} \cost_t(\state_t, \ctrl_t) +   \nabla_{\state_t} \dyn_t(\state_t, \ctrl_t) \costate_{t+1} \quad \mbox{for} \ t\in\{\horizon-1,\ldots,1\},
		\end{gather}
		we have, for $0\leq t\leq \horizon-1$,
\[
			\qua_t(\diffstate_t, \diffctrl_t) = \qua_{\cost_t}^{\state_t, \ctrl_t} (\diffstate_t, \diffctrl_t) + \frac{1}{2} \nabla^2 \dyn_t(\state_t, \ctrl_t)[\cdot, \cdot, \lambda_{t+1}](\diffstate_t, \diffctrl_t)+ \frac{\reg}{2} \|\diffctrl_t\|_2^2,
\] 
where for $f:\reals^\dimstate \times \reals^\dimctrl \rightarrow
\reals^\dimstate$, $\state\in \reals^\dimstate$, $\ctrl \in \reals^\dimctrl$,
$\costate\in \reals^\dimstate$, we define
\begin{align}\nonumber
\nabla^2\dyn(\state, \ctrl)[\cdot, \cdot,  \costate]: (\diffstate, \diffctrl) \rightarrow & \nabla_{\state \state}^2\dyn(\state,\ctrl )[\diffstate, \diffstate, \costate ] + 2\nabla_{\state \ctrl}^2\dyn(\state,\ctrl )[\diffstate, \diffctrl, \costate]  {+} \nabla_{\ctrl \ctrl}^2\dyn(\state,\ctrl )[\diffctrl, \diffctrl, \costate ].
\end{align}
	\end{enumerate}
\end{lemm} 
\begin{proof} In the following, we denote for simplicity $\auxdyn(\ctrls) =
\traj(\initstate, \ctrls)$. The optimization oracles can be rewritten as
follows. \begin{enumerate} \item The gradient oracle~\eqref{eq:grad_step} is
given by \begin{align} \diffctrl^* & = \argmin_{\varaux \in \reals^\dimvar}
\bigg\{ \nabla \cost(\augtrajfunc(\fixedvar))^\top  \nabla
\augtrajfunc(\fixedvar)^\top\varaux + \frac{\reg}{2}\|\varaux\|_2^2
\bigg\}.\label{eq:grad_step_detailed} \end{align} \item The Gauss-Newton
oracle~\eqref{eq:gn_step} is given by \begin{align} \diffctrl^* & =
\argmin_{\varaux \in \reals^\dimvar} \bigg\{ \frac{1}{2}\varaux^\top \nabla
\augtrajfunc(\fixedvar) \nabla^2 \cost(\augtrajfunc(\ctrls)) \nabla
\augtrajfunc(\fixedvar)^\top \varaux + \nabla \cost(\augtrajfunc(\ctrls))^\top
\nabla \augtrajfunc(\fixedvar)^\top\varaux + \frac{\reg}{2}\|\varaux\|_2^2
\bigg\}.\label{eq:gn_step_detailed} \end{align} \item The Newton
oracle~\eqref{eq:newton_step} is given by \begin{align}
\hspace{-25pt}\diffctrl^* & {=} \argmin_{\varaux \in \reals^\dimvar}  \bigg\{
\frac{1}{2}\varaux^\top \nabla \augtrajfunc(\fixedvar) \nabla^2
\cost(\augtrajfunc(\ctrls)) \nabla \augtrajfunc(\fixedvar)^\top \varaux {+}
\frac{1}{2}\nabla^2 \augtrajfunc(\fixedvar)[ \varaux, \varaux,
\nabla\cost(\augtrajfunc(\ctrls))]   
    {+} \nabla \cost(\augtrajfunc(\ctrls))^\top  \nabla
    \augtrajfunc(\fixedvar)^\top\varaux {+} \frac{\reg}{2}\|\varaux\|_2^2
    \bigg\}. \label{eq:newton_step_detailed} \end{align} \end{enumerate}

We have, denoting $\states = \auxdyn(\ctrls)$, \begin{align*} \nabla
\cost(\augtrajfunc(\ctrls))^\top  \nabla \augtrajfunc(\fixedvar)^\top\varaux &  
  {=} \nabla_\states \cost(\states, \fixedvar)^\top \nabla
  \auxdyn(\fixedvar)^\top \auxctrls 
  + \nabla_\ctrls \cost(\states, \fixedvar)^\top \auxctrls \\
    \varaux^\top \nabla \augtrajfunc(\fixedvar) \nabla^2
    \cost(\augtrajfunc(\ctrls)) \nabla \augtrajfunc(\fixedvar)^\top \varaux  &
    {=}  \auxctrls^\top \nabla \auxdyn(\fixedvar) \nabla^2_{\states \states}
    \cost(\states, \fixedvar) \nabla \auxdyn(\fixedvar)^\top \auxctrls
    {+}\auxctrls^\top \nabla^2_{\ctrls \ctrls} \cost(\states,
    \fixedvar)\auxctrls {+} 2\auxctrls^\top \nabla \auxdyn(\fixedvar)
    \nabla^2_{\states \ctrls} \cost(\states, \fixedvar)\auxctrls \\
    \nabla^2 \augtrajfunc(\fixedvar)[ \varaux, \varaux,
    \nabla\cost(\augtrajfunc(\ctrls))] & {=} \nabla^2
    \auxdyn(\fixedvar)[\auxctrls, \auxctrls, \nabla_\states \cost(\states,
    \fixedvar)]. \end{align*}

For $\auxctrls= (\auxctrl_0;\ldots;\auxctrl_{\horizon-1}) \in
\reals^{\horizon\dimctrl}$, denoting $\auxstates = \nabla \auxdyn(\ctrls)^\top
\auxctrls = (\auxstate_1;\ldots;\auxstate_\horizon)$, with $\auxstate_0 = 0$,
we have  then \begin{equation}\label{eq:decomp_lin_costs} \nabla
\cost(\augtrajfunc(\ctrls))^\top  \nabla \augtrajfunc(\fixedvar)^\top\varaux
{=} \sum_{t=0}^{\horizon-1} \left[\nabla_{\state_t}\cost_t(\state_t,
\ctrl_t)^\top \auxstate_t {+} \nabla_{\ctrl_t}\cost_t(\state_t, \ctrl_t)^\top
\auxctrl_t\right] {+} \nabla \cost_\horizon(\state_\horizon)^\top
\auxstate_\horizon {=} \sum_{t=0}^{\horizon-1} \lin_{\cost_t}^{\state_t,
\ctrl_t}(\auxstate_t, \auxctrl_t) {+}
\lin_{\cost_\horizon}^{\state_\horizon}(\auxstate_\horizon). \end{equation}
Following the proof of Lemma~\ref{lem:grad_hess_detailed}, we have that
$\auxstates = \nabla \auxdyn(\ctrls)^\top \auxctrls =
(\auxstate_1;\ldots;\auxstate_\horizon)$ satisfies
\begin{equation}\label{eq:linear_chain_grad} \auxstate_{t+1} =
\nabla_{\state_t}\dyn_t(\state_t, \ctrl_t)^\top \auxstate_t + \nabla_{\ctrl_t}
\dyn_t(\state_t,\ctrl_t)^\top \auxctrl_t = \lin_{\dyn_t}^{\state_t, \ctrl_t}
(\auxstate_t, \auxctrl_t), \quad \mbox{for} \ t \in \{0, \ldots, \horizon-1\},
\end{equation} with $\auxstate_0=0$. Hence, plugging
Eq.~\eqref{eq:decomp_lin_costs} and Eq.~\eqref{eq:linear_chain_grad} into
Eq.~\eqref{eq:grad_step_detailed} we get the claim for the gradient oracle. 

The Hessians of the total cost are block diagonal with, e.g.,
$\nabla_{\ctrls\ctrls}^2 \cost(\states, \ctrls) $ being composed of $\horizon$
diagonal blocks of the form $\nabla_{\ctrl_t\ctrl_t}^2 \cost_t(\state_t,
\ctrl_t)$ for $t\in \{0, \ldots, \horizon-1\}$. Therefore, we have
\begin{align} & \frac{1}{2}\varaux^\top \nabla \augtrajfunc(\fixedvar)
\nabla^2 \cost(\augtrajfunc(\ctrls)) \nabla \augtrajfunc(\fixedvar)^\top
\varaux \nonumber\\
  &  = \sum_{t=0}^{\horizon-1} \left[ \frac{1}{2} \auxstate_t^\top
  \nabla_{\state_t\state_t}^2\cost_t(\state_t, \ctrl_t) \auxstate_t {+}
  \frac{1}{2}\auxctrl_t^\top \nabla_{\ctrl_t\ctrl_t}^2\cost_t(\state_t,
  \ctrl_t) \auxctrl_t {+} \auxstate_t^\top
  \nabla_{\state_t\ctrl_t}^2\cost_t(\state_t, \ctrl_t) \auxctrl_t    \right]
{+}  
\frac{1}{2}\auxstate_\horizon^\top \nabla^2 \cost_\horizon(\state_\horizon)
\auxstate_\horizon. \nonumber \end{align} The linear quadratic approximation
in~\eqref{eq:gn_step_detailed} can then be written as \begin{align}
\frac{1}{2}\varaux^\top \nabla \augtrajfunc(\fixedvar) \nabla^2
\cost(\augtrajfunc(\ctrls)) \nabla \augtrajfunc(\fixedvar)^\top \varaux 
   + \nabla \cost(\augtrajfunc(\ctrls))^\top  \nabla
     \augtrajfunc(\fixedvar)^\top\varaux  = \sum_{t=0}^{\horizon-1}
     \qua_{\cost_t}^{\state_t, \ctrl_t} (\auxstate_t, \auxctrl_t) +
     \qua_{\cost_\horizon}^{\state_\horizon}(\auxstate_\horizon).
     \label{eq:decomp_quad_costs} \end{align} Hence, plugging
     Eq.~\eqref{eq:decomp_quad_costs} and Eq.~\eqref{eq:linear_chain_grad}
     into Eq.~\eqref{eq:gn_step_detailed} we get the claim for the
     Gauss-Newton oracle.

For the Newton oracle, denoting $\dualvars {= }
\nabla_\states\cost(\states,\var) {=} ( \nabla_{\state_1} \cost_1(\state_1,
\ctrl_1);\ldots; \nabla_{\state_{\horizon-1}}
\cost_{\horizon-1}(\state_{\horizon-1}, \ctrl_{\horizon-1}); \nabla
\cost_\horizon(\state_\horizon))$,  and defining adjoint variables
$\costate_t$ as \begin{align*} \costate_\horizon &
=\nabla\cost_\horizon(\state_\horizon) \qquad \costate_{t} =
\nabla_{\state_t} \cost_t(\state_t, \ctrl_t) + \nabla_{\state_t}
\dyn_{t}(\state_{t}, \ctrl_{t}) \costate_{t+1} \qquad \mbox{for} \
t\in\{1,\ldots, \horizon-1\}, \end{align*} we have, as in the proof of
Lemma~\ref{lem:grad_hess_detailed}, \begin{align} \nabla^2 \auxdyn(\var)
[\varaux, \varaux, \nabla_\states\cost(\states,\var)] & =
\sum_{t=0}^{\horizon-1}  \nabla^2 \auxdyn_{t+1}(\ctrls)[\auxctrls, \auxctrls,
\dualvar_{t+1}]  \nonumber\\
  &  = \sum_{t=0}^{\horizon-1} \Big(
  \nabla^2_{\state_t\state_t}\dyn_t(\state_t, \ctrl_t)[\auxstate_{t},
  \auxstate_{t}, \costate_{t+1}] 
  + \nabla^2_{\ctrl_t\ctrl_t}\dyn_t(\state_t, \ctrl_t)[\auxctrl_t, \auxctrl_t,
    \costate_{t+1} ] \nonumber   \\    & \hspace{34pt}
  + 2 \nabla^2_{\state_t\ctrl_t}\dyn_t(\state_t,
    \ctrl_t)[\auxstate_{t},\auxctrl_t, \costate_{t+1} ] \Big).
    \label{eq:newton_decomp} \end{align} Hence, plugging
    Eq.~\eqref{eq:decomp_quad_costs}, Eq.~\eqref{eq:newton_decomp} and
    Eq.~\eqref{eq:linear_chain_grad} into Eq.~\eqref{eq:newton_step_detailed}
    we get the claim for the Newton oracle. \end{proof}

  \section{Line-search}\label{app:line_searches}
	So far, we defined procedures that, given a command and some regularization
parameter, output a direction that minimizes an approximation of the objective
or approximately minimizes a shifted objective. Given access to such procedures,
the next command can be computed in several ways. The main criterion is to
ensure that the value of the objective decreases along the iterations, which is
generally done by a line-search. 

In the following, we only consider oracles based on linear quadratic or
quadratic approximations of the objective such as Gauss-Newton and Newton, and
refer the reader to~\citet{nocedal2006numerical} for classical line-searches for
gradient descent. 

\subsection{Rule}
We start by considering the implementation of line-searches for classical
optimization oracles which can again exploit the dynamical structure of the
problem and are mimicked by differential dynamic programming approaches. We
consider, as in Sec.~\ref{sec:classical_optim}, that we have access to an oracle
for an objective $\obj$, that, given a command $\ctrls \in \reals^{\horizon
\dimctrl}$ and any regularization $\reg \geq 0$, outputs
\begin{equation}\label{eq:generic_oracle}
\Oracle_\reg(\obj)(\var) = \argmin_{\diffctrls \in \reals^{\horizon\dimctrl}} \model_{\obj}^\ctrls(\diffctrls) + \frac{\reg}{2} \|\diffctrls\|_2^2,
\end{equation}
where $\model_\obj^\ctrls$ is a linear quadratic or quadratic expansion of the
objective $\obj$ around $\ctrls$ s.t. $\obj(\ctrls +\diffctrls) \approx
\obj(\ctrls) + \model_\obj^\ctrls(\diffctrls)$. Given such an oracle, we can
define a new candidate command that decreases the value of the objective in
several ways.\\

{\bf Directional Step.}
The next iterate can be defined along the direction provided by the oracle, as
long as this direction is a descent direction. Namely, the next iterate can be
computed as 
\begin{equation}\label{eq:descent_dir}
\ctrls^{\nxt} = \ctrls + \stepsize \diffctrls, \quad \mbox{with} \ \diffctrls = \Oracle_\reg(\obj)(\var)\ \mbox{for}\ \reg\geq 0 \ \mbox{s.t.} \ \nabla \obj(\ctrls)^\top \diffctrls<0,
\end{equation}
where the stepsize $\stepsize$ is chosen to satisfy, e.g., an Armijo
condition~\citep[Chapter 3]{nocedal2006numerical}, that is,
\begin{equation}\label{eq:linesearch_armijo}
\obj(\ctrls + \stepsize\diffctrls) \leq \obj(\ctrls) + \frac{\stepsize}{2}\nabla \obj(\ctrls)^\top \diffctrls.
\end{equation}
In this case, the search is usually initialized at each step with $\stepsize=1$.
If condition~\eqref{eq:linesearch_armijo} is not satisfied for $\stepsize=1$,
the stepsize is decreased by a factor $\rho_{\dec}<1$ until
condition~\eqref{eq:linesearch_armijo} is satisfied. If a stepsize $\stepsize=1$
is accepted, then the linear quadratic or quadratic algorithms may exhibit a
quadratic local convergence~\citep[Chapter 3, 10]{nocedal2006numerical}.
Alternative line-search criterions such as Wolfe's curvature condition or
trust-region methods can also be implemented~\citep[Chapter
3]{nocedal2006numerical}.\\

{\bf Regularized Step.}
Given a current iterate $\ctrls \in \reals^{\horizon\dimctrl}$,  we can find a
regularization such that the current iterate plus the direction output by the
oracle decreases the objective.  Namely, the next command can  be computed as
\begin{equation}\label{eq:reg_step}
	\ctrls^{\nxt}   = \ctrls + \diffctrls^\stepsize, \quad \mbox{where} \quad \diffctrls^\stepsize = \Oracle_{1/\stepsize}(\obj)(\ctrls) =  \argmin_{\diffctrls\in \reals^{\horizon \dimctrl}} \model_\obj^\ctrls( \diffctrls) + \frac{1}{2\stepsize}\|\diffctrls\|_2^2 , \qquad 
\end{equation}
where the parameter $\stepsize>0 $ acts as a stepsize that controls how large
should be the step (the smaller the parameter $\stepsize$, the smaller the step
$\diffctrls^\stepsize$). The stepsize $\stepsize$ can then be chosen to satisfy
\begin{equation}\label{eq:linesearch_reg}
	\obj(\ctrls+\diffctrls^\stepsize) \leq \obj(\ctrls) + \model_\obj^\ctrls( \diffctrls^\stepsize) + \frac{1}{2\stepsize} \|\diffctrls^\stepsize\|_2^2,
\end{equation}
which ensures a sufficient decrease of the objective to, e.g., prove convergence
to stationary points~\citep{roulet2019iterative}. In practice, as for the
line-search on the descent direction, given an initial stepsize for the
iteration, the stepsize is either selected or reduced by a factor $\rho_{\dec}$
until condition~\eqref{eq:linesearch_reg} is satisfied. However, here, we
initialize the stepsize at each iteration as $\rho_{\inc}\stepsize_{prev}$ where
$\stepsize_{prev}$ is the stepsize selected at the previous iteration and
$\rho_{\inc}>1$ is an increasing factor. By trying a larger stepsize at each
iteration, we may benefit from larger steps in some  regions of the optimization
path. Note that such an approach is akin to trust region methods which increase
the radius of the trust region at each iteration depending on the success of
each iteration~\citep{nocedal2006numerical}.

In practice, we observed that, when using regularized steps, acceptable
stepsizes for condition~\eqref{eq:linesearch_reg} tend to be arbitrarily large
as the iterations increase. Namely, we tried choosing $\rho_{\inc}=10$ and
observed that the acceptable stepsizes tended to plus infinity with such a
procedure. To better capture this tendency, we consider regularizations that may
depend on the current state and of the form $\reg(\ctrls) \propto \|\nabla
\cost(\states, \ctrls)\|_2$, i.e., stepsizes of the form $\stepsize(\ctrls)=
\stepsizescaled/\|\nabla \cost(\states, \ctrls)\|_2$. The line-search is then
performed on $\stepsizescaled$ only. Intuitively, as we are getting closer to a
stationary point, quadratic models are getting more accurate to describe the
objective. By scaling the regularization with respect to $\|\nabla
\cost(\states, \ctrls)\|_2$, which is a measure of stationarity, we may better
capture such behavior. Note that for $\reg=0$, we retrieve the iteration with a
descent direction of stepsize $\stepsize=1$ described above.

\subsection{Implementation}~
{\bf Directional Step.}
The Armijo condition~\eqref{eq:linesearch_armijo} can  be computed directly from
the knowledge of a gradient oracle and the chosen oracle (such as Gauss-Newton
or Newton). We present here the implementation of the line-search in terms of
the dynamical structure of the problem. Denote
\[
(\pi_t)_{t=0}^{\horizon-1}, \costogo_0 = \Backward((\model_{\dyn_t}^{\state_t, \ctrl_t})_{t=0}^{\horizon-1}, (\model_{\cost_t}^{\state_t, \ctrl_t})_{t=0}^{\horizon-1}, \model_{\cost_\horizon}^{\state_\horizon}, \reg)
\]
the policies and the value of the cost-to-go function output by the backward
pass of the considered oracle, i.e.,  Gauss-Newton or Newton.

By definition, $\costogo_0(0)$ is the minimum of the corresponding linear
quadratic control problem~\eqref{eq:lin_quad_oracle}. Moreover, the linear
quadratic control problem can be summarized as a quadratic problem of the form
$\min_{\diffctrls} \model_\obj(\diffctrls) + \frac{\reg}{2}\|\diffctrls\|_2^2 =
\min_{\diffctrls}\frac{1}{2}\diffctrls^\top (Q+\reg\idm)\diffctrls +  \nabla
\obj(\ctrls)^\top \diffctrls$ with $Q $ a quadratic that is either the Hessian
of $\obj$ for a Newton oracle or an approximation of it for a Gauss-Newton
oracle. Therefore, we have that, for a Newton or a Gauss-Newton oracle
$\diffctrls = \Oracle_{1/\stepsize}(\obj)$,
\[
\frac{1}{2}\nabla \obj(\ctrls)^\top \diffctrls = - \frac{1}{2} \nabla \obj(\ctrls)^\top(Q+\reg \idm)^{-1}\nabla \obj(\ctrls)= \min_{\diffctrls \in \reals^{\horizon \dimctrl}} \model_\obj(\diffctrls) + \frac{\reg}{2}\|\diffctrls\|_2^2 =  \costogo_0(0).
\]
Therefore, the right-hand part of condition~\eqref{eq:linesearch_armijo} can be
given by the value of the cost-to-go function $\costogo_0(0)$. On the other
hand,  sequences of controllers of the form $\gamma \diffctrls$ can be defined
by modifying the policies output in the backward pass as shown in the following
lemma adapted from~\citet[Theorem 1]{liao1992advantages}. 
\begin{lemm}
	Given a sequence of affine  policies $(\pi_t)_{t=0}^{\horizon-1}$, linear
	dynamics $(\lin_t)_{t=0}^{\horizon-1}$ and an initial state $\diffstate_0=0$,
	denote $\diffctrls^* = \Roll(\diffstate_0, (\pi_t)_{t=0}^{\horizon-1},
	(\lin_t)_{t=0}^{\horizon-1}) $  and $\pi_t^{\stepsize}: \diffstate \rightarrow
	\stepsize\pi_t(0) + \nabla \pi_t(0)^\top \diffstate$ for $t = 0, \ldots,
	\horizon-1$. We have that
	\begin{align*}
		\gamma \diffctrls^* &= \diffctrls^\stepsize, 
	\quad 	\mbox{where} \quad  \diffctrls^\stepsize =  \Roll(\diffstate_0,  (\pi_t^{\stepsize})_{t=0}^{\horizon-1}, (\lin_t)_{t=0}^{\horizon-1}).	
\end{align*}
\end{lemm}
\begin{proof}
	Define $(\diffstate_t^\stepsize)_{t=0}^{\horizon-1}$ as
	$\diffstate_{t+1}^\stepsize = \lin_t(\diffstate_t^\stepsize,
	\pi_t(\diffstate_t^\stepsize))$ for $t\in\{0, \ldots, \horizon-1\}$ with
	$\diffstate^\stepsize_0 = 0$. We have that $\diffstate^\stepsize_1$ is linear
	w.r.t. $\stepsize$. Proceeding by induction, we have that
	$\diffstate_t^\stepsize$ is linear w.r.t. $\stepsize$ using the form of
	$\pi_t^\stepsize$ and the fact that $\lin_t$ is linear. Therefore,
	$\diffctrl_t^\stepsize = \pi_t^\stepsize(\diffstate_t^\stepsize)$ is linear
	w.r.t. $\stepsize$ which gives the claim.   
\end{proof}

Therefore, computing the next sequence of controllers by moving along a descent
direction as in~\eqref{eq:descent_dir} according to an Armijo
condition~\eqref{eq:linesearch_armijo} amounts to computing, with
Algo.~\ref{algo:line_search},
\begin{align*}
	\ctrls^{\nxt} & = \linesearch(\ctrls, (\cost_t)_{t=0}^\horizon,  (\dyn_t)_{t=0}^{\horizon-1}, (\lin_{\dyn_t}^{\state_t, \ctrl_t})_{t=0}^{\horizon-1}, \Pol),\\
	\mbox{where}\quad &  \Pol: \stepsize \rightarrow \left(\begin{array}{ll}
		(\pi_t^{\stepsize}: &\diffstate \rightarrow \stepsize\pi_t(0) + \nabla \pi_t(0)^\top \diffstate)_{t=0}^{\horizon-1} \\
		\hspace{6pt} \costogo_0^\stepsize: &\diffstate \rightarrow \stepsize \costogo_0(\diffstate)
	\end{array}\right) \\
	& 	(\pi_t)_{t=0}^{\horizon-1}, \costogo_0 = 
	\Backward((\model_{\dyn_t}^{\state_t, \ctrl_t})_{t=0}^{\horizon-1}, (\model_{\cost_t}^{\state_t, \ctrl_t})_{t=0}^{\horizon-1}, \model_{\cost_\horizon}^{\state_\horizon}, \reg), \  \mbox{for} \ \reg\geq 0 \ \mbox{s.t.}\  \costogo_0(0)<0,
\end{align*}
where $\Backward \in \{\Backward_{\operatorname{GN}},
\Backward_{\operatorname{NE}}\}$ is given in Algo.~\ref{algo:backward_gn} or
Algo.~\ref{algo:backward_newton}. 

In practice, in our implementation of the backward passes in
Algo.~\ref{algo:backward_gn}, Algo.~\ref{algo:backward_newton}, the returned
initial cost-to-go function is either negative if the step is well-defined or
infinite if it is not. To find a regularization that ensures a descent
direction, i.e., $\costogo_0(0)<0$, it suffices thus to find a feasible step. In
our implementation, we first try to compute a descent direction without
regularization ($\reg=0$), then try a small regularization $\reg=10^{-6}$, which
we increase by 10 until a finite negative cost-to-go function $c_0(0)$ is
returned. See Algo.~\ref{algo:gn_algo} for an instance of such implementation.

From the above discussion, it is clear that one iteration of the Iterative
Linear Quadratic Regulator algorithm described in Sec.~\ref{ssec:nonlin_algo}
uses a Gauss-Newton oracle without regularization to move along the direction of
the oracle by using an Armijo condition. The overall iteration is given in
Algo.~\ref{algo:gn_algo}, where we added a procedure to ensure that the output
direction is a descent direction. All other algorithms, with or without
regularization can be written in a similar way using a forward, a backward pass,
and multiple roll-out phases until the next sequence of controllers is found.\\

{\bf Regularized Step.}
For regularized steps, the line-search~\eqref{eq:linesearch_reg} requires
computing $\model_\obj^\ctrls( \diffctrl^\stepsize) + \frac{1}{2\stepsize}
\|\diffctrl^\stepsize\|_2^2$. This is by definition the minimum of the
sub-problem that is computed by dynamic programming. This minimum can therefore
be accessed as $\model_\obj^\ctrls( \diffctrls^\stepsize) + \frac{1}{2\stepsize}
\|\diffctrls^\stepsize\|_2^2 = \costogo_0(0)$ for $\costogo_0$  output by the
backward pass with a regularization $\reg =1/\stepsize$. Overall, the next
sequence of controls is then provided through the line-search procedure given in
Algo.~\ref{algo:line_search} as 
\begin{align*}
	\ctrls^{\nxt} & = \linesearch(\ctrls, (\cost_t)_{t=0}^\horizon,  (\dyn_t)_{t=0}^{\horizon-1}, (\lin_{\dyn_t}^{\state_t, \ctrl_t})_{t=0}^{\horizon-1}, \Pol),\\
	\mbox{where}\quad &  \Pol: \stepsize \rightarrow \Backward((\model_{\dyn_t}^{\state_t, \ctrl_t})_{t=0}^{\horizon-1}, (\model_{\cost_t}^{\state_t, \ctrl_t})_{t=0}^{\horizon-1}, \model_{\cost_\horizon}^{\state_\horizon}, 1/\stepsize),
\end{align*}
where $\Backward \in \{\Backward_{\operatorname{GN}},
\Backward_{\operatorname{NE}}\}$ is given in Algo.~\ref{algo:backward_gn} or
Algo.~\ref{algo:backward_newton}.\\

{\bf Line-search for Differential Dynamic Programming Approaches.}
The line-search for DDP approaches as presented by, e.g., \citet[Sec.
2.2]{liao1992advantages} based on~\citet{jacobson1970differential}, mimics the
one done for the classical optimization oracles except that the policies are
rolled out on the original dynamics. Namely, the usual line-search consists in
applying Algo.~\ref{algo:line_search} as follows
\begin{align*}
	\ctrls^{\nxt} & = \linesearch(\ctrls, (\cost_t)_{t=0}^\horizon,  (\dyn_t)_{t=0}^{\horizon-1}, (\diff_{\dyn_t}^{\state_t, \ctrl_t})_{t=0}^{\horizon-1}, \Pol)\\
	\mbox{where}\quad &  \Pol: \stepsize \rightarrow \left(\begin{array}{ll}
		(\pi_t^{\stepsize}: &\diffstate \rightarrow \stepsize\pi_t(0) + \nabla \pi_t(0)^\top \diffstate)_{t=0}^{\horizon-1}, \\
		\hspace{6pt} c_0^\stepsize: &\diffstate \rightarrow \stepsize c_0(\diffstate)
	\end{array}\right) \\
	& 	(\pi_t)_{t=0}^{\horizon-1}, c_0 = 
	\Backward((\model_{\dyn_t}^{\state_t, \ctrl_t})_{t=0}^{\horizon-1}, (\model_{\cost_t}^{\state_t, \ctrl_t})_{t=0}^{\horizon-1}, \model_{\cost_\horizon}^{\state_\horizon}, \reg)  \  \mbox{for} \ \reg\geq 0 \ \mbox{s.t.}\  \costogo_0(0)<0,
\end{align*}
where $\Backward \in \{\Backward_{\operatorname{GN}},
\Backward_{\operatorname{DDP}}\}$ is given by Algo.~\ref{algo:backward_gn} or
Algo.~\ref{algo:backward_ddp}. As for the classical optimization oracles, a
direction is first computed without regularization and if the resulting
direction is not a descent direction a small regularization is added to ensure
that $\costogo_0(0)<0$.

We also consider line-searches based on selecting an appropriate regularization.
Namely, we consider line-searches of the form 
\begin{align*}
	\ctrls^{\nxt} & = \linesearch(\ctrls, (\cost_t)_{t=0}^\horizon,  (\dyn_t)_{t=0}^{\horizon-1}, (\diff_{\dyn_t}^{\state_t, \ctrl_t})_{t=0}^{\horizon-1}, \Pol),\\
	\mbox{where}\quad &  \Pol: \stepsize \rightarrow \Backward((\model_{\dyn_t}^{\state_t, \ctrl_t})_{t=0}^{\horizon-1}, (\model_{\cost_t}^{\state_t, \ctrl_t})_{t=0}^{\horizon-1}, \model_{\cost_\horizon}^{\state_\horizon}, 1/\stepsize),
\end{align*}
where $\Backward \in \{\Backward_{\operatorname{GN}},
\Backward_{\operatorname{DDP}}\}$ is given by Algo.~\ref{algo:backward_gn} or
Algo.~\ref{algo:backward_ddp}.

  \section{Detailed Computational Scheme}\label{app:summary}
\begin{figure}
	\begin{center}
		\includegraphics[width=0.9\linewidth]{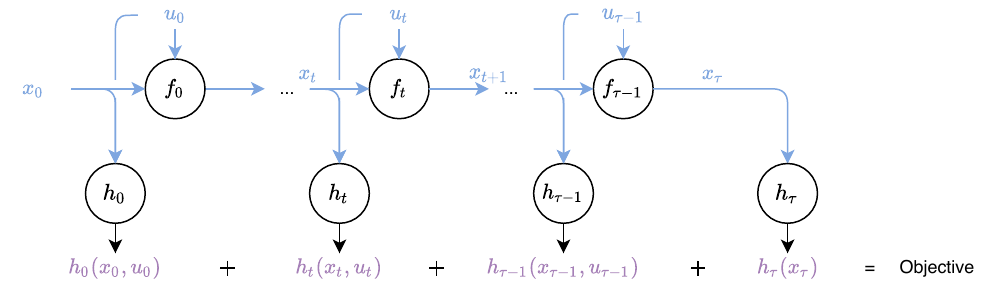}
	\end{center}
	\caption{Computational scheme of the discrete time control
	problem~\eqref{eq:discrete_ctrl_pb}. \label{fig:graph_comput_traj}}
\end{figure}

	We detail here the algorithms presented in Figure~\ref{fig:taxonomy}. Recall
that our objective is 
\begin{align*}
	\obj(\ctrls) = &  \sum_{t=0}^{\horizon-1} \cost_t(\state_t, \ctrl_t) +\cost_\horizon(\state_\horizon) \\
	& \mbox{s.t.} \quad \state_{t+1} = f_t(\state_t, \ctrl_t) \quad \mbox{for} \ t \in \{0, \ldots, \horizon-1\}, \quad \state_0 = \initstate,
\end{align*}
that can be summarized as $\obj(\ctrls) = \cost(\augtraj(\ctrls))$, where, for
$\ctrls = (\ctrl_0;\ldots;\ctrl_{\horizon-1})$, $\states
=(\state_1;\ldots;\state_\horizon)$,
\begin{align}\nonumber
 \cost(\states, \ctrls) = \sum_{t=0}^{\horizon-1} \cost_t(\state_t, \ctrl_t) +\cost_\horizon(\state_\horizon), \ \augtraj(\ctrls) = (\traj(\initstate, \ctrls), \ctrls),  \
	\traj(\state_0, \ctrls)  &= ( \state_1;\ldots;\state_\horizon) \nonumber\\[-10pt]
	\mbox{s.t.} \quad \state_{t+1} &= \dyn_t(\state_t, \ctrl_t) \quad \mbox{for} \ t \in \{0,\ldots, \horizon-1\}. \nonumber
\end{align}
The computational graph of the objective is illustrated in
Figure~\ref{fig:graph_comput_traj}.

We present nonlinear control algorithms from a functional viewpoint by
introducing finite difference, linear and quadratic expansions of the dynamics
and the costs presented in the notations in Eq.~\eqref{eq:lin_approx}.

For a function $f:\reals^\dimstate \times \reals^\dimctrl \rightarrow\reals^p$,
with $p=1$ (for the costs) or $p=\dimstate$ (for the dynamics), these expansions
read for $\state, \ctrl \in \reals^\dimstate \times \reals^\dimctrl$, 
\begin{gather}
		\diff_f^{\state, \ctrl}: \auxstate, \auxctrl \rightarrow 
		f(\state+ \auxstate, \ctrl +\auxctrl) - f(\state, \ctrl), \qquad 
		\lin_f^{\state, \ctrl}: \auxstate, \auxctrl \rightarrow
		\nabla_\state f(\state, \ctrl)^\top \auxstate + \nabla_\ctrl f(\state, \ctrl)^\top \auxctrl \label{eq:reminder_approx}\\ 
		\qua_f^{\state, \ctrl}: \auxstate, \auxctrl \rightarrow
		\nabla_\state f(\state, \ctrl)^\top \auxstate + \nabla_\ctrl f(\state, \ctrl)^\top \auxctrl 
		+ \frac{1}{2}  \nabla_{\state \state}^2\dyn(\state,\ctrl )[\diffstate, \diffstate, \cdot ] 
		+\frac{1}{2}\nabla_{\ctrl \ctrl}^2\dyn(\state,\ctrl )[\diffctrl, \diffctrl, \cdot ] 
		+ \nabla_{\state \ctrl}^2\dyn(\state,\ctrl )[\diffstate, \diffctrl, \cdot]   \nonumber
\end{gather}
For $\costate \in \reals^p$, we denote shortly
\[
\frac{1}{2}\nabla^2\dyn(\state, \ctrl)[\cdot, \cdot,  \costate]: (\diffstate, \diffctrl) \rightarrow  \frac{1}{2}\nabla_{\state \state}^2\dyn(\state,\ctrl )[\diffstate, \diffstate, \costate ] +\frac{1}{2} \nabla_{\ctrl \ctrl}^2\dyn(\state,\ctrl )[\diffctrl, \diffctrl, \costate ] + \nabla_{\state \ctrl}^2\dyn(\state,\ctrl )[\diffstate, \diffctrl, \costate].
\]

In the algorithms, we consider storing in memory linear or quadratic functions
by storing the associated vectors, matrices or tensors defining the linear or
quadratic functions. For example, to store the linear expansion $\lin_f^x$ or
the quadratic expansion $q_f^x$ of a function $f:\reals^d\rightarrow\reals^p$
around a point $x$, we consider storing $\nabla f(x) \in \reals^{d\times p}$ and
$\nabla^2 f(x) \in \reals^{d\times d \times p}$. In the backward or roll-out
passes, we consider that having access to the linear or quadratic functions,
means having  access to the associated matrices/tensors defining the operations
as presented in, e.g., Algo.~\ref{algo:BellLQ}. The functional viewpoint helps
to isolate the main technical operations in the procedures $\BellLQ$ in
Algo.~\ref{algo:BellLQ} or $\BellL$  in Algo.~\ref{algo:bp} and to identify the
discrepancies between, e.g., the Newton oracle in Algo.~\ref{algo:newton} and a
DDP oracle with quadratic approximations presented in Algo.~\ref{algo:ddp_q}.
For a presentation of the algorithms in a purely algebraic viewpoint, we refer
the reader to, e.g.,~\citet{wright1990solution, liao1992advantages,
sideris2005efficient}.

In
Algo.~\ref{algo:backward_gn},~\ref{algo:backward_newton},~\ref{algo:backward_ddp},
we a priori need to check whether the subproblems defined by the Bellman
recursion are strongly convex or not. Namely, in
Algo.~\ref{algo:backward_gn},~\ref{algo:backward_newton},~\ref{algo:backward_ddp},
we need to check that $\qua_t(\state, \cdot) + \costogo_{t+1}(\lin_t(\state,
\cdot))$ is strongly convex for any $x$. With the notations of
Algo.~\ref{algo:BellLQ}, this amounts checking that $\Q + \B^\top \J_{t+1}\B
\succ 0$. This can be done by checking the positivity of the minimum eigenvalue
of $\Q + \B^\top \J_{t+1}\B $.  In our implementation, we simply check that 
\begin{equation}\label{eq:weak_feasibility_condition}
	\jcst_t - \jcst_{t+1} = -\frac{1}{2}(\q + \B^\top \j_{t+1})^\top (\Q + \B^\top \J_{t+1}\B)^{-1}(\q + \B^\top \j_{t+1}) <0.
\end{equation}
If condition~\eqref{eq:weak_feasibility_condition} is not satisfied then
necessarily $\Q + \B^\top \J_{t+1}\B \not \succeq 0$. We chose to use
condition~\eqref{eq:weak_feasibility_condition} since this quantity is directly
available and computing the eigenvalues of $\Q + \B^\top \J_{t+1}\B \succ 0$ can
slow down the computations. Moreover, if
criterion~\eqref{eq:weak_feasibility_condition} is satisfied for all $t\in \{0,
\ldots, \horizon-1\}$, this means that, for the Gauss-Newton and the Newton
methods, the resulting direction is a descent direction for the objective.
Algo.~\ref{algo:check} details the aforementioned verification step.

  \begin{algorithm}[t]\caption{Dynamic programming procedure\label{algo:dyn_prog}  \\\hspace{\textwidth}
		$  \left[\DynProg: (\dyn_t)_{t=0}^{\horizon-1}, (\cost_t)_{t=0}^\horizon, \initstate, \Bell \rightarrow (u^*_0; \ldots; u^*_{\horizon-1})\right]$.}
	\begin{algorithmic}[1]
		\State{{\bf Inputs}: Dynamics $(\dyn_t)_{t=0}^{\horizon-1}$, costs $(\cost_t)_{t=0}^\horizon$, initial state $\initstate$, procedure $\Bell$}
		\State{Initialize $\costogo_\horizon = \cost_\horizon$}
		\For{$t =\horizon-1, \ldots, 0$}
		\State{Compute $\costogo_t, \pi_t = \Bell(\dyn_t, \cost_t, \costogo_{t+1})$, store $\pi_t$ \label{line:dyn_prog_backward}} \Comment{{\it \small $\Bell=\BellLQ$ (Algo.~\ref{algo:BellLQ}) for linear quadratic control}}
		\EndFor
		\State{Initialize $\state_0^* = \initstate$}
		\For{$t=0, \ldots, \horizon-1$}
		\State{Compute $\ctrl_t^* = \pi_t(\state_t^*)$,\quad  $\state_{t+1}^* = \dyn_t(\state_t^*, \ctrl_t^*)$}
		\EndFor
		\State{{\bf Output:} Optimal command $\ctrls =( \ctrl_0^*; \ldots;\ctrl_{\horizon-1}^*)$ for problem~\eqref{eq:discrete_ctrl_pb}}
	\end{algorithmic}
\end{algorithm}

\begin{algorithm}\caption{Analytic solution of Bellman's equation~\eqref{eq:bellman_lq} for linear dynamics, quadratic costs \\\hspace{\textwidth}
  $\left[\BellLQ: \lin_t, \qua_t, \costogo_{ t+1} \rightarrow \costogo_t, \pi_t\right]$
  \label{algo:BellLQ}}
\begin{algorithmic}[1]
  \State{{\bf Inputs:}
    \begin{enumerate}
      \item  Linear function $\lin_t$ parameterized as	$\lin_t(\state, \ctrl) = \A_t \state + \B _t\ctrl$
      \item Quadratic function $\qua_t$ parameterized as $\qua_t(\state, \ctrl) = \frac{1}{2}\state^\top \H_t \state+  \frac{1}{2}\ctrl^\top \G_t \ctrl + \state^\top \R_t \ctrl+ \h_t^\top \state  +  \g_t^\top \ctrl$
      \item Quadratic function $\costogo_{ t+1}$ parameterized as $
      \costogo_{t+1}(\state) = \frac{1}{2}\state^\top \J_{t+1}\state + \j_{t+1}^\top \state + \jcst_{t+1}$
  \end{enumerate} }
  \State{Define the cost-to-go function $\costogo_t: \state \rightarrow 	\frac{1}{2}\state^\top \J_t\state + \j_t^\top \state + \jcst_t$ with 
    \begin{align}
      \J_t & = \H _t+ \A_t^\top  \J_{t+1}\A_t - (\R_t + \A_t^\top  \J_{t+1}\B_t) (\G_t + \B_t^\top \J_{t+1}\B_t)^{-1} (\R_t^\top + \B_t^\top  \J_{t+1} \A_t)\nonumber \\
      \j_t & = \h_t + \A_t^\top   \j_{t+1} - (\R _t+ \A_t^\top  \J_{t+1}\B_t) (\G_t + \B_t^\top  \J_{t+1}\B_t)^{-1} (\g_t+ \B_t^\top  \j_{t+1}),  \nonumber \\
      \jcst_t & =\jcst_{t+1} - \frac{1}{2}(\g_t + \B_t^\top \j_{t+1})^\top (\G_t + \B_t^\top\J_{t+1}\B_t)^{-1} ( \g_t + \B_t^\top \j_{t+1}) \nonumber
    \end{align}
  }		
  \State{Define the policy $\pi_t: \state \rightarrow K_t\state + k_t$ with 
    \begin{align*}
      K_t  & = -(\G_t + \B_t^\top\J_{t+1}\B_t)^{-1} (\R_t^\top + \B_t^\top \J_{t+1} \A_t), \qquad k_t  = -(\G_t + \B_t^\top\J_{t+1}\B_t)^{-1} (\g_t + \B_t^\top \j_{t+1})  \nonumber
    \end{align*}
  }
  \State{{\bf Output:} Cost-to-go $\costogo_t$ and policy $\pi_t$ at time $t$}
\end{algorithmic}
\end{algorithm}

\begin{algorithm}\caption{Analytic solution of Bellman's equation~\eqref{eq:bellman_l} for linear dynamics, linear regularized costs \label{algo:bp}
  \\\hspace{\textwidth}
  $\left[\BellL: \ell_t^\dyn, \ell_t^\cost, \costogo_{ t+1}, \nu \rightarrow \costogo_t, \pi_t\right]$
}
\begin{algorithmic}[1]
  \State {\bf Inputs:} 
  \begin{enumerate}
    \item Linear function $\ell_\dyn$ parameterized as	$\ell_t^\dyn(\state, \ctrl) = \A_t \state + \B_t \ctrl$
    \item Linear function $\ell_\cost$ parameterized as $\ell_t^\cost(\state, \ctrl) =   \h_t^\top \state  +  \g_t^\top \ctrl $
    \item Affine function $\costogo_{t+1}$ parameterized as $
    \costogo_{t+1}(\state) = \j_{t+1}^\top \state + \jcst_{t+1}$
    \item Regularization $\nu\geq 0$
  \end{enumerate} 
  \State{Define $\costogo_t: \state \rightarrow 	 \j_t^\top \state + \jcst_t$ with 
    $
    \j_t = \h_t + \A_t^\top   \j_{t+1} ,  \ 
    \jcst_t  =\jcst_{t+1} - \|\g_t + \B_t^\top \j_{t+1}\|_2^2/(2\reg). \nonumber
    $
  }		
  \State{Define $\pi_t: \state \rightarrow k_t$ with 
    $
    k_t  = - (\g_t + \B_t^\top \j_{t+1})/\reg.  
    $
  }
  \State{{\bf Output:} Cost-to-go $\costogo_t$ and policy $\pi_t$ at time $t$}
\end{algorithmic}
\end{algorithm}

\begin{algorithm}\caption{Check if subproblems given by $\qua_t(\auxstate, \cdot) + \costogo_{t+1}(\lin_t(\auxstate, \cdot))$ are valid for solving Bellman's equation~\eqref{eq:bellman_lq} \label{algo:check}
  \\\hspace{\textwidth}
  $\left[\CheckSubPb: \lin_t, \qua_t, \costogo_{ t+1} \rightarrow \valid \in \{\True, \False\}\right]$
}
\begin{algorithmic}[1]
\State{\bf Option:} Check strong convexity of subproblems or check only if the result gives a descent direction
  \State{{\bf Inputs:}
\begin{enumerate}
  \item  Linear function $\lin_t$ parameterized as	$\lin_t(\state, \ctrl) = \A_t \state + \B _t\ctrl$, 
  \item Quadratic function $\qua_t$ parameterized as $\qua_t(\state, \ctrl) = \frac{1}{2}\state^\top \H_t \state+  \frac{1}{2}\ctrl^\top \G_t \ctrl + \state^\top \R_t \ctrl+ \h_t^\top \state  +  \g_t^\top \ctrl$
  \item Quadratic function $\costogo_{ t+1}$ parameterized as $
  \costogo_{t+1}(\state) = \frac{1}{2}\state^\top \J_{t+1}\state + \j_{t+1}^\top \state + \jcst_{t+1}.$
\end{enumerate} }
\If{check strong convexity}
\State Compute the eigenvalues $\lambda_1 \leq \ldots\leq \lambda_{\dimctrl}$ of $Q_t +B_t^\top J_{t+1}B_t$
\IfThenElse{$\lambda_1>0$}{$\valid=\True$}{$\valid=\False$}
\ElsIf{check descent direction}
\State Compute $\jcst_t - \jcst_{t+1} =  -\frac{1}{2}(\q_t + \B_t^\top \j_{t+1})^\top (\Q_t + \B_t^\top \J_{t+1}\B_t)^{-1}(\q_t + \B_t^\top \j_{t+1})$
\IfThenElse{$\jcst_t - \jcst_{t+1} <0$}{$\valid=\True$}{$\valid=\False$}
\EndIf 
\State {\bf Output:} $\valid$
\end{algorithmic}	
\end{algorithm}

\begin{algorithm}
\caption{Forward pass \label{algo:forward} \\\hspace{\textwidth}
  $\left[\Forward: \ctrls,(\dyn_t)_{t=0}^{\horizon-1}, (\cost_t)_{t=0}^\horizon, \initstate, o_\dyn, o_\cost  \rightarrow   \obj(\ctrls), (\model_{\dyn_t}^{\state_t, \ctrl_t})_{t=0}^{\horizon-1}, (\model_{\cost_t}^{\state_t, \ctrl_t})_{t=0}^{\horizon-1}, \model_{\cost_\horizon}^{\state_\horizon} \right]$}
\begin{algorithmic}[1]
  \State{{\bf Inputs:} Command $\ctrls=  (\ctrl_0; \ldots; \ctrl_{\horizon-1})$, dynamics $(\dyn_t)_{t=0}^{\horizon-1}$, costs $(\cost_t)_{t=0}^\horizon$, initial state $\initstate$, order of the information to collect on the dynamics $o_\dyn \in \{0, 1, 2\}$  and the costs $o_\cost \in \{0,1 , 2\}$}
  \State{Initialize $\state_0 = \initstate$, $\obj(\ctrls) = 0$}
  \For{$t=0, \ldots \horizon-1$}
  \State{Compute $\cost_t(\state_t, \ctrl_t)$, update $\obj(\ctrls) \leftarrow \obj(\ctrls) + \cost_t(\state_t, \ctrl_t)$}
  \IfThen{$o_\cost \geq 1$}{Compute and store $\nabla \cost_t(\state_t, \ctrl_t)$ defining $\lin_{\cost_t}^{\state_t, \ctrl_t}$ as in~\eqref{eq:reminder_approx}}
  \IfThen{$o_\cost = 2$}{Compute and store $\nabla^2 \cost_t(\state_t, \ctrl_t)$ defining, with $\nabla \cost_t(\state_t, \ctrl_t)$, $\qua_{\cost_t}^{\state_t, \ctrl_t}$ as in~\eqref{eq:reminder_approx}}
  \State{Compute $\state_{t+1} = \dyn_t(\state_t, \ctrl_t)$}
  \IfThen{$o_\dyn \geq 1$}{Compute and store $\nabla \dyn_t(\state_t, \ctrl_t)$ defining $\lin_{\dyn_t}^{\state_t, \ctrl_t}$ as in~\eqref{eq:reminder_approx}}
  \IfThen{$o_\dyn = 2$}{Compute and store $\nabla^2 \dyn_t(\state_t, \ctrl_t)$ defining, with $\nabla \dyn_t(\state_t, \ctrl_t)$, $\qua_{\dyn_t}^{\state_t, \ctrl_t}$ as in~\eqref{eq:reminder_approx}}
  \EndFor
  \State{Compute $\cost_\horizon(\state_\horizon)$, update $\obj(\ctrls) \leftarrow \obj(\ctrls) + \cost_\horizon(\state_\horizon)$}
  \IfThen{$o_\cost \geq 1$}{Compute and store $\nabla \cost_\horizon(\state_\horizon)$ defining $\lin_{\cost_\horizon}^{\state_\horizon}$ as in~\eqref{eq:reminder_approx}}
  \IfThen{$o_\cost = 2$}{Compute and store $\nabla^2\cost_\horizon(\state_\horizon)$ defining, with $\nabla \cost_\horizon(\state_\horizon)$,  $\qua_{\cost_\horizon}^{\state_\horizon}$ as in~\eqref{eq:reminder_approx}}
  \State{{\bf Output:} Total cost $\obj(\ctrls)$}
  
  \State{{\bf Stored:} (if $o_f$ and $o_h$ non-zeros) Approximations $(\model_{\dyn_t}^{\state_t, \ctrl_t})_{t=0}^{\horizon-1}, (\model_{\cost_t}^{\state_t, \ctrl_t})_{t=0}^{\horizon-1}, \model_{\cost_\horizon}^{\state_\horizon}$ defined by 
    \[
    \model_{\dyn_t}^{\state_t, \ctrl_t} = \begin{cases} 
      \lin_{\dyn_t}^{\state_t, \ctrl_t} &\mbox{if} \ o_\dyn = 1 \\
      \qua_{\dyn_t}^{\state_t, \ctrl_t} & \mbox{if} \ o_\dyn = 2
    \end{cases}, \quad  	\model_{\cost_t}^{\state_t, \ctrl_t} = \begin{cases}
      \lin_{\cost_t}^{\state_t, \ctrl_t} &\mbox{if} \ o_\cost = 1 \\
      \qua_{\cost_t}^{\state_t, \ctrl_t} & \mbox{if} \ o_\cost = 2
    \end{cases}, \quad \model_{\cost_\horizon}^{\state_\horizon} = \begin{cases}
      \lin_{\cost_\horizon}^{\state_\horizon} &\mbox{if} \ o_\cost = 1 \\
      \qua_{\cost_\horizon}^{\state_\horizon} & \mbox{if} \ o_\cost = 2
    \end{cases}
    \]
  }
\end{algorithmic}
\end{algorithm}

\begin{algorithm}\caption{Backward pass for gradient oracle \label{algo:backward_grad} \\\hspace{\textwidth}
  $\left[ \Backward_{\operatorname{GD}}: (\lin_{\dyn_t}^{\state_t, \ctrl_t})_{t=0}^{\horizon-1}, (\lin_{\cost_t}^{\state_t, \ctrl_t})_{t=0}^{\horizon-1}, \lin_{\cost_\horizon}^{\state_\horizon}, \reg) \rightarrow (\pi_t)_{t=0}^{\horizon-1},\costogo_0\right]$}
\begin{algorithmic}[1]
  \State{{\bf Inputs}: Linear expansions of the dynamics $(\lin_{\dyn_t}^{\state_t, \ctrl_t})_{t=0}^{\horizon-1}$, linear expansions of the costs $(\lin_{\cost_t}^{\state_t, \ctrl_t})_{t=0}^{\horizon-1}, \lin_{\cost_\horizon}^{\state_\horizon}$, regularization 
    $\reg >0$}
  \State{Initialize $\costogo_\horizon = \lin_{\cost_\horizon}^{\state_\horizon}$
  }
  \For{$t=\horizon-1, \ldots 0$}
  \State Define $\lin_t = \lin_{\dyn_t}^{\state_t, \ctrl_t}$,
  \quad $\qua_t: \diffstate_t, \diffctrl_t \rightarrow  \lin_{\cost_t}^{\state_t, \ctrl_t}(\diffstate_t, \diffctrl_t) + \frac{\reg}{2} \|\diffctrl_t\|_2^2$
  \State Compute $\costogo_t, \pi_t = \BellLQ(\lin_t, \qua_t, \costogo_{ t+1} ) = \BellL (\lin_{\dyn_t}^{\state_t, \ctrl_t}, \lin_{\cost_t}^{\state_t, \ctrl_t}, \costogo_{ t+1}, \nu)$ where $\BellL$ is given in Algo.~\ref{algo:bp}
  \EndFor
  \State{\bf Outputs:} Policies $ (\pi_t)_{t=0}^{\horizon-1}$, cost-to-go function at initial time $\costogo_0$
\end{algorithmic}
\end{algorithm}

\begin{algorithm}\caption{Backward pass for Gauss-Newton oracle \label{algo:backward_gn} \\\hspace{\textwidth}
  $	\left[ \Backward_{\operatorname{GN}}: (\lin_{\dyn_t}^{\state_t, \ctrl_t})_{t=0}^{\horizon-1}, (\qua_{\cost_t}^{\state_t, \ctrl_t})_{t=0}^{\horizon-1}, \qua_{\cost_\horizon}^{\state_\horizon}, \reg) \rightarrow (\pi_t)_{t=0}^{\horizon-1},\costogo_0\right]$
}
\begin{algorithmic}[1]
  \State{{\bf Inputs}: Linear expansions of the dynamics $(\lin_{\dyn_t}^{\state_t, \ctrl_t})_{t=0}^{\horizon-1}$, quadratic expansions of the costs $(\qua_{\cost_t}^{\state_t, \ctrl_t})_{t=0}^{\horizon-1}, \qua_{\cost_\horizon}^{\state_\horizon}$, regularization 
    $\reg \geq0$}
  \State{Initialize $\costogo_\horizon = \qua_{\cost_\horizon}^{\state_\horizon}$
  }
  \For{$t=\horizon-1, \ldots 0$}
  \State Define $\lin_t = \lin_{\dyn_t}^{\state_t, \ctrl_t}$, \quad $\qua_t: \diffstate_t, \diffctrl_t \rightarrow \qua_{\cost_t}^{\state_t, \ctrl_t}(\diffstate_t, \diffctrl_t) + \frac{\reg}{2}\|\diffctrl_t\|_2^2$,  
  \If{ $\CheckSubPb(\lin_t, \qua_t, \costogo_{ t+1}) \ \mbox{is} \ \True$}
  \State Compute $\costogo_t, \pi_t = \BellLQ(\lin_t, \qua_t, \costogo_{ t+1} )$ with $\BellLQ$ given in Algo.~\ref{algo:BellLQ}
  \Else
  \State $\pi_s:\state\rightarrow 0$ for $s\leq t$,  $\costogo_0: \state \rightarrow  -\infty$, {\bf break}
  \EndIf
  \EndFor
  \State{\bf Outputs:} Policies $ (\pi_t)_{t=0}^{\horizon-1}$, cost-to-go function at initial time $\costogo_0$
\end{algorithmic}
\end{algorithm}
\begin{algorithm}\caption{Backward pass for Newton oracle \label{algo:backward_newton}
  \\\hspace{\textwidth}
  $	\left[ \Backward_{\operatorname{NE}}: (\qua_{\dyn_t}^{\state_t, \ctrl_t})_{t=0}^{\horizon-1}, (\qua_{\cost_t}^{\state_t, \ctrl_t})_{t=0}^{\horizon-1}, \qua_{\cost_\horizon}^{\state_\horizon}, \reg) \rightarrow (\pi_t)_{t=0}^{\horizon-1},\costogo_0\right]$
}
\begin{algorithmic}[1]
  \State{{\bf Inputs}: Quadratic expansions of the dynamics $(\qua_{\dyn_t}^{\state_t, \ctrl_t})_{t=0}^{\horizon-1}$, quadratic expansions of the costs $(\qua_{\cost_t}^{\state_t, \ctrl_t})_{t=0}^{\horizon-1}, \qua_{\cost_\horizon}^{\state_\horizon}$,  regularization 
    $\reg \geq0$}
  \State{Initialize $\costogo_\horizon = \qua_{\cost_\horizon}^{\state_\horizon}$, $\costate_\horizon = \nabla \cost_\horizon(\state_\horizon)$
  }
  \For{$t=\horizon-1, \ldots 0$}
  \State Define $\lin_t = \lin_{\dyn_t}^{\state_t, \ctrl_t}$, \quad $q_t: (\diffstate_t, \diffctrl_t) \rightarrow 
  \qua_{\cost_t}^{\state_t, \ctrl_t}(\diffstate_t, \diffctrl_t) 
  + \frac{\reg}{2} \|\diffctrl_t\|_2^2  + \frac{1}{2}\nabla^2\dyn_t(\state_t, \ctrl_t)[\cdot, \cdot,  \costate_{t+1}](\diffstate_t, \diffctrl_t)$
  \State{Compute $\costate_t =\nabla_{\state_t} \cost_t(\state_t, \ctrl_t) +   \nabla_{\state_t} \dyn_t(\state_t, \ctrl_t) \costate_{t+1} $}
  \If{ $\CheckSubPb(\lin_t, \qua_t, \costogo_{ t+1}) \ \mbox{is} \ \True$}
  \State Compute $\costogo_t, \pi_t = \BellLQ(\lin_t, \qua_t, \costogo_{ t+1} )$ with $\BellLQ$ given in Algo.~\ref{algo:BellLQ}
  \Else
  \State $\pi_s:\state\rightarrow 0$ for $s\leq t$,  $\costogo_0: \state \rightarrow  -\infty$, {\bf break}
  \EndIf
  \EndFor
  \State{\bf Outputs:} Policies $ (\pi_t)_{t=0}^{\horizon-1}$, cost-to-go function at initial time $\costogo_0$
\end{algorithmic}
\end{algorithm}

\begin{algorithm}\caption{Backward pass for a DDP approach with quadratic approximations \label{algo:backward_ddp} \\\hspace{\textwidth}
  $	\left[ \Backward_{\operatorname{DDP}}: (\qua_{\dyn_t}^{\state_t, \ctrl_t})_{t=0}^{\horizon-1}, (\qua_{\cost_t}^{\state_t, \ctrl_t})_{t=0}^{\horizon-1}, \qua_{\cost_\horizon}^{\state_\horizon}, \reg) \rightarrow (\pi_t)_{t=0}^{\horizon-1},\costogo_0\right]$
}
\begin{algorithmic}[1]
  \State{{\bf Inputs}: Quadratic expansions on the dynamics $(\qua_{\dyn_t}^{\state_t, \ctrl_t})_{t=0}^{\horizon-1}$, quadratic expansions on the costs $(\qua_{\cost_t}^{\state_t, \ctrl_t})_{t=0}^{\horizon-1}, \qua_{\cost_\horizon}^{\state_\horizon}$,  regularization 
    $\reg \geq0$}
  \State{Initialize $\costogo_\horizon = \qua_{\cost_\horizon}^{\state_\horizon}$
  }
  \For{$t=\horizon-1, \ldots 0$}
  \State{Define $\lin_t = \lin_{\dyn_t}^{\state_t, \ctrl_t}$, \quad  $q_t: \auxstate_t, \auxctrl_t \rightarrow \qua_{\cost_t}^{\state_t, \ctrl_t}(\auxstate_t, \auxctrl_t ) +\frac{\reg}{2} \|\auxstate_t\|_2^2+\frac{1}{2} \nabla^2\dyn_t(\state_t, \ctrl_t)[\cdot, \cdot,  \nabla\costogo_{t+1}(0)](\auxstate_t, \auxctrl_t) $ }
  \If{ $\CheckSubPb(\lin_t, \qua_t, \costogo_{ t+1}) \ \mbox{is} \ \True$}
  \State Compute $\costogo_t, \pi_t = \BellLQ(\lin_t, \qua_t, \costogo_{ t+1} )$ with $\BellLQ$ given in Algo.~\ref{algo:BellLQ}
  \Else
  \State $\pi_s:\state\rightarrow 0$ for $s\leq t$,  $\costogo_0: \state \rightarrow  -\infty$, {\bf break}
  \EndIf
  \EndFor
  \State{{\bf Outputs:} Policies $ (\pi_t)_{t=0}^{\horizon-1}$, cost-to-go function at initial time $\costogo_0$}
\end{algorithmic}
\end{algorithm}

\begin{algorithm}	\caption{Backward pass for Newton oracle with function storage \label{algo:backward_newton_simplified} }
\begin{algorithmic}[1]
  \State{ {\bf Inputs}: 
    Stored functions $(\dyn_t)_{t=0}^{\horizon-1}$, costs $(\cost_t)_{t=0}^\horizon$, inputs $(\ctrl_t)_{t=0}^{\horizon-1}$ with associated trajectory $(\state_t)_{t=0}^\horizon$
  }
  \State Compute the quadratic expansion $\qua_{\cost_\horizon}^{\state_\horizon}$ of the final cost and the derivative $\nabla \cost_\horizon(\state_\horizon)$ of the final cost  on $\state_\horizon$
  \State Set  $\costogo_\horizon = \qua_{\cost_\horizon}^{\state_\horizon}$, $\costate_\horizon = \nabla \cost_\horizon(\state_\horizon)$
  \For{$t=\horizon-1, \ldots 0$}
  \State Compute the linear approximation $ \lin_{\dyn_t}^{\state_t, \ctrl_t}$  of the dynamic around $\state_t, \ctrl_t$
  \State Compute the quadratic approximation $\qua_{\cost_t}^{\state_t, \ctrl_t}$ of the cost around $\state_t, \ctrl_t$
  \State Compute the Hessian of $\state_t, \ctrl_t \rightarrow \dyn_t(\state_t, \ctrl_t)^\top \costate_{t+1}$ on $\state_t, \ctrl_t$ which gives $\frac{1}{2}\nabla^2\dyn_t(\state_t, \ctrl_t)[\cdot, \cdot,  \costate_{t+1}]$.  
  \State Define $\lin_t = \lin_{\dyn_t}^{\state_t, \ctrl_t}$, \quad $q_t: (\diffstate_t, \diffctrl_t) \rightarrow 
  \qua_{\cost_t}^{\state_t, \ctrl_t}(\diffstate_t, \diffctrl_t) 
  + \frac{\reg}{2} \|\diffctrl_t\|_2^2  + \frac{1}{2}\nabla^2\dyn_t(\state_t, \ctrl_t)[\cdot, \cdot,  \costate_{t+1}](\diffstate_t, \diffctrl_t)$
  \State{Compute $\costate_t =\nabla_{\state_t} \cost_t(\state_t, \ctrl_t) +   \nabla_{\state_t} \dyn_t(\state_t, \ctrl_t) \costate_{t+1} $}
  \If{ $\CheckSubPb(\lin_t, \qua_t, \costogo_{ t+1}) \ \mbox{is} \ \True$}
  \State{Compute $\costogo_t, \pi_t = \BellLQ(\lin_t, \qua_t, \costogo_{ t+1} )$}
  \Else
  \State $\pi_s:\state\rightarrow 0$ for $s\leq t$,  $\costogo_0: \state \rightarrow  -\infty$, {\bf break}
  \EndIf
  \EndFor
  \State{\bf Outputs:} Policies $ (\pi_t)_{t=0}^{\horizon-1}$, cost-to-go function at initial time $\costogo_0$
\end{algorithmic}
\end{algorithm}

\begin{algorithm}\caption{Roll-out on dynamics \label{algo:roll_out} \\\hspace{\textwidth}
  $\left[ 
  \Roll:
  \auxstate_0, (\pi_t)_{t=1}^{\horizon-1},(\auxdyn_t) _{t=0}^{\horizon-1} \rightarrow  \auxctrls 
  \right]$}
\begin{algorithmic}[1]
  \State {\bf Inputs}:  Initial state $\auxstate_0$, sequence of policies $(\pi_t)_{t=0}^{\horizon-1}$, dynamics to roll-on $(\auxdyn_t) _{t=0}^{\horizon-1}$
  \For{$t=0, \ldots, \horizon-1$}
  \State Compute and store
  $
  \auxctrl_t= \pi_t(\auxstate_t), \ \auxstate_{t+1}= \auxdyn_t(\auxstate_t, \auxctrl_t).
  $
  \EndFor
  \State {\bf Output:} 
  Sequence of controllers $\auxctrls = (\auxctrl_0;\ldots; \auxctrl_{\horizon-1})$
\end{algorithmic}
\end{algorithm}

\begin{algorithm}\caption{Gradient oracle \label{algo:gd}
  \\\hspace{\textwidth}
  $\left[ 
  \operatorname{GD}:
  \ctrls, (\dyn_t)_{t=0}^{\horizon-1}, (\cost_t)_{t=0}^\horizon, \initstate, \reg \rightarrow  \auxctrls 
  \right]$}
\begin{algorithmic}[1]
  \State {\bf Inputs:} Command $\ctrls{=} (\ctrl_0; \ldots;\ctrl_{\horizon-1})$,  dynamics $(\dyn_t)_{t=0}^{\horizon-1}$, costs $(\cost_t)_{t=0}^\horizon$, initial state $\initstate$, regularization $\reg{>}0$
  \State Compute with Algo.~\ref{algo:forward}
  \[
  \obj(\ctrls),  (\lin_{\dyn_t}^{\state_t, \ctrl_t})_{t=0}^{\horizon-1}, (\lin_{\cost_t}^{\state_t, \ctrl_t})_{t=0}^{\horizon-1}, \lin_{\cost_\horizon}^{\state_\horizon} = \Forward(\ctrls,(\dyn_t)_{t=0}^{\horizon-1}, (\cost_t)_{t=0}^\horizon, \initstate, o_\dyn=1, o_\cost=1)
  \]
  \State Compute with Algo.~\ref{algo:backward_grad}
  \[
  (\pi_t)_{t=0}^{\horizon-1}, \costogo_0 = \Backward_{\operatorname{GD}}((\lin_{\dyn_t}^{\state_t, \ctrl_t})_{t=0}^{\horizon-1}, (\lin_{\cost_t}^{\state_t, \ctrl_t})_{t=0}^{\horizon-1}, \qua_{\cost_\horizon}^{\state_\horizon}, \reg)
  \]
  \State Compute with Algo.~\ref{algo:roll_out}
  \[
  \auxctrls = \Roll(0, (\pi_t)_{t=0}^{\horizon-1}, (\lin_{\dyn_t}^{\state_t, \ctrl_t})_{t=0}^{\horizon-1})
  \]
  \State {\bf Output:} Gradient direction 
  $\auxctrls =\argmin_{\tilde \diffctrls \in \reals^{\horizon\dimctrl}}
  \left\{	\lin_{\cost \circ \augtrajfunc}^{\ctrls}(\tilde \diffctrls) +\frac{\reg}{2}\|\tilde \diffctrls\|_2^2\right\}
  = -\reg^{-1} \nabla (\cost\circ \augtraj)(\ctrls)$
\end{algorithmic}
\end{algorithm}

\begin{algorithm}\caption{Gauss-Newton oracle (ILQR)\label{algo:gn}
  \\\hspace{\textwidth}
  $\left[ 
  \operatorname{GN}:
  \ctrls, (\dyn_t)_{t=0}^{\horizon-1}, (\cost_t)_{t=0}^\horizon, \initstate, \reg \rightarrow  \auxctrls 
  \right]$
}
\begin{algorithmic}[1]
  \State {\bf Inputs:} Command $\ctrls{= }(\ctrl_0; \ldots;\ctrl_{\horizon-1})$,  dynamics $(\dyn_t)_{t=0}^{\horizon-1}$, costs $(\cost_t)_{t=0}^\horizon$, initial state $\initstate$, regularization $\reg{\geq} 0$
  \State Compute with Algo.~\ref{algo:forward}
  \[
  \obj(\ctrls),  (\lin_{\dyn_t}^{\state_t, \ctrl_t})_{t=0}^{\horizon-1}, (\qua_{\cost_t}^{\state_t, \ctrl_t})_{t=0}^{\horizon-1}, \qua_{\cost_\horizon}^{\state_\horizon} = \Forward(\ctrls,(\dyn_t)_{t=0}^{\horizon-1}, (\cost_t)_{t=0}^\horizon, \initstate, o_\dyn=1, o_\cost=2)
  \]
  \State Compute with Algo.~\ref{algo:backward_gn}
  \[
  (\pi_t)_{t=0}^{\horizon-1}, \costogo_0 = \Backward_{\operatorname{GN}}((\lin_{\dyn_t}^{\state_t, \ctrl_t})_{t=0}^{\horizon-1}, (\qua_{\cost_t}^{\state_t, \ctrl_t})_{t=0}^{\horizon-1}, \qua_{\cost_\horizon}^{\state_\horizon}, \reg)
  \]
  \State Compute with Algo.~\ref{algo:roll_out}
  \[
  \auxctrls = \Roll(0, (\pi_t)_{t=0}^{\horizon-1}, (\lin_{\dyn_t}^{\state_t, \ctrl_t})_{t=0}^{\horizon-1})
  \]
  \State {\bf Output:} If $\costogo_0(0) = +\infty$, returns $\infeasible$, otherwise returns Gauss-Newton direction 
  
  $\auxctrls =\argmin_{\tilde \diffctrls \in \reals^{\horizon\dimctrl}}
  \left\{	\qua_{\cost}^{\augtrajfunc(\ctrls)}(\lin_\augtrajfunc^\ctrls( \tilde\diffctrls))  + \frac{\reg}{2}\|\tilde\diffctrls\|_2^2 \right\}
  = -(\nabla \augtraj(\ctrls)\nabla^2\cost(\states, \ctrls) \nabla \augtraj(\ctrls)+ \reg\idm)^{-1} \nabla (\cost\circ \augtraj)(\ctrls)$
\end{algorithmic}
\end{algorithm}

\begin{algorithm}\caption{Newton oracle \label{algo:newton}
  \\\hspace{\textwidth}
  $\left[ 
  \operatorname{NE}:
  \ctrls, (\dyn_t)_{t=0}^{\horizon-1}, (\cost_t)_{t=0}^\horizon, \initstate, \reg \rightarrow  \auxctrls 
  \right]$
}
\begin{algorithmic}[1]
  \State {\bf Inputs:} Command $\ctrls{= }(\ctrl_0; \ldots;\ctrl_{\horizon-1})$,  dynamics $(\dyn_t)_{t=0}^{\horizon-1}$, costs $(\cost_t)_{t=0}^\horizon$, initial state $\initstate$, regularization $\reg{\geq} 0$
  \State Compute with Algo.~\ref{algo:forward}
  \[
  \obj(\ctrls),  (\qua_{\dyn_t}^{\state_t, \ctrl_t})_{t=0}^{\horizon-1}, (\qua_{\cost_t}^{\state_t, \ctrl_t})_{t=0}^{\horizon-1}, \qua_{\cost_\horizon}^{\state_\horizon} = \Forward(\ctrls,(\dyn_t)_{t=0}^{\horizon-1}, (\cost_t)_{t=0}^\horizon, \initstate, o_\dyn=2, o_\cost=2)
  \]
  \State Compute with Algo.~\ref{algo:backward_newton}
  \[
  (\pi_t)_{t=0}^{\horizon-1}, \costogo_0 = \Backward_{\operatorname{NE}}((\qua_{\dyn_t}^{\state_t, \ctrl_t})_{t=0}^{\horizon-1}, (\qua_{\cost_t}^{\state_t, \ctrl_t})_{t=0}^{\horizon-1}, \qua_{\cost_\horizon}^{\state_\horizon}, \reg)
  \]
  \State Compute with Algo.~\ref{algo:roll_out}
  \[
  \auxctrls = \Roll(0, (\pi_t)_{t=0}^{\horizon-1}, (\lin_{\dyn_t}^{\state_t, \ctrl_t})_{t=0}^{\horizon-1})
  \]
  \State {\bf Output:} 	
  If $\costogo_0(0) = +\infty$, returns $\infeasible$, otherwise returns  Newton direction 
  
  $\auxctrls =\argmin_{\tilde \diffctrls \in \reals^{\horizon\dimctrl}}
  \left\{\qua_{\cost\circ \augtrajfunc}^\ctrls(\tilde \diffctrls)  + \frac{\reg}{2}\|\tilde\diffctrls \|_2^2 \right\}
  = -(\nabla^2 (\cost\circ \augtraj)(\ctrls) + \reg \idm)^{-1}\nabla (\cost\circ \augtraj)(\ctrls)$
\end{algorithmic}
\end{algorithm}

\begin{algorithm}\caption{Differential dynamic  programming oracle with linear quadratic approximations (iLQR) \label{algo:ddp_lq}
  \\\hspace{\textwidth}
      $\left[ 
  \operatorname{DDP-LQ}:
  \ctrls, (\dyn_t)_{t=0}^{\horizon-1}, (\cost_t)_{t=0}^\horizon, \initstate, \reg \rightarrow  \auxctrls 
  \right]$	
}
\begin{algorithmic}[1]
  \State {\bf Inputs:} Command $\ctrls{= }(\ctrl_0; \ldots;\ctrl_{\horizon-1})$,  dynamics $(\dyn_t)_{t=0}^{\horizon-1}$, costs $(\cost_t)_{t=0}^\horizon$, initial state $\initstate$, regularization $\reg{\geq} 0$
  \State Compute with Algo.~\ref{algo:forward}
  \[
  \obj(\ctrls),  (\lin_{\dyn_t}^{\state_t, \ctrl_t})_{t=0}^{\horizon-1}, (\qua_{\cost_t}^{\state_t, \ctrl_t})_{t=0}^{\horizon-1}, \qua_{\cost_\horizon}^{\state_\horizon} = \Forward(\ctrls,(\dyn_t)_{t=0}^{\horizon-1}, (\cost_t)_{t=0}^\horizon, \initstate, o_\dyn=1, o_\cost=2)
  \]
  \State Compute with Algo.~\ref{algo:backward_gn}
  \[
  (\pi_t)_{t=0}^{\horizon-1}, \costogo_0 = \Backward_{\operatorname{GN}}((\lin_{\dyn_t}^{\state_t, \ctrl_t})_{t=0}^{\horizon-1}, (\qua_{\cost_t}^{\state_t, \ctrl_t})_{t=0}^{\horizon-1}, \qua_{\cost_\horizon}^{\state_\horizon}, \reg)
  \]
  \State Compute with Algo.~\ref{algo:roll_out}, for $\diff_{\dyn_t}^{\state_t, \ctrl_t}(\auxstate_t, \auxctrl_t) = \dyn(\state_t +\auxstate_t , \ctrl_t + \auxctrl_t) - \dyn(\state_t, \ctrl_t)$,
  \[
  \auxctrls = \Roll(0, (\pi_t)_{t=0}^{\horizon-1}, (\diff_{\dyn_t}^{\state_t, \ctrl_t})_{t=0}^{\horizon-1})
  \]
  \State {\bf Output:} If $\costogo_0(0) = +\infty$, returns $\infeasible$, otherwise returns DDP oracle with linear-quadratic approximations $\auxctrls$
\end{algorithmic}
\end{algorithm}

\begin{algorithm}\caption{Differential dynamic programming oracle with quadratic approximations (DDP) \label{algo:ddp_q}
  \\\hspace{\textwidth}
  $\left[ 
\operatorname{DDP-Q}:
\ctrls, (\dyn_t)_{t=0}^{\horizon-1}, (\cost_t)_{t=0}^\horizon, \initstate, \reg \rightarrow  \auxctrls 
\right]$
}
\begin{algorithmic}[1]
  \State {\bf Inputs:} Command $\ctrls{= }(\ctrl_0; \ldots;\ctrl_{\horizon-1})$,  dynamics $(\dyn_t)_{t=0}^{\horizon-1}$, costs $(\cost_t)_{t=0}^\horizon$, initial state $\initstate$, regularization $\reg{\geq} 0$
  \State Compute with Algo.~\ref{algo:forward}
  \[
  \obj(\ctrls),  (\qua_{\dyn_t}^{\state_t, \ctrl_t})_{t=0}^{\horizon-1}, (\qua_{\cost_t}^{\state_t, \ctrl_t})_{t=0}^{\horizon-1}, \qua_{\cost_\horizon}^{\state_\horizon} = \Forward(\ctrls,(\dyn_t)_{t=0}^{\horizon-1}, (\cost_t)_{t=0}^\horizon, \initstate, o_\dyn=2, o_\cost=2)
  \]
  \State Compute with Algo.~\ref{algo:backward_ddp}
  \[
  (\pi_t)_{t=0}^{\horizon-1}, \costogo_0 = \Backward_{\operatorname{DDP}}((\qua_{\dyn_t}^{\state_t, \ctrl_t})_{t=0}^{\horizon-1}, (\qua_{\cost_t}^{\state_t, \ctrl_t})_{t=0}^{\horizon-1}, \qua_{\cost_\horizon}^{\state_\horizon}, \reg)
  \]
  \State Compute with Algo.~\ref{algo:roll_out},  for $\diff_{\dyn_t}^{\state_t, \ctrl_t}(\auxstate_t, \auxctrl_t) = \dyn(\state_t +\auxstate_t , \ctrl_t + \auxctrl_t) - \dyn(\state_t, \ctrl_t)$,
  \[
  \auxctrls = \Roll(0, (\pi_t)_{t=0}^{\horizon-1}, (\diff_{\dyn_t}^{\state_t, \ctrl_t})_{t=0}^{\horizon-1})
  \]
  \State {\bf Output:} 	
  If $\costogo_0(0) = +\infty$, returns $\infeasible$, otherwise returns  DDP oracle with quadratic approximations  $\auxctrls$
\end{algorithmic}
\end{algorithm}

\begin{algorithm}\caption{Line-search\label{algo:line_search}
      \\\hspace{\textwidth}
      $\left[\linesearch:\ctrls, (\cost_t)_{t=0}^\horizon,  (\dyn_t)_{t=0}^{\horizon-1}, (\auxdyn_t)_{t=0}^{\horizon-1}, (\Pol:\stepsize \rightarrow (\pi_t^\stepsize)_{t=0}^{\horizon-1},\costogo_0^\stepsize )   \rightarrow \ctrls^{\nxt}\right]$}
\begin{algorithmic}[1]
  \State {\bf Option:} directional step or regularized step
  \State{{\bf Inputs:}  Current controls $\ctrls$, costs $(\cost_t)_{t=0}^\horizon$, initial state $\initstate$, original dynamics $(\dyn_t)_{t=0}^{\horizon-1}$, dynamics to roll out on $(\auxdyn_t)_{t=0}^{\horizon-1}$, family of policies and corresponding costs given by  $\stepsize \rightarrow (\pi_t^\stepsize)_{t=0}^{\horizon-1},\costogo_0^\stepsize$,  decreasing factor $\rho_{\dec} \in (0,1)$, increasing factor $\rho_{\inc}>1$, previous stepsize $\stepsize_{\prev}$ }
  \State{Compute $\obj(\ctrls) = \Forward(\ctrls,(\dyn_t)_{t=0}^{\horizon-1}, (\cost_t)_{t=0}^\horizon, \initstate, o_\dyn=0, o_\cost=0)$}
  \If{directional step}
  \State Initialize $\stepsize=1$
  \ElsIf{regularized step}
  \State Compute $\nabla \cost(\states, \ctrls)$ for $\states = \traj(\initstate, \ctrls)$
  \State  Initialize $\stepsize= \rho_{\inc} \stepsize_{\prev}/\|\nabla \cost(\states, \ctrls)\|_2$
  \EndIf
  \State{Initialize $\diffstate_0=0$, $\text{accept}=\text{False}$, minimal stepsize $\stepsize_{\min}=10^{-12}$}
  \While{not accept}
  \State{Get $\pi_t^{\stepsize}, \costogo_0^{\stepsize} = \Pol(\stepsize)$}
  \State{Compute $\diffctrls^\gamma = \Roll(\diffstate_0, (\pi_t^{\stepsize})_{t=1}^{\horizon-1},(\auxdyn_t) _{t=0}^{\horizon-1})$}
  \State{Set $\ctrls^{\nxt} = \ctrls + \diffctrls^{\stepsize}$}
  \State{Compute  $\obj(\ctrls^{\nxt}) = \Forward(\ctrls^{\nxt},(\dyn_t)_{t=0}^{\horizon-1}, (\cost_t)_{t=0}^\horizon, \initstate, o_\dyn=0, o_\cost=0)$}
  \IfThenElse{$\obj(\ctrls^{\nxt})  - \obj(\ctrls) \leq \costogo_0^{\stepsize}(0)$}{set $\text{accept}=\text{True}$}{set $\stepsize \rightarrow \rho_{\dec}\stepsize$}
  \IfThen{$\stepsize \leq \stepsize_{\min}$}{\textbf{break}}
  \EndWhile
  \IfThen{regularized step}{$\stepsize := \stepsize \|\nabla \cost(\states, \ctrls)\|_2$}
  \State {\bf Output:} Next sequence of controllers $\ctrls^{\nxt}$, store value of the stepsize selected $\stepsize$
\end{algorithmic}
\end{algorithm}

\begin{algorithm}\caption{Iterative Linear Quadratic Regulator/Gauss-Newton step with line-search on  descent directions \label{algo:gn_algo} }
\begin{algorithmic}[1]
  \State{{\bf Inputs:} Command $\ctrls$, dynamics $(\dyn_t)_{t=0}^{\horizon-1}$, costs $(\cost_t)_{t=0}^{\horizon-1}$, initial state $\initstate$}
  \State{Compute with Algo.~\ref{algo:forward}
    \[
    \obj(\ctrls), (\lin_{\dyn_t}^{\state_t, \ctrl_t})_{t=0}^{\horizon-1}, (\qua_{\cost_t}^{\state_t, \ctrl_t})_{t=0}^{\horizon-1}, \qua_{\cost_\horizon}^{\state_\horizon} = \Forward(\ctrls,(\dyn_t)_{t=0}^{\horizon-1}, (\cost_t)_{t=0}^\horizon, \initstate, o_\dyn=1, o_\cost=2)
    \]}
  \State{Compute with Algo.~\ref{algo:backward_gn}
    \[
    (\pi_t)_{t=0}^{\horizon-1}, \costogo_0 = \Backward_{\operatorname{GN}}((\lin_{\dyn_t}^{\state_t, \ctrl_t})_{t=0}^{\horizon-1}, (\qua_{\cost_t}^{\state_t, \ctrl_t})_{t=0}^{\horizon-1}, \qua_{\cost_\horizon}^{\state_\horizon}, 0)
    \]}
  \State Set $\reg=\reg_{\operatorname{init}}$ {with, e.g., $\reg_{\operatorname{init}} = 10^{-6}$}
  \While{$c_0(0) = +\infty$}
  \State Compute $(\pi_t)_{t=0}^{\horizon-1}, \costogo_0 = \Backward_{\operatorname{GN}}((\lin_{\dyn_t}^{\state_t, \ctrl_t})_{t=0}^{\horizon-1}, (\qua_{\cost_t}^{\state_t, \ctrl_t})_{t=0}^{\horizon-1}, \qua_{\cost_\horizon}^{\state_\horizon}, \reg)$
  \State Set $\reg \rightarrow \rho_{\operatorname{inc}}\reg$  {with, e.g., $\rho_{\operatorname{inc}} =10$}
  \EndWhile
  \State{Define $\Pol: \stepsize \rightarrow \left(\begin{array}{ll}
      (\pi_t^{\stepsize}: &\diffstate \rightarrow \stepsize\pi_t(0) + \nabla \pi_t(0)^\top \diffstate)_{t=0}^{\horizon-1}, \\
      \hspace{6pt} c_0^\stepsize: &\diffstate \rightarrow \stepsize c_0(\diffstate)
    \end{array}\right)$}
  \State{Compute with Algo.~\ref{algo:line_search}
    \[
    \ctrls^{\nxt} = \linesearch(\ctrls, (\cost_t)_{t=0}^\horizon,  (\dyn_t)_{t=0}^{\horizon-1}, (\lin_{\dyn_t}^{\state_t, \ctrl_t})_{t=0}^{\horizon-1}, \Pol)
    \]}
  \State{{\bf Output:} Next sequence of controllers $\ctrls^{\nxt}$}
\end{algorithmic}

\end{algorithm}

  \clearpage
\begin{figure}
	\begin{center}
		\includegraphics[width=\linewidth]{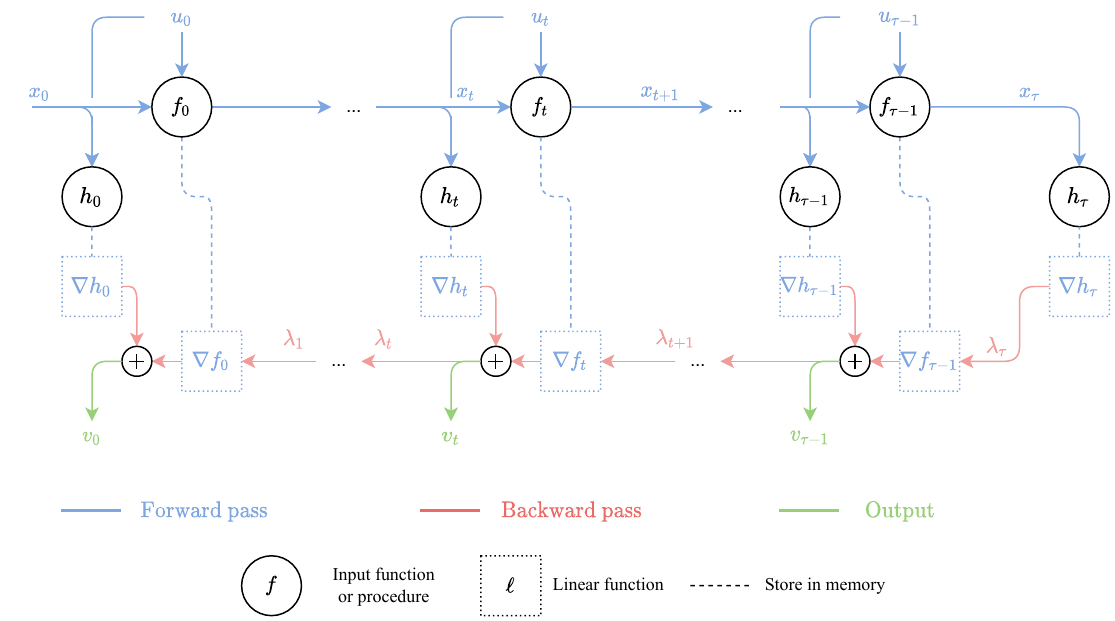}
		\caption{Computational scheme of a gradient oracle.  
			\label{fig:grad}}
	\end{center}
\end{figure}

\begin{figure}
	\begin{center}
		\includegraphics[width=\linewidth]{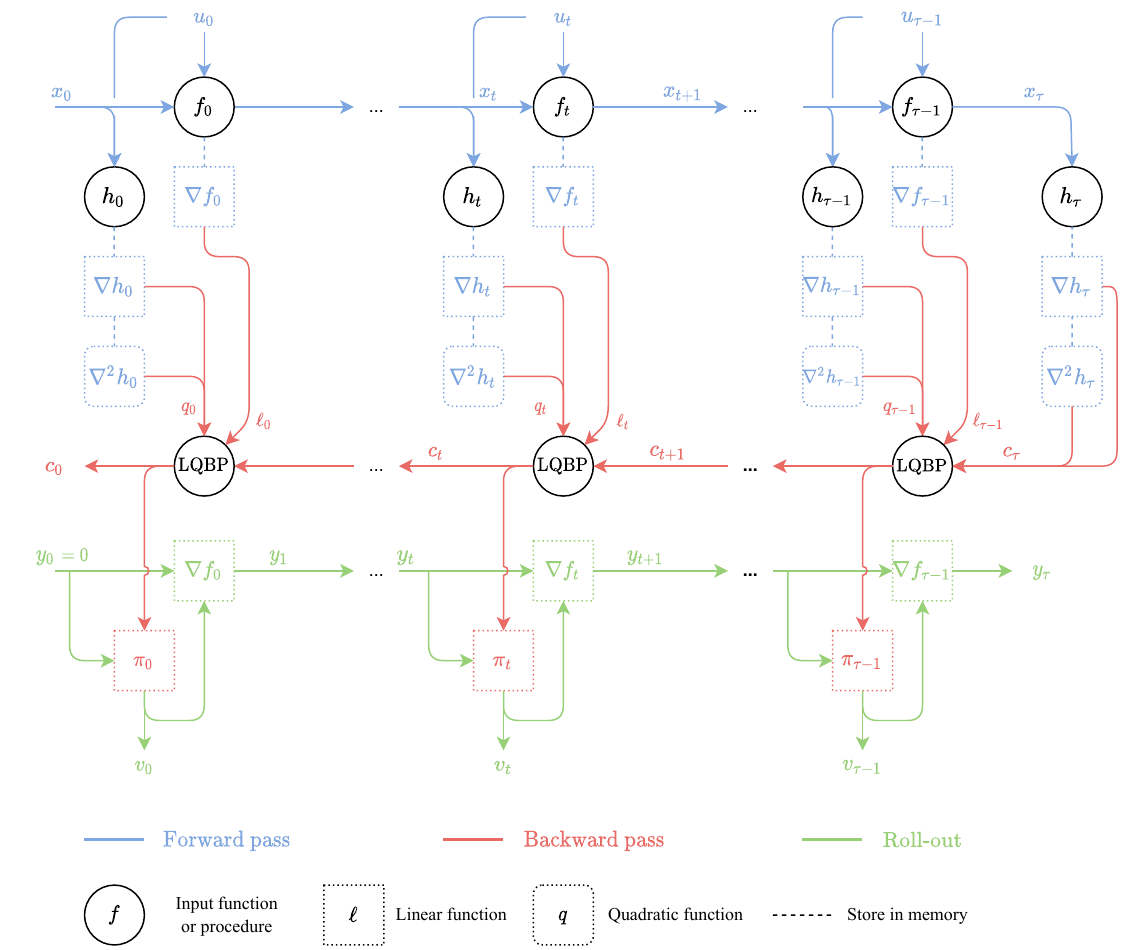}
	\end{center}
	\caption{Computational scheme of a Gauss-Newton oracle (ILQR).
		\label{fig:ilqr}}
\end{figure}

\begin{figure}
	\begin{center}
		\includegraphics[width=\linewidth]{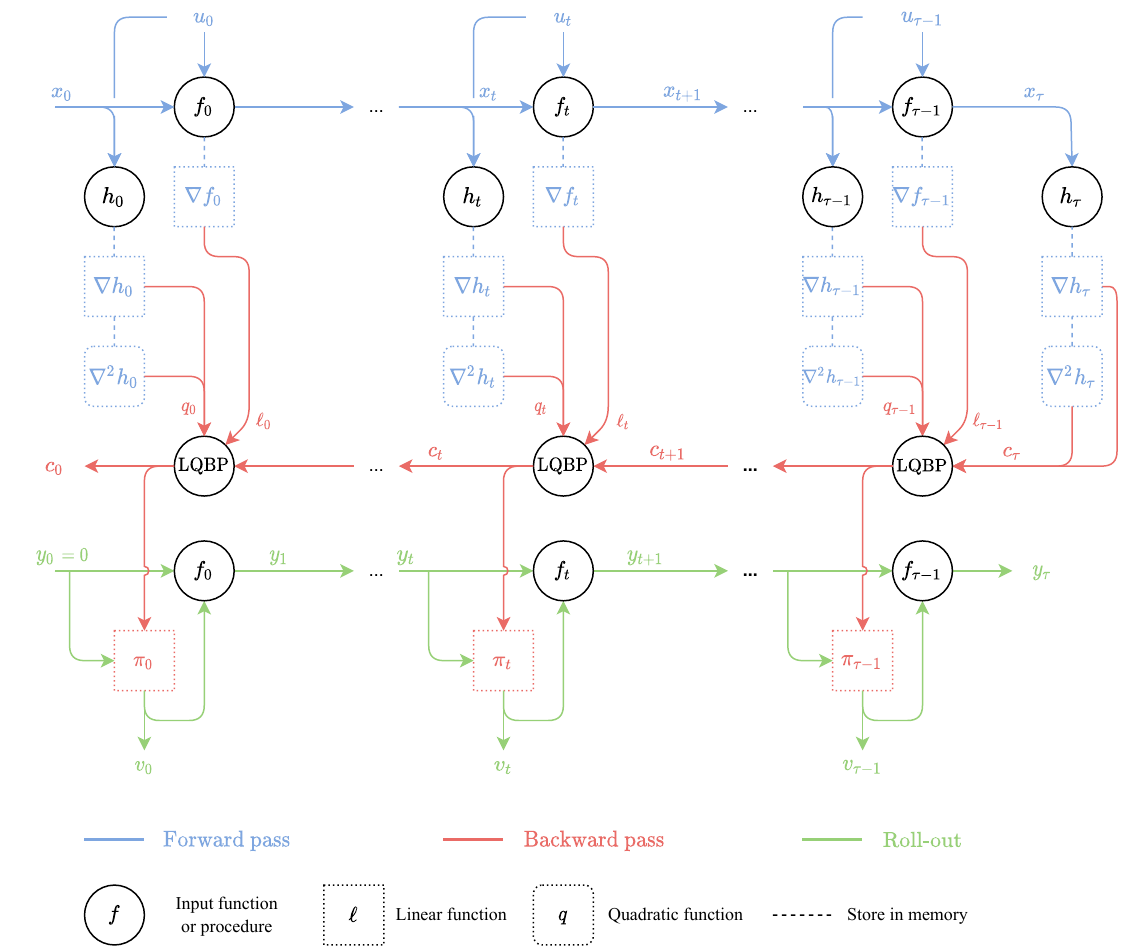}
		\caption{Computational scheme of a DDP oracle with linear quadratic approximations (iLQR).
			\label{fig:ddplq}}
	\end{center}
\end{figure}

\begin{figure}
	\begin{center}
		\includegraphics[width=\linewidth]{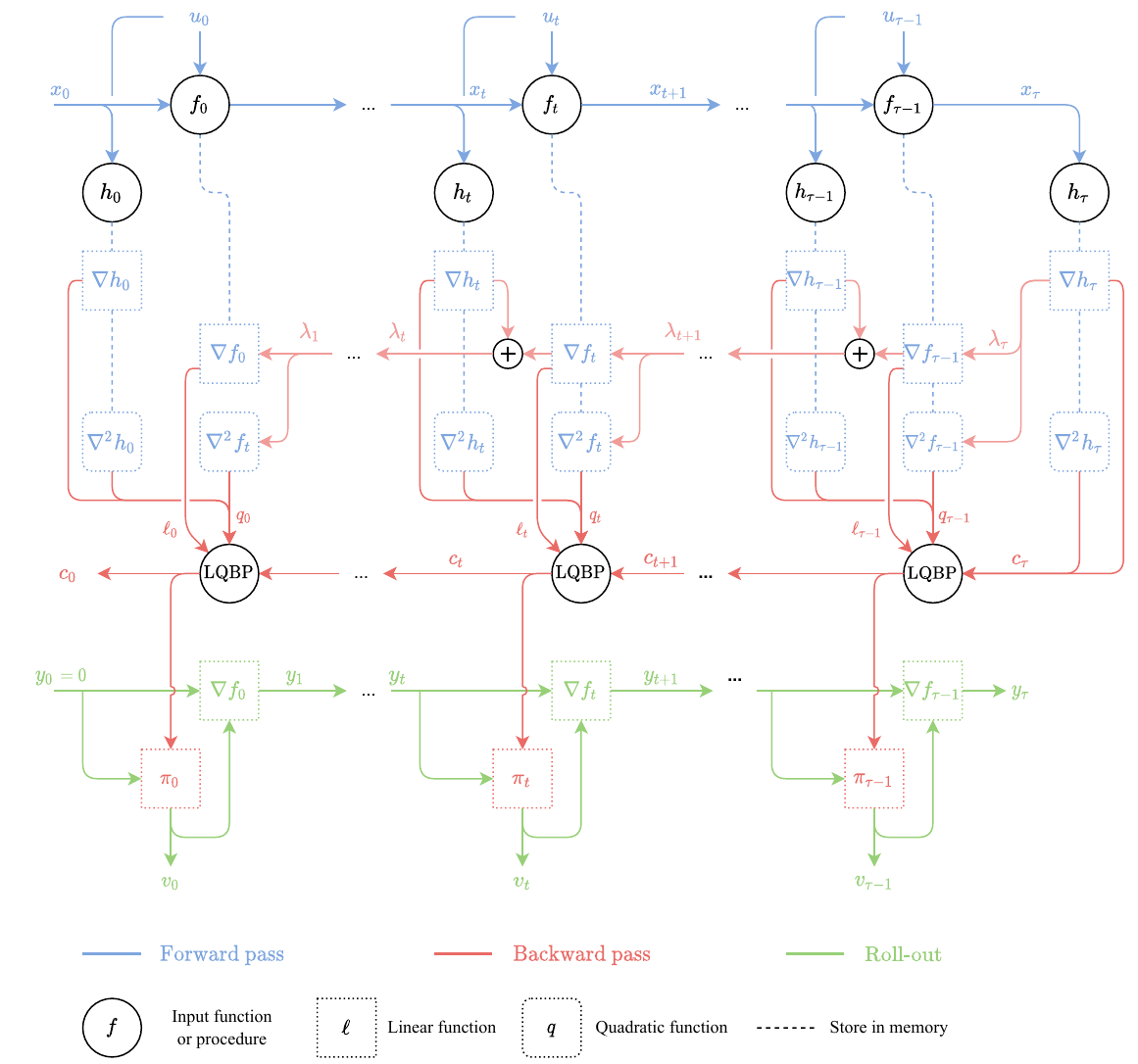}
	\end{center}
	\caption{Computational scheme of a Newton oracle. 
		 \label{fig:newton}}
\end{figure}

\begin{figure}
	\begin{center}
		\includegraphics[width=\linewidth]{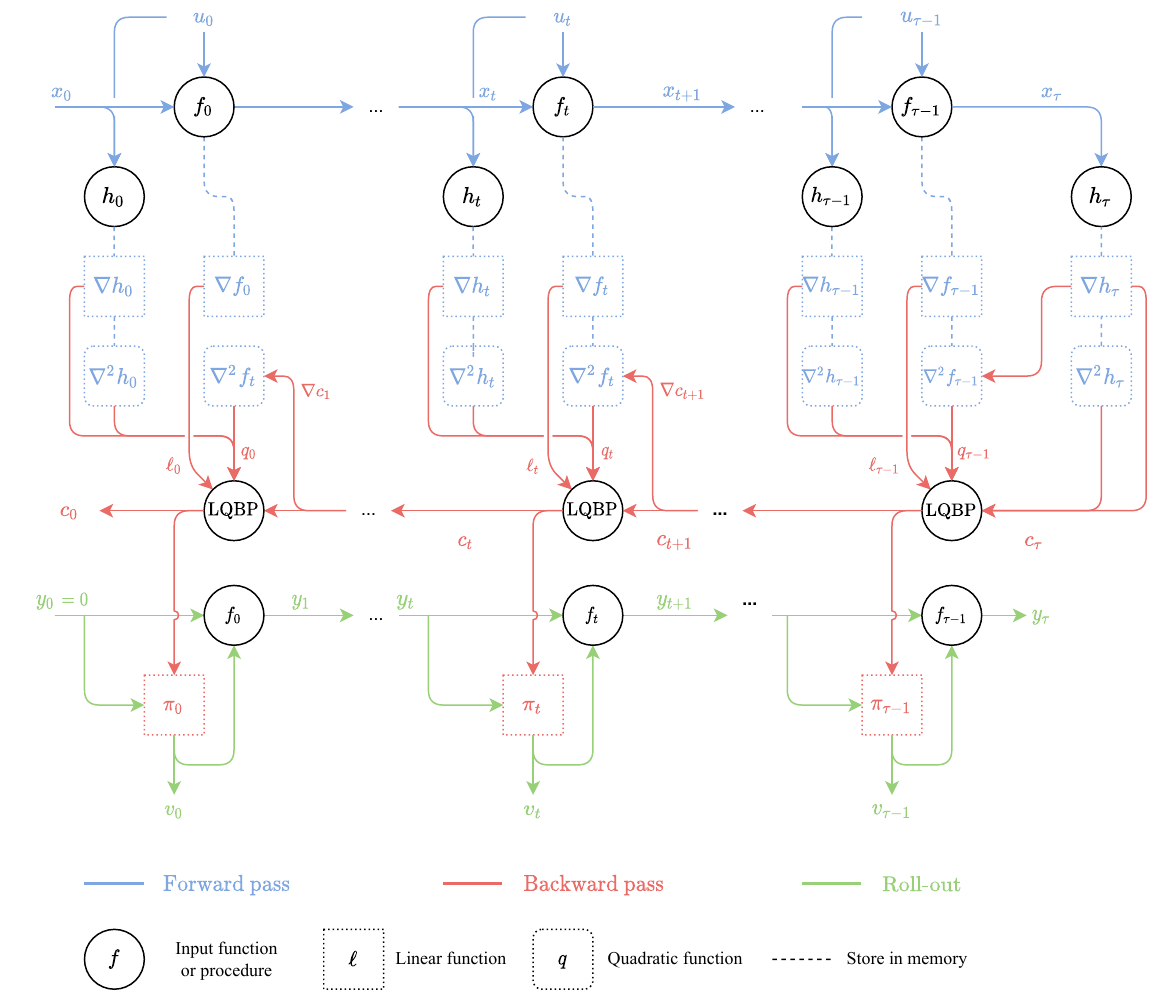}
		\caption{Computational scheme of a DDP oracle with quadratic approximations. (DDP)
			 \label{fig:ddpq}}
	\end{center}
\end{figure}

	\clearpage

  \section{Computational Complexity in a Differentiable Programming Framework}\label{app:cplxity}
  We detail here how to alleviate intermediate storing of second order information
to lower the computational cost of oracles based on quadratic approximations.

The time complexities of the forward pass presented in
Section~\ref{sec:comput_cplxity}, corresponding to the computations of the
gradients of the dynamics or the costs and Hessians of the costs, are then
incurred during the backward pass. A major difference lies in the computation of
the quadratic information of the dynamic required in quadratic oracles such as a
Newton oracle or a DDP oracle with quadratic approximations. Indeed, a closer
look at Algo.~\ref{algo:backward_newton} and Algo.~\ref{algo:backward_ddp} show
that only the Hessians of scalar functions of the form $\state, \ctrl
\rightarrow \dyn(\state, \ctrl)^\top \lambda$ need to be computed, which comes
at a cost $(\dimstate+\dimctrl)^2$. In comparison, the cost of computing the
second order information of $\dyn$ is $O((\dimstate+\dimctrl)^2\dimstate)$. As
an example, Algo.~\ref{algo:backward_newton_simplified} presents an
implementation of a Newton step using stored functions and inputs.

The computational complexities of the oracles when the dynamics and the costs
functions are stored in memory are presented in
Table~\ref{tab:efficient_complexities}. We
consider for simplicity that the memory cost of storing the information
necessary to evaluate a function $f:\reals^d\rightarrow \reals^n$ is $nd$ as it
is the case for a linear function $f$.

\begin{table}
	\begin{center}
		\bgroup
		\def\arraystretch{1.5}
		
		\begin{tabularx}{\linewidth}{p{60pt}|p{280pt}}
			\multicolumn{1}{c}{}	  & \multicolumn{1}{l}{Time complexities of the forward pass}\\
			\toprule
			All cases
			& $\horizon\Big(
			\underbrace{\dimstate^2 {+} \dimstate\dimctrl}_{\dyn_t} 
			{+} \underbrace{\dimstate {+} \dimctrl}_{\cost_t}
			\Big){=}O(\horizon(\dimstate^2 {+} \dimstate\dimctrl))$
			\\
			\bottomrule
		\end{tabularx}
		
		\vspace{1ex}
		
		\begin{tabularx}{\linewidth}{p{60pt}|p{280pt}}
			\multicolumn{1}{c}{}	  & \multicolumn{1}{l}{Space complexities of the forward pass}\\
			\toprule
			Function eval.
			& $0$
			\\
			All other cases
			& $\horizon\Big(
			\underbrace{\dimstate^2 {+} \dimstate\dimctrl}_{ \dyn_t }
			{+} \underbrace{\dimstate {+} \dimctrl}_{ \cost_t}
			\Big){=} O(\horizon(\dimstate^2 {+} \dimstate\dimctrl))$
			\\
			\bottomrule
		\end{tabularx}
		
		\vspace{1ex}
		\begin{tabularx}{\linewidth}{p{60pt}|p{280pt}}
			\multicolumn{1}{c}{} & \multicolumn{1}{l}{Time complexities of the backward passes}\\
			\toprule
			GD 
			& $\horizon\Big( 	\underbrace{\dimstate^2 {+} \dimstate\dimctrl}_{\nabla \dyn_t }
			{+} \underbrace{\dimstate {+} \dimctrl}_{\nabla \cost_t}{+}
			\underbrace{\dimstate^2 {+} \dimstate\dimctrl}_{\textrm{Roll}} 
			{+} \underbrace{\dimstate^2 {+} \dimstate\dimctrl
			}_{\BellL}
			\Big){=}O(\horizon(\dimstate^2 {+} \dimstate\dimctrl))$ 
			\\
			GN/DDP-LQ
			& $ \horizon\Big(	\underbrace{\dimstate^2 {+} \dimstate\dimctrl}_{\nabla \dyn_t} 
			{+} \underbrace{\dimstate {+} \dimctrl}_{\nabla \cost_t} {+} 
			\underbrace{\dimstate^2 {+} \dimctrl^2 {+} \dimstate\dimctrl}_{\nabla^2 \cost_t}\Big)\newline
			+	\horizon\Big(	\underbrace{\dimstate^2 {+} \dimstate\dimctrl}_{\textrm{Roll}} 
			{+} \underbrace{\dimstate^3 {+} \dimctrl^3 {+} \dimctrl^2\dimstate}_{\BellLQ}
			\Big) {=} O(\horizon(\dimstate {+} \dimctrl)^3)$
			\\
			NE/DDP-Q 
			& $ \horizon\Big(	\underbrace{\dimstate^2 {+} \dimstate\dimctrl}_{\nabla \dyn_t} 
			{+} \underbrace{\dimstate {+} \dimctrl}_{\nabla \cost_t} {+} 
			\underbrace{\dimstate^2 {+} \dimctrl^2 {+} \dimstate\dimctrl}_{\nabla^2 \cost_t}  +\underbrace{\dimstate^2 {+} \dimctrl^2 {+} \dimstate\dimctrl}_{\nabla^2(\dyn_t^\top \lambda)}\Big)\newline
			+ \horizon\Big(
			\underbrace{\dimstate^2 {+} \dimstate\dimctrl}_{\textrm{Roll}} 
			{+} \underbrace{\dimstate^3 {+} \dimctrl^3 {+} \dimctrl^2\dimstate}_{\BellLQ}\Big) {=} O(\horizon(\dimstate {+} \dimctrl)^3)$
			\\
			\bottomrule
		\end{tabularx}

		\egroup
		\caption{Space and time complexities of the oracles when storing functions as in, e.g., Algo.~\ref{algo:backward_newton_simplified}. 
			\label{tab:efficient_complexities}}
	\end{center}
\end{table}

  \section{Alternative Resolution of Linear-Quadratic Control Problem}\label{app:sparse_solvers}
  We presented the implementation of all algorithms in a unified viewpoint with
dynamic programming as the core subroutine. For classical optimization steps
such as Gauss-Newton or Newton, once the problem has been instantiated, as done
in Lemma~\ref{eq:lin_quad_oracle}, the resulting quadratic optimization
subproblem can be solved in several other ways. We present such
alternatives for completeness.

\subsection{Block Band Diagonal Underlying Structure}
The subproblems we are interested to solve are linear quadratic control problems
of the form
\begin{align*}
\min_{
  \substack{
    \state_0, \ldots, \state_\horizon \in \reals^\dimstate \\ 
    \ctrl_0, \ldots, \ctrl_{\horizon-1} \in \reals^\dimctrl
    }
  }
& \sum_{t=0}^{\horizon-1} \left(
    \frac{1}{2}\state_t^\top \H_t \state_t 
    + \frac{1}{2}\ctrl_t^\top \G_t \ctrl_t 
    + \state_t^\top \R_t \ctrl_t
    + \h_t^\top \state_t  
    + \g_t^\top \ctrl_t
  \right) 
  + \frac{1}{2}\state_\horizon^\top \H_\horizon \state_\horizon
  + \h_\horizon^\top \state_\horizon 
  \nonumber \\
\mbox{subject to} \quad 
& \state_{t+1} = A_t \state_t + B_t\ctrl_t,  
\quad \mbox{for} \ 
t \in \{0, \ldots, \horizon-1\}, \quad 
\state_0 = \initstate. \nonumber
\end{align*}
By introducing Lagrange multipliers $(\costate_t)_{t=0}^\horizon$ for the
constraints, the problem can be stated as follows.
\begin{align*}
& \min_{
  \substack{
    \state_0, \ldots, \state_\horizon \in \reals^\dimstate \\ 
    \ctrl_0, \ldots, \ctrl_{\horizon-1} \in \reals^\dimctrl
    }
  }
\sup_{
  \costate_0, \ldots, \costate_\horizon \in \reals^\dimstate
} 
\mathcal{L}(
  (\state_t)_{t=0}^\horizon,
  (\ctrl_t)_{t=0}^{\horizon-1},
  (\costate_t)_{t=0}^\horizon
  )
\\
\mbox{for} \ & \mathcal{L}(
  (\state_t)_{t=0}^\horizon,
  (\ctrl_t)_{t=0}^{\horizon-1},
  (\costate_t)_{t=0}^\horizon
  ) \\
& =
\sum_{t=0}^{\horizon-1} \left(
    \frac{1}{2}\state_t^\top \H_t \state_t 
    + \frac{1}{2}\ctrl_t^\top \G_t \ctrl_t 
    + \state_t^\top \R_t \ctrl_t
    + \h_t^\top \state_t  
    + \g_t^\top \ctrl_t
    + \costate_{t+1}^\top 
    (\state_{t+1} - A_t \state_t - B_t\ctrl_t)
  \right) \\
  & \hspace*{10pt} 
  + \costate_0^\top(\state_0 - \initstate)
  + \frac{1}{2}\state_\horizon^\top \H_\horizon \state_\horizon
  + \h_\horizon^\top \state_\horizon.
  \nonumber 
\end{align*}
The optimality conditions, a.k.a. KKT conditions, are  
\begin{align*}
  x_0 - \initstate 
  & = 0
  & (\partial_{\costate_0} \mathcal{L} = 0)\ \\
  P_t\state_t + R_t\ctrl_t + p_t - A_t^\top \costate_{t+1} + \costate_t
  & = 0 
  & t \in \{0, \ldots, \horizon -1\} \ (\partial_{\state_t} \mathcal{L} = 0)\ \\
  Q_t \ctrl_t + R_t^\top \state_t + q_t - B_t^\top \costate_{t+1}
  & = 0
  & t \in \{0, \ldots, \horizon -1\} \ (\partial_{\ctrl_t} \mathcal{L} = 0)\ \\
  x_{t+1} - A_t\state_t - B_t \ctrl_t
  & = 0
  & t \in \{0, \ldots, \horizon -1\} \ (\partial_{\costate_{t+1}} \mathcal{L} = 0)\ \\
  P_\horizon \state_\horizon + p_\horizon + \costate_\horizon
  & = 0
  & (\partial_{\state_\horizon} \mathcal{L} = 0).
\end{align*}
As noted by \citet{wright1991partitioned}, these equations can be ordered as
\begin{align*}
  x_0 
  & = \initstate 
  & (\partial_{\costate_0} \mathcal{L} = 0)\ \\
  \costate_0 + P_0 \state_0 + R_0 \ctrl_0 - A_0^\top \costate_1 
  & = - p_0 
  & (\partial_{\state_0} \mathcal{L} = 0)\ \\
  R_0^\top \state_0 + Q_0 \ctrl_0 - B_0^\top \costate_1
  &  = - q_0
  & (\partial_{\ctrl_0} \mathcal{L} = 0)\ \\
  - A_0 \state_0 - B_0 \ctrl_0 + \state_1
  & = 0
  & (\partial_{\costate_1} \mathcal{L} = 0)\ \\
  \costate_1 + P_1 \state_1 + R_1 \ctrl_1 - A_1^\top \costate_2 
  & = - p_1 
  & (\partial_{\state_1} \mathcal{L} = 0)\ \\
  R_1^\top \state_1 + Q_1 \ctrl_1 - B_1^\top \costate_2
  &  = - q_1
  & (\partial_{\ctrl_1} \mathcal{L} = 0)\ \\
  \vdots & & \\
  \costate_\horizon + P_\horizon \state_\horizon 
  & = - p_\horizon 
  & (\partial_{\state_\horizon} \mathcal{L} = 0).
\end{align*}
Written in matrix form the system to be solved is
\begin{align*}
  \begin{pmatrix}
    0 & I & & & & & & & & \\
    I & P_0 & R_0 & -A_0 & & & & & & \\
    & R_0^\top & Q_0 & -B_0^\top & & & & & & \\
    & -A_0 & - B_0 & 0 & I & & & & & \\
    & & & I & P_1 & R_1 & -A_1 & & & \\
    & & & & R_1^\top & Q_1 & -B_1^\top & & & \\
    & & & & -A_1 & -B_1 & 0 & \ddots & & \\
    & & & & & & \ddots & \ddots & & \\
    & & & & & & & & \ddots & I \\
    & & & & & & & & I & P_\horizon 
  \end{pmatrix}
  \begin{pmatrix}
    \costate_0 \\
    \state_0 \\
    \ctrl_0 \\
    \costate_1 \\
    \state_1 \\
    \ctrl_1 \\
    \costate_2 \\
    \vdots \\
    \costate_\horizon \\
    \state_\horizon
  \end{pmatrix}
  =
  \begin{pmatrix}
    - s_0 \\
    - p_0 \\
    - q_0 \\
    - s_1 \\
    - p_1 \\
    - q_1 \\
    - s_2 \\
    \vdots \\
    - s_\horizon \\
    - p_\horizon
  \end{pmatrix},
\end{align*}
where $s_0 = - \initstate$ and $s_t = 0$ are simply introduced for readability.

The system above is band block diagonal, which hints why it can be solved
efficiently by various methods. If all blocks were of size 1, that is,
$\dimstate = \dimctrl = 1$, the system would amount to a band diagonal matrix
$M$ with bandwidth $\sup\{|i-j|: M_{ij} > 0\} = 2$. Gaussian eliminations of
band-diagonal $n\times n$ matrices of bandwidth $k$ are well-known to have a
complexity of the order $O(n k^2)$. In our case, since the blocks are not of
size one, implementations of Gaussian elimination-like algorithms would incur an
$O(\dim_\state^3)$ or $O(\dim_\ctrl^3)$ to inverse each block.

\subsection{Riccati-Based Implementation}
~
{\bf Implementation.}
The system of equations presented above suggest some elimination
strategies~\citep{polak1971computational,wright1991partitioned}. For example,
the control variables $u_t$ can be eliminated from the system of equations as we
have 
\[
  \ctrl_t = -Q_t^{-1}R_t^\top \state_t - Q_t^{-1} q_t + Q_t^{-1} B_t^\top \costate_{t+1}. 
\]
After eliminating the control variables, the optimality conditions read
\begin{align*}
  x_0
  & = \initstate 
  & (\partial_{\costate_0} \mathcal{L} = 0)\ \\
  C_t\state_t 
  - D_t^\top \costate_{t+1} 
  + \costate_t
  & = c_t
  & t \in \{0, \ldots, \horizon -1\} \ (\partial_{\state_t} \mathcal{L} = 0)\ \\
  x_{t+1} 
  - D_t \state_t 
  - E_{t+1} \costate_{t+1} 
  & = e_{t+1}
  & t \in \{0, \ldots, \horizon -1\} \ (\partial_{\costate_{t+1}} \mathcal{L} = 0)\ \\
  P_\horizon \state_\horizon + \costate_\horizon
  & = - p_\horizon
  & (\partial_{\state_\horizon} \mathcal{L} = 0).
\end{align*}
for 
\begin{align*}
  C_t & = P_t - R_t Q_t^{-1}R_t^\top, \\
  D_t & = A_t + B_t Q_{t-1} R_t^\top, \\
  E_{t+1} & = B_t Q_t^{-1} B_t^\top, \\
  c_t & = - p_t + R_t Q_t^{-1} q_t,  \\
  e_{t+1} & = - B_t Q_t^{-1} q_t.
\end{align*}
The corresponding system of equations to solve is then band diagonal of the
following form, denoting $e_0 = \initstate, c_\horizon = - p_\horizon$,
\begin{align*}
  \begin{pmatrix}
    0 & I & & & & & \\
    I & C_0 & -D_0^\top & & & & \\
    & -D_0 & -E_1 & I & & & \\
    & & I & C_1 & \ddots & & \\
    & & & \ddots & \ddots & & \\
    & & & & & -E_\horizon & I \\
    & & & & & I & P_\horizon
  \end{pmatrix}
  \begin{pmatrix}
    \costate_0 \\
    \state_0 \\
    \costate_1 \\ 
    \state_1 \\
    \vdots \\
    \costate_\horizon \\
    \state_\horizon 
  \end{pmatrix}
  =
  \begin{pmatrix}
    e_0 \\
    c_0 \\
    e_1 \\
    c_1 \\
    \vdots \\
    e_\horizon \\
    c_\horizon
  \end{pmatrix}.
\end{align*}

We can show by induction that the Lagrange multipliers necessarily satisfy
\[
  \costate_t = F_t \state_t + f_t 
  \ \mbox{for all} \ t \in \{1, \ldots, \horizon-1\},
\]
for some matrices $F_t$ and vectors $f_t$. For $t= \horizon$, we already know
that $\costate_\horizon = - P_\horizon \state_\horizon - p_\horizon$. Assume the
property is true at time $t+1$, then 
\begin{align*}
\costate_{t+1} 
& = F_{t+1}(A_t \state_t + B_t\ctrl_t) + f_{t+1}
= F_{t+1}(A_t \state_t 
- B_t Q_t^{-1}R_t^\top \state_t 
- B_t Q_t^{-1} q_t 
+ B_tQ_t^{-1} B_t^\top \costate_{t+1}) + f_{t+1}.
\end{align*}
Rearranging the terms, we get that 
\begin{align*}
  \costate_{t+1} & = (I - F_{t+1}E_{t+1})^{-1}(F_{t+1}P_t x_t  + F_{t+1}e_{t+1} + f_{t+1}).
\end{align*}
Injecting this expression in the optimality conditions associated to $x_t$ (that
is the line $\partial_{\state_t} \mathcal{L} = 0$), we get
\[
  \costate_t 
  = (D_t^\top (I - F_{t+1}E_{t+1})^{-1}F_{t+1}P_t - C_t)\state_t
  + c_t + D_t^\top(I - F_{t+1}E_{t+1})^{-1}(F_{t+1}e_{t+1} + f_{t+1}).
\]
Hence, we can express $\costate_t = F_t \state_t + f_t$ with 
\begin{equation}\label{eq:backward_riccati}
  F_t = (D_t^\top (I - F_{t+1}E_{t+1})^{-1}F_{t+1}P_t - C_t), \qquad 
  f_t = c_t + D_t^\top(I - F_{t+1}E_{t+1})^{-1}(F_{t+1}e_{t+1} + f_{t+1}).
\end{equation}
Similarly, given $F_{t+1}, f_{t+1}$ such that $\costate_{t+1} =
F_{t+1}\state_{t+1} + f_{t+1}$, we can compute an expression of the optimal
$\state_{t+1}$ in terms of $\state_t$ from the optimality condition on
$\costate_{t+1}$. Namely, we have
\[
  \state_{t+1} - D_t \state_t - E_{t+1}F_{t+1} \state_{t+1} + E_{t+1}f_{t+1} = e_{t+1},
\]
and so 
\begin{equation}
  \state_{t+1} = (I - E_{t+1}F_{t+1})^{-1}(D_t\state_t - E_{t+1}f_{t+1} + e_{t+1}).
\label{eq:riccati}
\end{equation}
The whole resolution consists then in
\begin{enumerate}
  \item Computing $F_t, f_t$ from $t=\horizon$ to $0$ using Eq.~\eqref{eq:backward_riccati} starting from $F_\horizon= P_\horizon$, $f_\horizon =p_\horizon$. 
  \item Computing the optimal $\state_0, \ldots, \state_\horizon$ starting from
  $\state_0 = \initstate$ and using Eq.~\eqref{eq:riccati} from $t = 0, \ldots,
  \horizon-1$.
\end{enumerate}

{\bf Computational Complexity.}
Compared to the implementation by dynamic programming, we retrieve a linear
complexity with respect to the horizon $\tau$ (only two passes on the dynamics),
and cubic in the control and state dimensions. One finds that the computational complexity of
the method presented above, taking into account the symmetry of some matrices, \citep{wright1991partitioned}
is of the order
\[
  \horizon(7 \dimstate^3 + 4\dimstate^2 \dimctrl + 4\dimstate \dimctrl^2 + \frac{1}{3} \dimctrl^3).
\]
In comparison, the computational
complexity of a dynamic programming-based approach is \citep{wright1991partitioned}
\[
  \horizon(3\dimstate^3 + 5\dimstate^2 \dimctrl + 3 \dimstate \dimctrl^3 + \frac{1}{3} \dimctrl^3) 
  + O(\horizon(\dimstate^2 + \dimctrl^2)).
\]
While the method presented in this section may be slightly more computationally
expansive than a dynamic programming approach, it may be easier to use in a
parallel context as recalled below.

\subsection{Parallel Implementation}
Rather than eliminating the set of control variables, one can consider
eliminating blocks of variables to enable parallel implementations of such
methods as presented by \citet{wright1991partitioned}. Briefly, the approach
consists in considering a system reduced to the variables at $L+1$ time steps,
i.e., $(\costate_{t_i}, \state_{t_i}, \ctrl_{t_i})_{i=0}^L$ for $t_0 = 0$ and
$t_L = \horizon$. Intermediate variables between time-steps, that is
$(\costate_{t_i+j}, \state_{t_i + j}, \ctrl_{t_i+ j})_{j=1}^{t_{i+1}-1}$ are
eliminated by appropriate computations to reduce the system as a set of $3(P+1)
- 1$ equations, akin to the original system, 
\begin{align*}
  \begin{pmatrix}
    0 & I & & & & & & & & \\
    I & \tilde P_0 & \tilde R_0 & -\tilde A_0 & & & & & & \\
    & \tilde R_0^\top & \tilde Q_0 & -\tilde B_0^\top & & & & & & \\
    & -\tilde A_0 & - \tilde B_0 & 0 & I & & & & & \\
    & & & I & \tilde P_1 & \tilde R_1 & -\tilde A_1 & & & \\
    & & & & \tilde R_1^\top & \tilde Q_1 & -\tilde B_1^\top & & & \\
    & & & & -\tilde A_1 & - \tilde B_1 & 0 & \ddots & & \\
    & & & & & & \ddots & \ddots & & \\
    & & & & & & & & \ddots & I \\
    & & & & & & & & I & \tilde P_L 
  \end{pmatrix}
  \begin{pmatrix}
    \costate_{t_0} \\
    \state_{t_0} \\
    \ctrl_{t_0} \\
    \costate_{t_1} \\
    \state_{t_1} \\
    \ctrl_{t_1} \\
    \costate_{t_2} \\
    \vdots \\
    \costate_{t_L} \\
    \state_{t_L}
  \end{pmatrix}
  =
  \begin{pmatrix}
    - \tilde s_{t_0} \\
    - \tilde p_{t_0} \\
    - \tilde q_{t_0} \\
    - \tilde s_{t_1} \\
    - \tilde p_{t_1} \\
    - \tilde q_{t_1} \\
    - \tilde s_{t_2} \\
    \vdots \\
    - \tilde s_{t_L} \\
    - \tilde p_{t_L}
  \end{pmatrix},
\end{align*}
The matrices $\tilde M_j$ for $M \in \{A, B, P, Q, R\}$ can be computed 
as functions of the intermediate results at that stage, that is a function of 
$M_{t_j +1}, \ldots,  M_{t_j -1}$ for $M \in \{A, B, P, Q, R\}$, see 
\citet{wright1991partitioned} for detailed expressions.
Solving the reduced system above is naturally less computationally expensive
than computing the original system, while the computations of the reduced system,
that is, the computations of $\tilde M_j$ for $M \in \{A, B, P, Q, R\}$ can be
done in parallel.

\subsection{Matrix-free Solver}
Finally, rather than considering computing Newton or Gauss-Newton steps by
exploiting the structure of the problem, one can directly use the access to
hessian-vector products in a differentiable programming framework.\\

{\bf Implementation.}
Consider the case of a Newton step, which requires computing
\[
\nabla^2 \obj(\ctrls)^{-1} \nabla \obj(\ctrls),
\]
for $\obj$ the objective defined in Section~\ref{sec:classical_optim}. 
Rather than computing the Hessian, and inverting it, this oracle can be
computed by solving for $\auxctrls$ such that 
\[
\nabla^2 \obj(\ctrls) \auxctrls = \nabla \obj(\ctrls),
\]
which can be done approximately by an iterative method such as a
conjugate gradient method or a generalized minimal residual
method~\citep{nocedal2006numerical}, provided 
that we have access only to the linear operator 
$\auxctrls \mapsto \nabla^2 \obj(\ctrls) \auxctrls$. This can be done efficiently in a 
differentiable programming framework as recalled below.

Automatic differentiation naturally gives access to the gradient 
$\nabla \obj(\ctrls)$ of the objective at some inputs
$\ctrls$ by means of the reverse mode of automatic differentiation. For a
function $g: \reals^n \rightarrow \reals^m$, its directional derivative at
$\ctrls$ along a direction $\auxctrls$, that is the derivative of $t\mapsto
g(\ctrls + t\auxctrls)$ at $t=0$, denoted $\partial g(\ctrls)[\auxctrls]$, can
be computed by forward mode automatic differentiation. The linear operator
$\auxctrls \mapsto \nabla^2 \obj(\ctrls) \auxctrls$ amounts to the directional
derivative of the gradient, that is, $ \nabla^2 \obj(\ctrls) \auxctrls =
\partial (\nabla \obj)(\ctrls)[\auxctrls]$. It can then be computed by forward
mode automatic differentiation on top of reverse mode automatic differentiation
at approximately twice the computational cost of the gradient~\citep{griewank2008evaluating}.

{\bf Computational Complexity.}
The overall computational cost of computing the Newton oracle depends then on the 
condition number of the Hessian 
$\kappa = \sigma_{\max}(\nabla^2 \obj(\ctrls))/\sigma_{\min}(\nabla^2\obj(\ctrls))$ as
\[
O\left(
\min
\left\{
\horizon \dimctrl,
\frac{\sqrt{\kappa} - 1}{\sqrt{\kappa} + 1}\log(\varepsilon^{-1})
\right\}
\mathcal{T}(\nabla^2 \obj(\ctrls)))\right) 
= 
O\left(
\min
\left\{
\horizon \dimctrl,
\frac{\sqrt{\kappa} - 1}{\sqrt{\kappa} + 1}\log(\varepsilon^{-1})
\right\}
\horizon(\dimstate^2 + \dimctrl^2)\right)
\] 
operations, where $\mathcal{T}(\nabla^2 \obj(\ctrls)))$ denotes the cost of computing
the Hessian-vector product $\auxctrls \mapsto \nabla^2 \obj(\ctrls))\auxctrls$ in 
a differentiable 
programming framework and can be approximated roughly as 
$
\mathcal{T}(\nabla^2 \obj(\ctrls))) = O(\horizon(\dimstate^2 + \dimctrl^2)).
$
Overall, such ``matrix-free'' methods, which circumvent the need to compute actually 
the Hessian, have a priori a quadratic complexity and not linear complexity w.r.t. the horizon 
if the matrix is ill-conditioned.
On the other hand, the complexity of such methods remains quadratic in the state
dimension. 

A similar approach can be used to compute Gauss-Newton steps by using the
Jacobian vector product of the function that at controls associate the
trajectory. For problems with a single final cost, Gauss-Newton methods can also
benefit from considering their dual formulation as shown
by~\citet{roulet2019iterative}.

We present in Appendix~\ref{app:exp_sup} numerical comparisons of such matrix-free
implementations to the implementation by dynamic programming. Note that by using matrix-free
solvers in a differentiable programming framework we can cast any nonlinear control as a generic 
numerical optimization problem amenable to solutions with off-the-shelf programs such as 
IPOPT~\citep{wachter2006implementation}.

  \section{Experimental Detail}\label{app:exp_details}
  We describe in detail the continuous time systems studied in the experiments.
The code is available at {\small \coderef}. Numerical constants are detailed at
the end for reference. All algorithms are run with double precision.

\subsection{Discretization}
In the following, we denote by $\statexp(t)$ the state of a system at time $t$.
Given a control $\ctrl(t)$ at time $t$, we consider time-invariant dynamical
systems governed by a differential equation of the form
\[
\dot \statexp(t) = \contdyn(\statexp(t), \ctrl(t)), \quad \mbox{for} \ t \in [0, T],
\]
where $\contdyn$ models the physics of the movement and is described below for
each model. 

Given a continuous time dynamic, the discrete time dynamics are given by a
discretization method such that the states follow dynamics of the form 
\[
\statexp_{t+1} = \dyn(\statexp_t, \ctrl_t) \quad \mbox{for} \ t\in \{0, \ldots \horizon-1\},
\]
for a sequence of controls $\ctrl_0, \ldots, \ctrl_{\horizon-1}$. One
discretization method is the Euler method, which, for a time-step $\Delta =
T/\horizon$, is  
\[
\dyn(\statexp_t, \ctrl_t)  = \statexp_t + \Delta \contdyn(\statexp_t, \ctrl_t).
\]
Alternatively, we can consider a Runge-Kutta method of order 4 that defines the
discrete-time dynamics as 
\begin{align*}
	\dyn(\statexp_t, \ctrl_t) & = \statexp_t + \frac{\Delta}{6} (k_1+k_2+k_3+k_4) \\
	\mbox{where} \quad k_1 & = \contdyn(\statexp_t, \ctrl_t) \hspace{60pt}
	k_2   = \contdyn(\statexp_t + \Delta k_1/2, \ctrl_t) \\
	k_3 & = \contdyn(\statexp_t + \Delta k_2/2, \ctrl_t) \qquad
	k_4  = \contdyn(\statexp_t + \Delta k_3, \ctrl_t),
\end{align*}
where we consider the controls to be piecewise constant, i.e., constant on time
intervals of size $\Delta$. We can also consider a Runge-Kutta method with
varying control inputs such that, for $\ctrl_t = (\auxctrl_t, \auxctrl_{t+1/3},
\auxctrl_{t+2/3})$,
\begin{align*}
	\dyn(\statexp_t, \ctrl_t) & = \statexp_t + \frac{\Delta}{6} (k_1+k_2+k_3+k_4) \\
	\mbox{where} \quad k_1 & = \contdyn(\statexp_t, \auxctrl_t) \hspace{78pt}
	k_2   = \contdyn(\statexp_t + \Delta k_1/2, \auxctrl_{t+1/3}) \\
	k_3 & = \contdyn(\statexp_t + \Delta k_2/2, \auxctrl_{t+1/3}) \qquad
	k_4  = \contdyn(\statexp_t + \Delta k_3, \auxctrl_{t+2/3}).
\end{align*}

\subsection{Swinging up a Pendulum}

\subsubsection{Fixed Pendulum}
\begin{figure}[t]
	\begin{minipage}{0.5\linewidth}
		\begin{center}
			\includegraphics[width=0.2\linewidth]{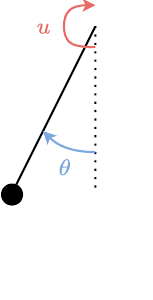}
			\caption{Fixed pendulum.\label{fig:pendulum}}
		\end{center}
	\end{minipage}~\hfill~
	\begin{minipage}{0.5\linewidth}
		\begin{center}
			\includegraphics[width=0.42\linewidth]{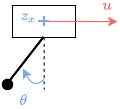}
			\caption{Pendulum on a cart. \label{fig:cart_pendulum}}
		\end{center}
	\end{minipage}
\end{figure}

We consider the problem of controlling a fixed pendulum such that it swings up
as illustrated in Fig.~\ref{fig:pendulum}. Namely, the dynamics of a pendulum
are given as
\begin{align*}
	ml^ 2\ddot \theta(t) & = -mlg \sin \theta(t) -\mu \dot \theta(t) + \ctrl(t),
\end{align*}
with $\theta$ the angle of the rod, $m$ the mass of the blob, $l$ the length of
the blob, $\mu$ a friction coefficient, $g$ the gravitational constant, and
$\ctrl$ a  torque applied to the pendulum (which defines the control we have on
the system). Denoting the angle speed $\omega = \dot \theta$ and the state of
the system $\state = (\theta; \omega)$, the continuous time dynamics are 
\[
\contdyn	:
	(	\state=(\theta; \omega), \ctrl)\rightarrow \left(\begin{matrix}
		\omega \\
		-\frac{g}{l} \sin \theta -\frac{\mu}{m l^2} \omega + \frac{1}{ml^2} \ctrl
	\end{matrix}\right),
\]
such that the continuous time system is defined by  $\dot \state(t) =
\contdyn(\state(t), \ctrl(t))$. After discretization by an Euler method, we get
discrete  time dynamics $\dyn_t(\state_t, \ctrl_t) = \dyn(\state_t, \ctrl_t) $
of the form,  for $\state_{t} = (\theta_{t}; \omega_{t})$ and  $\Delta$  the
discretization step,
\begin{align*}
	\dyn(\state_t, \ctrl_t)= \state_t + \Delta \contdyn(\state_t, \ctrl_t) = \left(\begin{matrix}
		\theta_t + \Delta \omega_t  \\
		\omega_t + \Delta \left( -\frac{g}{l} \sin \theta_t - \frac{\mu}{ml^2} \omega_t + \frac{1}{ml^2}u_t\right)
	\end{matrix}\right).
\end{align*}

A classical task is  to enforce the pendulum to swing up and stop without using
too much torque at each time step, i.e., for $\initstate = (0; 0)$, the costs we
consider are, for some non-negative parameters $\lambda \geq 0, \rho \geq 0$,
\[
\cost_t(\state_t, \ctrl_t) = \lambda\|\ctrl_t\|_2^2 \quad \mbox{for $t \in \{0, \ldots, \horizon-1\}$},\quad  \cost_\horizon(\state_\horizon) = (\pi- \theta_\horizon)^ 2 + \rho \|\omega_\horizon\|_2^ 2.
\]

\subsubsection{Pendulum on a Cart}
We consider here controlling a pendulum on a cart as illustrated in
Fig.~\ref{fig:cart_pendulum}. This system is described by the angle $\theta$ of
the pendulum with the vertical and the position $z_x$ of the cart on the
horizontal axis. Contrary to the previous example, here we do not control
directly the angle of the pendulum we only control the system with a force $u$
that drives the acceleration of the cart. The dynamics of the system satisfy
(see~\citet{magdy2019modeling} for detailed derivations)
\begin{align}
		(M+m)\ddot z_x +  ml \cos \theta  \ddot \theta &  = -b\dot z_x +ml\dot\theta^2 \sin \theta + u \nonumber \\
	ml \cos \theta \ddot z_x + (I + ml^2) \ddot \theta   &= -mgl \sin \theta, \label{eq:cart_pendulum} 
\end{align}
where  $M$ is the mass of the cart, $m$ is the mass of the pendulum rod, $I$ is
the pendulum rod moment of inertia, $l$ is the length of the rod, and $b$ is the
viscous friction coefficient of the cart. The system of equations can be written
in matrix form and solved to express the angle and position accelerations as 
\begin{align*}
\left(\begin{matrix}
	\ddot z_x\\
\ddot \theta 
\end{matrix}\right) 
& = \left(\begin{matrix}
	M+m & ml \cos \theta \\
	ml \cos \theta & I  + ml^2 
\end{matrix}\right)^{-1} 
\left(\begin{matrix}
	-b \dot z_x + ml\dot \theta^2 \sin \theta + u\\
	-mgl\sin \theta 
\end{matrix}\right) \\
& = \frac{1}{I(M+m) + ml^2 M + m^2 l^2\sin^2 \theta}
\left(
\begin{matrix}
I  + ml^2  &  -ml \cos \theta \\
- ml \cos \theta & 	M+m
\end{matrix}\right)
\left(\begin{matrix}
	-b \dot z_x + ml\dot \theta^2 \sin \theta + u\\
	-mgl\sin \theta 
\end{matrix}\right).
\end{align*}
The discrete dynamical system follows using an Euler discretization scheme or a
Runge Kutta method. We consider the task of swinging up the pendulum and keeping
it vertical for a few time steps while constraining the movement of the cart on
the horizontal line. Formally, we consider the following cost, defined for $x_t
= (z_x, \theta, \zeta_x, \omega)$,  where $\zeta_x, \omega$ represent the
discretizations of $\dot z_x$ and $\dot \theta$ respectively, 
\revised{
\begin{align}
	h(x_t, u_t) &  = \begin{cases}
		\rho_2 (\max((z_x - \bar z_x^+)^3, 0 ) + \max((z_x + \bar z_x^-)^3, 0 ))  + \lambda u_t^2 + (\theta + \pi)^2 + \rho_1 \omega^2 & \mbox{if} \ t \geq \bar t\\
		\rho_2 (\max((z_x - \bar z_x^+)^3, 0 ) + \max((z_x + \bar z_x^-)^3, 0 ))  + \lambda u_t^2 & \mbox{if} \ t < \bar t,
	\end{cases} \nonumber 
\end{align}
}
where $\rho_1, \rho_2, \lambda$ are some non-negative parameters, $\bar t$ is a
time step after which the pendulum needs to stay vertically inverted and $\bar
z_x^+, \bar z_x^-$ are bounds that restrain the movement of the cart along the
whole horizontal line.

\subsection{Autonomous Car Racing}
We consider the control of a car on a track through two different dynamical
models: a simple one where the orientation of the car is directly controlled by
the steering angle, and a more realistic one that takes into account the tire
forces to control the orientation of the car. In the following, we present  the
dynamics,  a simple tracking cost, and  a contouring cost enforcing the car to
race the track at a reference speed or  as fast as possible. 

\begin{figure}[t]
\begin{minipage}{0.5\linewidth}
	\begin{center}
	\vspace{65pt}	
	\includegraphics[width=0.6\linewidth]{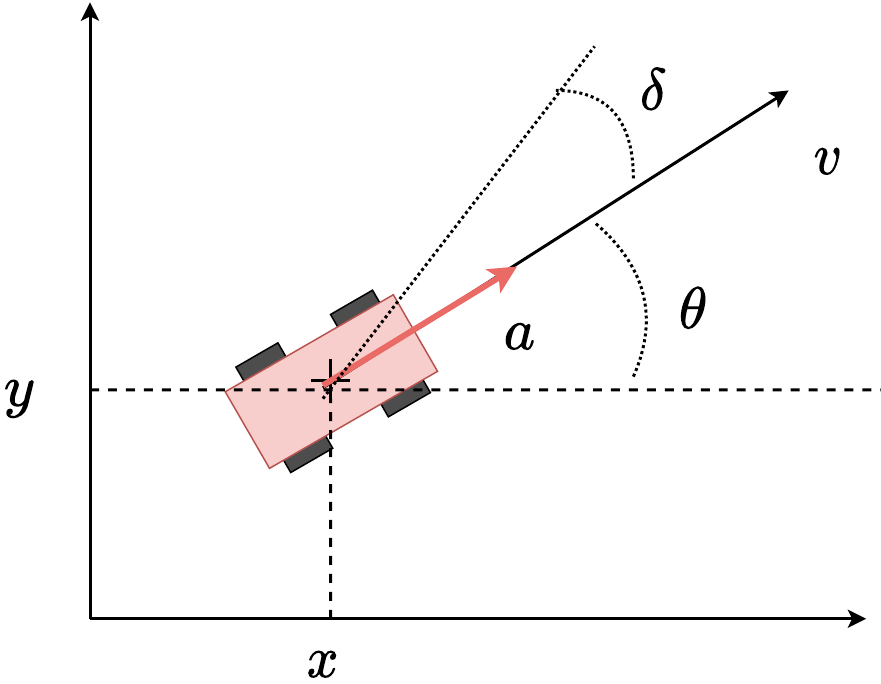}
	\caption{Simple model of a car.\label{fig:simple_model}}
	\end{center}
\end{minipage}~\hfill~
\begin{minipage}{0.5\linewidth}
\begin{center}
	\includegraphics[width=\linewidth]{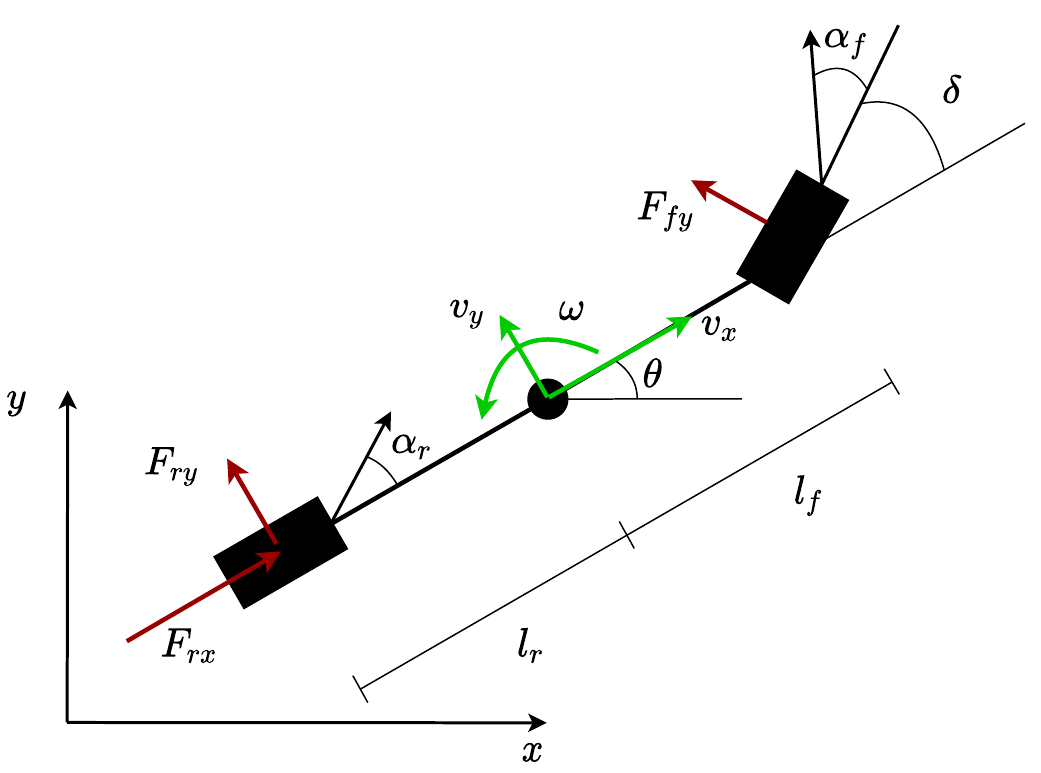}
	\caption{Bicycle model of a car.\label{fig:bicycle_model}}
\end{center}
\end{minipage}
\end{figure}

\subsubsection{Dynamic}

\paragraph{Simple Model}
A simple model of the car is described in Fig.~\ref{fig:simple_model}. The state
of the car is decomposed as $\statexp(t)  = (x(t), y(t), \yaw(t), v(t))$, where
(dropping the dependency w.r.t. time for simplicity)
\begin{enumerate}
	\item $x, y$ denote the position of the car on the plane,
	\item $\yaw$ denotes the angle between the orientation of the car and the
	horizontal axis, a.k.a. the yaw,
	\item $v$ denotes the longitudinal speed.
\end{enumerate}
\vspace{5pt}
The car is controlled through $\ctrl(t) = (a(t), \delta(t))$, where
\begin{enumerate}
	\item $a$ is the longitudinal acceleration of the car, 
	\item $\delta$ is the steering angle.
\end{enumerate}
For a car of length $L$, the continuous time dynamics are then 
\begin{align}\label{eq:simple_car}
	\dot x  = v \cos \yaw  \qquad 	
	\dot y   = v \sin\yaw \qquad  
	\dot \yaw   = v \tan(\delta)/L  \qquad
 	\dot v  =  a.
\end{align}

\paragraph{Bicycle Model}

\begin{figure}[t]
	\begin{center}
		\includegraphics[width=0.4\linewidth]{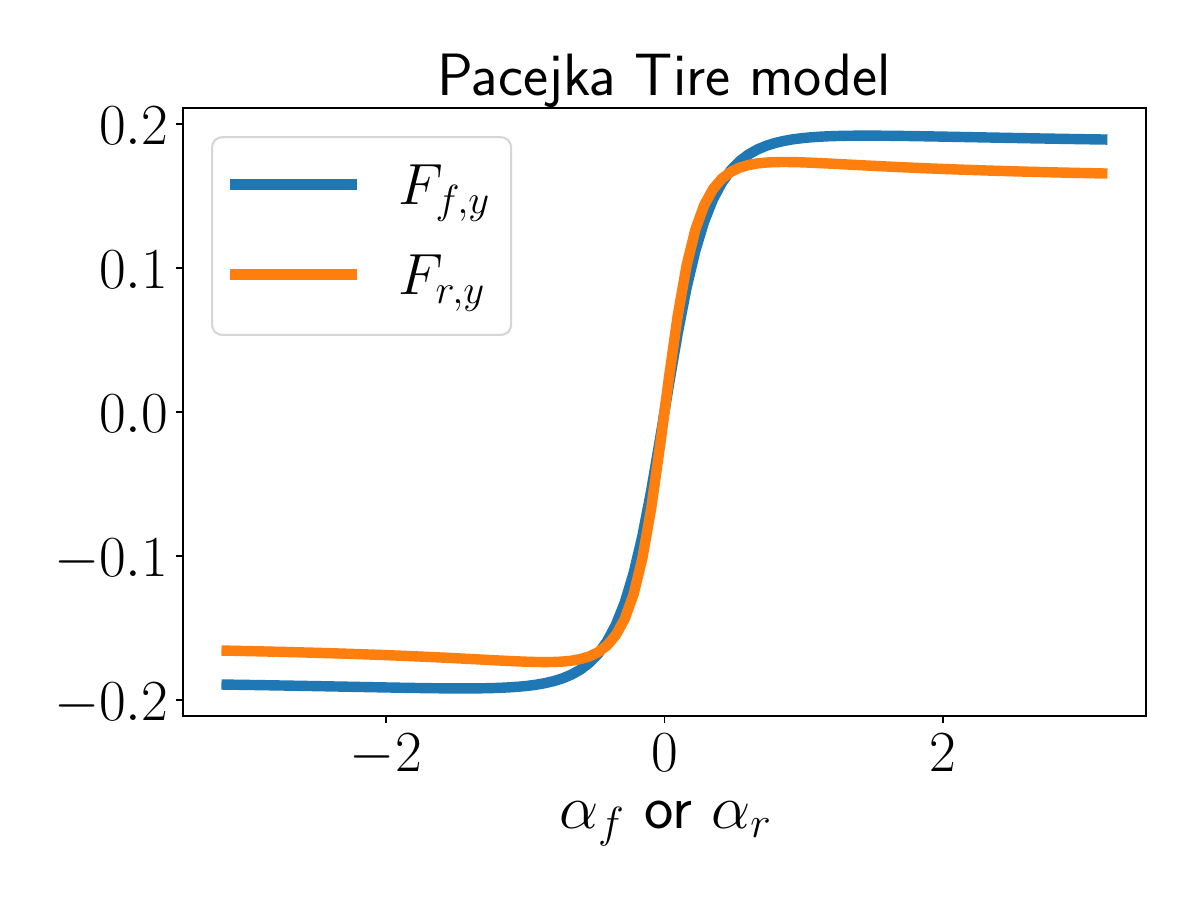}
	\end{center}
	\caption{Pacejka model of the friction on the tires  as a function of the slip
	angles\label{fig:pacejka_model}}
\end{figure}

We consider the model presented by \citet{liniger2015optimization} recalled
below and illustrated in Fig.~\ref{fig:bicycle_model}. In this model, the state
of the car at time $t$ is decomposed as $\statexp(t) = (x(t), y(t), \yaw(t),
v_x(t), v_y(t), \omega(t))$ where
\begin{enumerate}
	\item $x, y$ denote the position of the car on the plane,
	\item $\yaw$ denotes the angle between the orientation of the car and the
	horizontal axis, a.k.a. the yaw,
	\item $v_x$ denotes the longitudinal speed, 
	\item $v_y$ denotes the lateral speed, 
	\item $\omega$ denotes the derivative of the orientation of the car, a.k.a.
	the yaw rate.
\end{enumerate}
The control variables are analogous to the simple model, i.e., $\ctrl(t) =
(a(t), \delta(t))$, where
\begin{enumerate}
	\item $a$ is the PWM duty cycle of the car, this duty cycle can be negative to
	take into account braking,
	\item $\delta$ is the steering angle.
\end{enumerate}
These controls act on the state through the following forces.
\begin{enumerate}
\item A longitudinal force on the rear wheels, denoted $F_{r, x}$ modeled using
a motor model for the DC electric motor as well as a friction model for the
rolling resistance and the drag
$$
F_{r, x} = (C_{m1} -C_{m2} v_x) a - C_{r0} - C_{rd} v_x^2,
$$
where $C_{m1}, C_{m2}, C_{r0}, C_{rd}$ are constants estimated from experiments,
see Appendix~\ref{app:exp_details}. 
\item Lateral forces on the front and rear wheels, denoted $F_{f, y}, F_{y, r}$
respectively, modeled using a simplified Pacejka tire model 

\begin{align*}
	F_{f, y} & = D_f \sin(C_f \arctan(B_f \alpha_f)) \quad \textrm{where} \ \alpha_f = \delta - \operatorname{arctan2}\left(\frac{\omega l_f + v_y}{v_x}\right)  \\
	F_{r, y} & = D_r \sin(C_r \arctan(B_r \alpha_r)) \quad \textrm{where} \ \alpha_r = \operatorname{arctan2}\left(\frac{\omega l_r - v_y}{v_x}\right)  
\end{align*}

where $\alpha_f$, $\alpha_r$ are the slip angles on the front and rear wheels
respectively, $l_f, l_r$ are the distance from the center of gravity to the
front and the rear wheel respectively and the constants $B_r, C_r, D_r, B_f,
C_f, D_f$ define the exact shape of the semi-empirical curve, presented in
Fig.~\ref{fig:pacejka_model}.
\end{enumerate}

The continuous time dynamics are then 
\begin{align}
	\dot x & = v_x \cos \yaw - v_y\sin \yaw & \dot v_x & =\frac{1}{m}(F_{r,x} - F_{f, y} \sin \delta) + v_y \omega \label{eq:real_dyn}\\
	\dot y & = v_x\sin \yaw  + v_y\cos \yaw  & \dot v_y & = \frac{1}{m}(F_{r, y} + F_{f, y} \cos\delta) - v_x \omega  \nonumber\\
	\dot \yaw & = \omega & \dot \omega & = \frac{1}{I_z}(F_{f, y} l_f \cos\delta - F_{r, y} l_r),\nonumber
\end{align}
where $m$ is the mass of the car and $I_z$ is the inertia.

\subsubsection{Cost}
\begin{figure}[t]
	\begin{center}
		\includegraphics[width=0.4\linewidth]{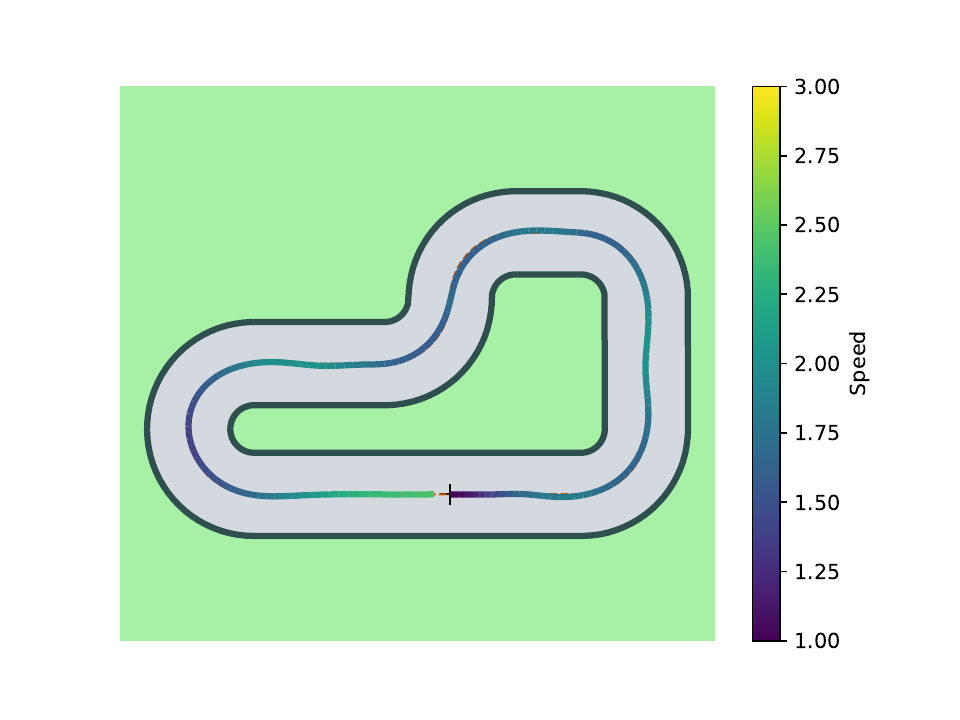}~
		\includegraphics[width=0.4\linewidth]{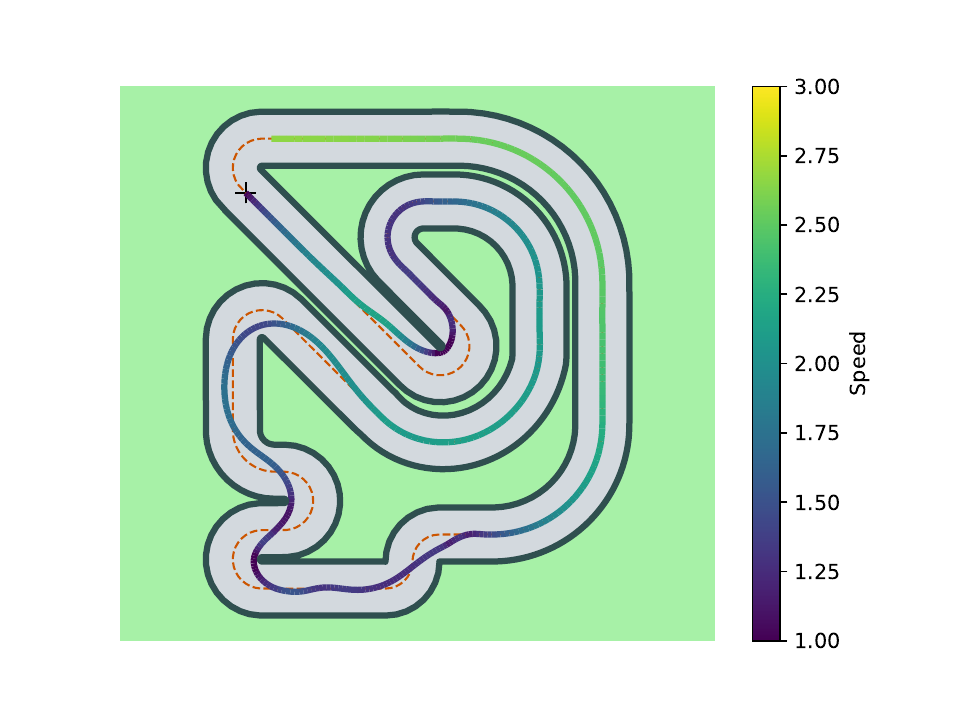}
		\caption{Simple and complex tracks used with  a trajectory computed on the
		bicycle model~\eqref{eq:real_dyn}.\label{fig:simple_track}}
	\end{center}
\end{figure}

\paragraph{Track}
We consider tracks that are given as a continuous curve, namely a cubic spline
approximating a set of points. As a result, for any time $t$, we have access to
the corresponding point $\hat x(t), \hat y(t)$ on the curve. The track we
consider is a simple track illustrated in Fig.~\ref{fig:simple_track}.

\paragraph{Tracking Cost}
A simple cost on the states is 
\begin{equation}\label{eq:tracking_cost}
\statecost_t(\statexp_t) = \|x_t - \hat x(\Delta v^{\textrm{ref}} t)\|_2^2 + \|y_t - \hat y(\Delta v^{\textrm{ref}}t)\|_2^2\quad \mbox{for} \ t =1, \ldots, \horizon,
\end{equation}
for $\statexp_t = (x_t, y_t)$, where $\Delta$ is some discretization step and
$v^{\textrm{ref}}$ is some reference speed. The cost above is the one we choose
for the simple model of a car. The disadvantage of such a cost is that it
enforces the car to follow the track at a constant speed which may not be
physically possible. We consider in the following a contouring cost as done by
\citet{liniger2015optimization}.

\begin{figure}[t]
	\begin{minipage}{0.33\linewidth}
		\begin{center}
		\includegraphics[width=\linewidth]{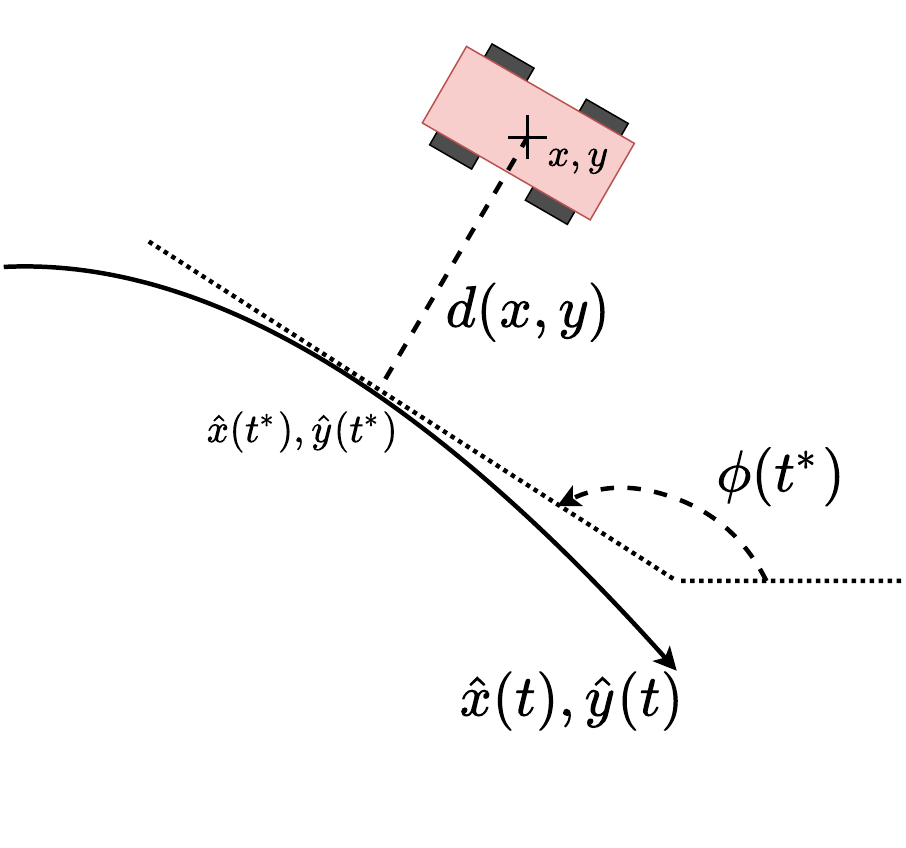}
		\caption{Distance to the track.\label{fig:true_dist}}
		\end{center}
	\end{minipage}~
	\begin{minipage}{0.33\linewidth}
		\begin{center}
		
		\vspace{40pt}
		\includegraphics[width=\linewidth]{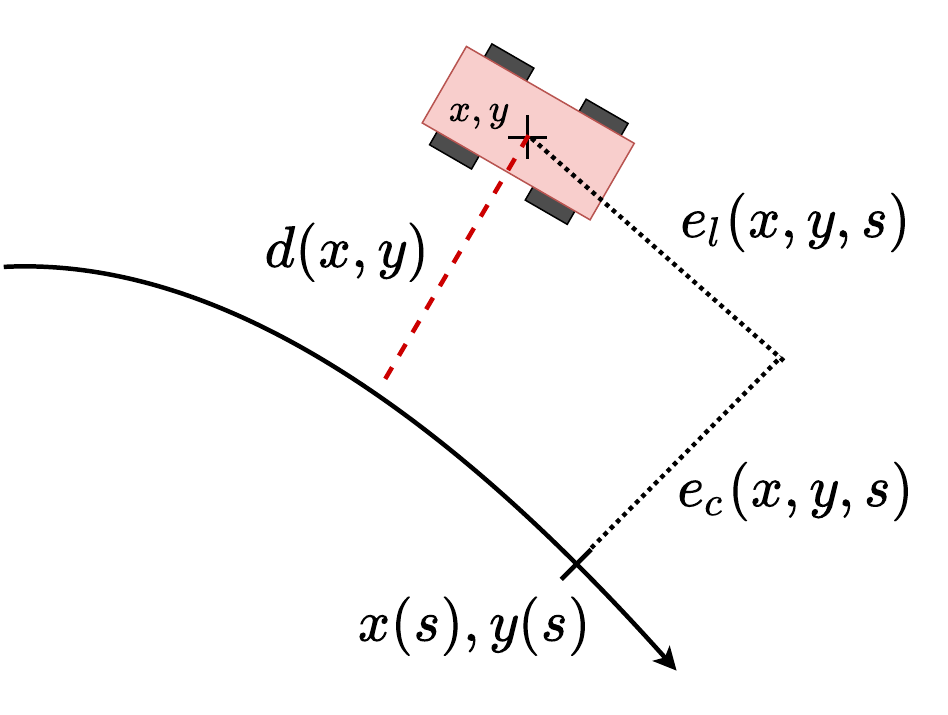}
		\caption{Approx. by contouring and lagging costs.
		\label{fig:contouring_cost_border_costs}}
		\end{center}	
	\end{minipage}~
	\begin{minipage}{0.33\linewidth}
		\begin{center}
		\includegraphics[width=0.9\linewidth]{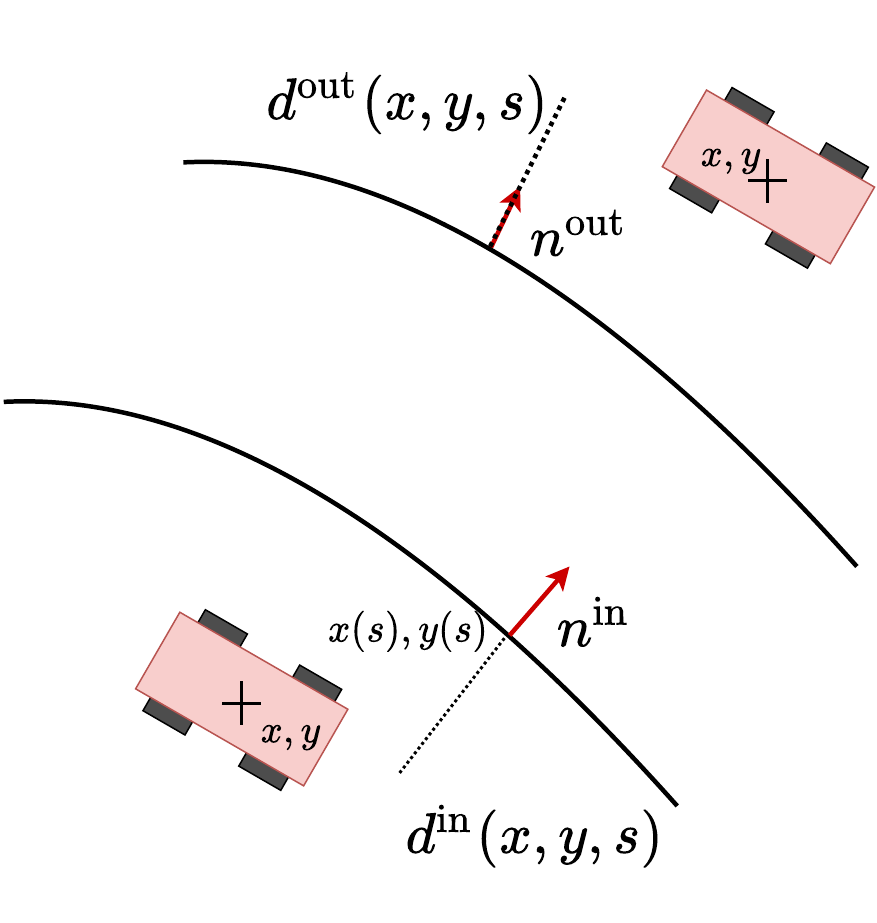}
		\caption{Border costs.\label{fig:border}}
		\end{center}
	\end{minipage}
\end{figure}

\paragraph{Ideal Cost}
Given a track  parameterized in continuous time, an ideal cost is to enforce the
car to be as close as possible to the track, while moving along the track as
fast as possible. Formally, define the distance from the car at position $(x,
y)$ to the track defined by the curve $\hat x(t), \hat y(t)$ as 
$$d(x, y) = \min_{t\in \reals}\quad \sqrt{((x - \hat x(t))^2 + (y - \hat
y(t))^2}.$$ Denoting $t^* = t(x, y)  = \argmin_{t\in \reals}\quad (x - \hat
x(t))^2 + (y - \hat y(t))^2,$ the reference time on the track for a car at
position $(x, y)$, the distance $d(x, y)$ can be expressed as 
\[d(x, y) = \sin(\yaw(t^*))\left(x - \hat x(t^*)\right) - \cos(\yaw(t^*))
\left(y - \hat y(t^*)\right),\] where  $\yaw(t) = \frac{\partial \hat
y(t)}{\partial \hat x(t)}$ is the angle of the track with the x-axis. The
distance $d(x, y)$ is illustrated in Fig.~\ref{fig:true_dist}. An ideal cost for
the problem is then defined as $ h(\statexp) = h(x, y) = d(x, y)^2 -t(x, y),$
which enforces the car to be close to the track by minimizing $d(x, y)^2$, and
also encourages the car to go as far as possible by adding the term $-t(x, y)$.

\paragraph{Contouring and Lagging Cost}
The computation of $t^*$ involves solving an optimization problem and is not
practical. As \citet{liniger2015optimization}, we rather augment the states with
a flexible reference time. Namely, we augment the state of the car by adding a
variable $s$ whose objective is to approximate the reference time $t^*$. The
cost is then decomposed into the \emph{contouring cost} and the \emph{lagging
cost} illustrated in Fig.~\ref{fig:contouring_cost_border_costs} and defined  as 
\begin{align*}
	e_c(x, y, s) & = \sin(\yaw(s))\left(x - \hat x(s)\right) - \cos(\yaw(s)) \left(y - \hat y(s)\right) \\
	e_l(x, y, s) & = -\cos(\yaw(s))\left(x - \hat x(s)\right) - \sin(\yaw(s)) \left(y - \hat y(s)\right).
\end{align*}

Rather than encouraging the car to make the most progress on the track, we
enforce them to keep a reference speed. Namely, we consider an additional penalty
of the form $\|\dot s - v^{\textrm{ref}}\|_2^2$ where $v^{\textrm{ref}}$ is a
parameter chosen in advance. For the reference time $s$ not to go backward in
time, we add a log-barrier term $-\varepsilon\log(\dot s)$ for
$\varepsilon=10^{-6}$.

Finally, we let the system control the reference time through its second order
derivative $\ddot s$. Overall this means that we augment the state variable by
adding the variables $s$ and $\nu := v_s$ and that we augment the control
variable by adding the variable $\alpha := a_s$ such that the discretized
problem is written for, e.g., the bicycle model, as
\begin{align*}
	\min_{\substack{(a_0, \delta_0, \alpha_0),  \ldots,  (a_{\tau-1}, \delta_{\tau-1},\alpha_{\tau-1}}) } \quad & 
	\sum_{t=0}^{\tau-1} 
	\rho_c e_c(x_t, y_t, s_t)^2 
	+ \rho_l e_l(x_t, y_t, s_t)^2  
	+ \rho_v \|v_{s, t} - v^{\textrm{ref}}\|_2^2 
	- \varepsilon \log \nu_t\\
	\textrm{s.t.} \quad & x_{t+1}, y_{t+1}, \yaw_{t+1}, v_{x, t+1}, v_{y, t+1}, \omega_{t+1} = \dyn(x_t, y_t, \yaw_, v_{x, t}, v_{y, t}, \omega_t, \delta_t, a_t) \\
	&  s_{t+1} = s_t + \Delta \nu_{t}, \quad   \nu_{t+1} = \nu_{t} + \Delta \alpha_{t} \\
	&   z_0 = \hat z_0  \quad s_0 =0 \quad \nu_0 = v^{\textrm{ref}},
\end{align*}
where $\dyn$ is a discretization of the continuous time dynamics, $\Delta$ is a
discretization step and $\hat z_0$ is a given initial state where $z_0$ regroups
all state variables at time 0 (i.e. all variables except $a_0, \delta_0$).

This cost is defined by the parameters $\rho_c, \rho_l, \rho_v,
v^{\textrm{ref}}$ which are fixed in advance. The larger the parameter $\rho_c$,
the closer the car to the track. The larger the parameter $\rho_l$, the closer
the car to its reference time $s$. In practice, we want the reference time to be
a good approximation of the ideal projection of the car on the track so $\rho_l$
should be chosen large enough. On the other hand, varying $\rho_c$ allows having
a car that is either conservative and potentially slow or a car that is fast but
inaccurate, i.e., far from the track. The most important aspect of the
trajectory is to ensure that the car remains inside the borders of the track
defined in advance. 

\paragraph{Border Cost}
To enforce the car to remain inside the track defined by some borders, we
penalize the approximated distance of the car to the border when it goes outside
the border as $e_b(x, y, s) = e_b^\textrm{in}(x, y, s) + e_b^\textrm{out}(x, y,
s)$ with
\revised{
\begin{align}\label{eq:border_costs}
	e_b^\textrm{in}(x, y, s) & =  \max((w + d^{\textrm{in}}(x, y, s))^3, 0) & d^{\textrm{in}}(x, y, s) & =  -(z - z^\textrm{in}(s))^\top n^{\textrm{in}}(s) \\
	e_b^\textrm{out}(x, y, s) & =  \max((w + d^{\textrm{out}}(x, y, s))^3, 0) & d^{\textrm{out}}(x, y, s)  & = (z - z^\textrm{out}(s))^\top n^{\textrm{out}}(s)\nonumber
\end{align}
}
for $\statexp = (x, y)$, where $n^{\textrm{in}}(s)$ and $n^{\textrm{out}}(s)$
denote the normal at the borders at time $s$ and $w$ is the width of the car. In
practice, we use a smooth approximation of the max function in
Eq.~\eqref{eq:border_costs}. The normals $n^{\textrm{in}}(s)$ and
$n^{\textrm{out}}(s)$ can easily be computed by differentiating the curves defining
the inner and outer borders. These costs are illustrated in
Fig.~\ref{fig:border}.

\paragraph{Constrained Control}
We constrain the steering angle to be between $[-\pi/3, \pi/3]$ by
parameterizing the steering angle as
\[
\delta(\tilde \delta) = \frac{2}{3} \arctan(\tilde \delta) \quad \mbox{for} \ \tilde \delta \in \reals.
\]
Similarly, we constrain the acceleration $a$ to be between $[c, d]$ (with
$c=-0.1, d=1.$), by parameterizing it as
\[
a (\tilde a) = (d-c)\operatorname{sig}(4 \tilde a /(d-c)) + c
\]
with $\operatorname{sig}: x \rightarrow 1/(1+e^{-x})$ the sigmoid function. The
final set of control variables is then $\tilde a, \tilde \delta, \alpha $.  

\paragraph{Control Cost}
For both trajectory costs, we add a square regularization on the control
variables of the system, i.e., the cost on the control variables is $ \lambda
\|\ctrl_t\|_2^2$ for some $\lambda \geq 0$ where $\ctrl_t$ are the control
variables at time $t$. 

\paragraph{Overall Contouring Cost}
The whole problem with contouring cost is then 
\begin{align}
	\min_{\substack{(\tilde a_0, \tilde \delta_0, \alpha_{0}), \ldots, (\tilde a_{\tau-1}, \tilde\delta_{\tau-1},\tilde \alpha_{\tau-1}}) } \quad & \sum_{t=0}^{\tau-1}
	\Big[ \rho_c e_c(x_t, y_t, s_t)^2 
	+ \rho_l e_l(x_t, y_t, s_t)^2  
	+ \rho_v \|v_{s, t} - v^{\textrm{ref}}\|_2^2 
	-\varepsilon\log(\nu_t) \nonumber\\
	& \qquad + \rho_b e_b(x_t, y_t, s_t)^2  
	+ \lambda(\tilde a_t^2 + \tilde\delta_t^2 + \alpha_{t}^2) \Big]\label{eq:total_contouring_cost}\\
	\textrm{s.t.} \quad & x_{t+1}, y_{t+1}, \yaw_{t+1}, v_{x, t+1}, v_{y, t+1}, \omega_{t+1} = \dyn(x_t, y_t, \yaw_t, v_{x, t}, v_{y, t}, \omega_t, \delta_t(\tilde \delta_t), a_t(\tilde a_t)) \nonumber\\
	&  s_{t+1} = s_t + \Delta \nu_{t}, \quad   \nu_{t+1} = \nu_{t} + \Delta \alpha_{t} \nonumber\\
&   z_0 = \hat z_0  \quad s_0 =0 \quad \nu_0 = v^{\textrm{ref}},\nonumber
\end{align}
with parameters $\rho_c, \rho_l, \rho_v, v^\textrm{ref}, \rho_b, \lambda$ and
$\dyn$ given in Eq.~\eqref{eq:real_dyn}.

\subsection{Numerical Constant}
The code is available at {\small \coderef}. We add for ease of reference, the
hyperparameters used for each setting.

\vspace{1em}

{\bf Pendulum}
\begin{enumerate}
	\item mass $m=1$, 
	\item gravitational constant $g=10$,
	\item length of the blob $l=1$,
	\item friction coefficient $\mu=0.01$,
	\item speed regularization $\lambda = 0.1$,
	\item control regularization $\rho=10^{-6}$,
	\item total time of the movement $T = 2$,  discretization step  $\Delta =
	T/\tau$ for varying $\tau$
	\item Euler discretization scheme.
\end{enumerate}

\vspace{1em}

{\bf Pendulum on a cart}
\begin{enumerate}
	\item mass of the rod $m=0.2$, 
	\item mass of the cart $M=0.5$,
	\item viscous coefficient $b=0.1$,
	\item moment of inertia $I =0.006$,
	\item length of the rod $0.3$,
	\item speed regularization $\lambda_1 = 0.1$,
	\item barrier parameter \revised{$\rho_2 = 10^{-6}.$},
	\item control regularization $\rho=10^{-6}$,
	\item total time of the movement $T = 2.5$, discretization step $\Delta =
	T/\tau$ for varying $\tau$,
	\item stay put time $\bar t = \horizon -\lfloor 0.6/\Delta\rfloor $,
	\item barriers $\bar z^+ =2$, $\bar z^- =-2$,
	\item Euler discretization scheme.
\end{enumerate}

\vspace{1em}

{\bf Simple car with tracking cost}
\begin{enumerate}
	\item length of the car $L=1$,
	\item reference speed $v^{ref}= 3$, 
	\item initial speed $v^{init} = 1$,
	\item control regularization $\lambda = 10^{-6}$, 
	\item total time of the movement $T = 2$, 
	\item simple track,
	\item Euler discretization scheme.
\end{enumerate}

\vspace{1em}

{\bf Bicycle model of a car with a contouring objective}
\begin{enumerate}
	\item $C_{m1} = 0.287$, $C_{m2} = 0.0545$,
	\item $C_{r0} = 0.0518$, $C_{rd} = 0.00035$,
	\item $B_r = 3.3852$, $C_r = 1.2691$, $D_r = 0.1737$, $l_r = 0.033$
	\item $B_f = 2.579$, $C_f =1.2$, $D_f = 0.192$, $l_f = 0.029$
	\item $m = 0.041$, $I_z = 27.8 \cdot 10^{-6}$
	\item contouring error penalty $\rho_c = 0.1$, 
	\item lagging error penalty $\rho_l =10$,
	\item reference speed penalty $\rho_v = 0.1$,
	\item barrier error penalty \revised{$\rho_b = 0.$,}
	\item reference speed $v^{\textrm{ref}} = 3$, 
	\item initial speed $v^{init} = 1$,
	\item control regularization $\lambda = 10^{-6}$, 
	\item total time of the movement $T=1$,
	\item simple track,
	\item Runge-Kutta discretization scheme.
\end{enumerate}

  \section{Additional Experiment}\label{app:exp_sup}
  \subsection{Time Comparison}
Figures~\ref{fig:conv_linquad_time} and~\ref{fig:conv_quad_time} present the
convergence of the algorithms presented in~\ref{fig:conv_linquad}
and~\ref{fig:conv_quad} in time rather than in iterations.

\subsection{Stepsize Selection}
\begin{figure}[t]
	\begin{center}
    \includegraphics[width=0.9\linewidth]{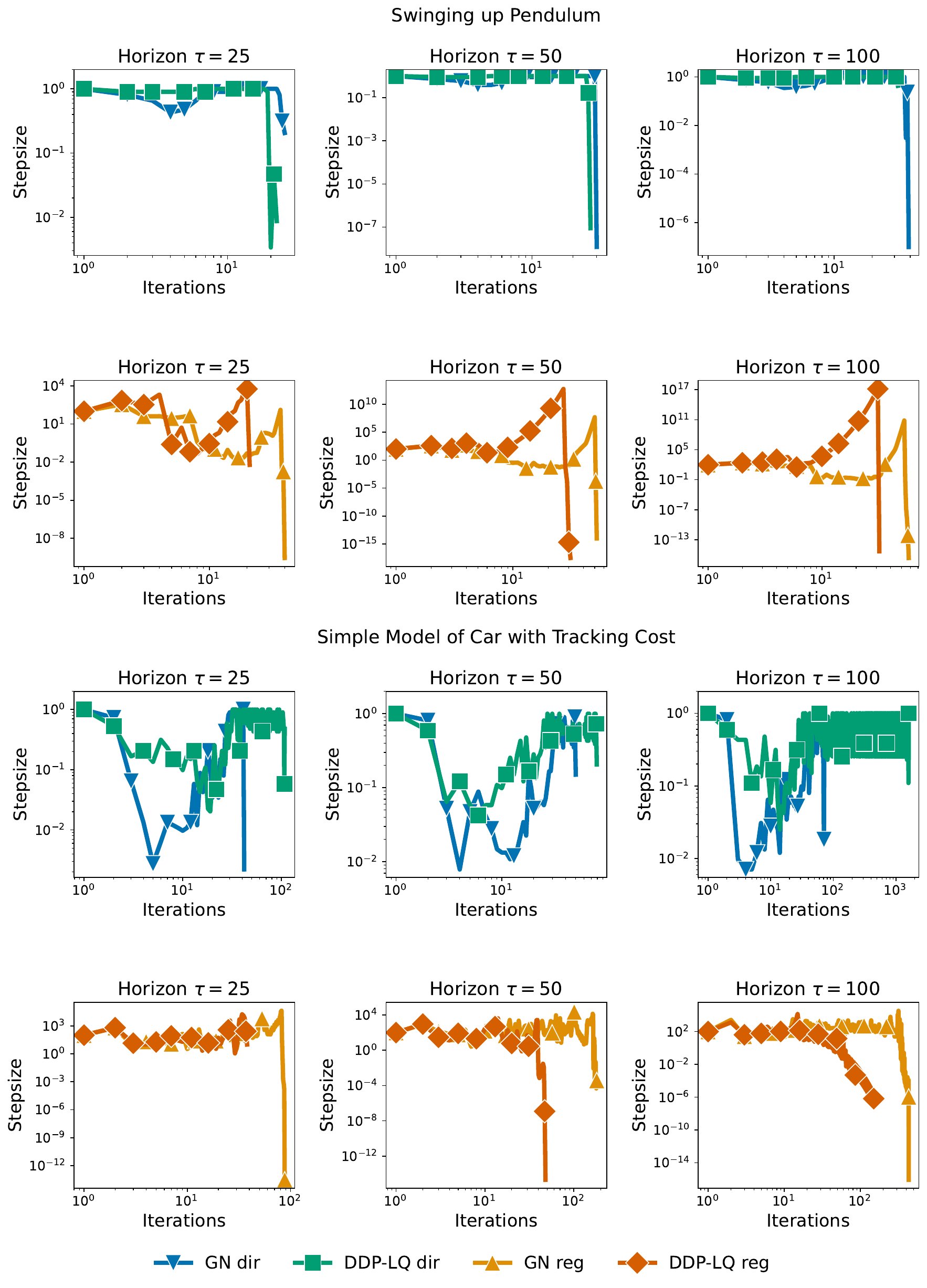}
    \caption{
    \revised{
		Stepsizes taken along the iterations on various control problems
		(see Appendix~\ref{app:exp_details}) with algorithms using linear (GD) or
		linear-quadratic approximations (GN, DDP-LQ, see Fig.~\ref{fig:taxonomy} for
		taxonomy details) and directional (dir \eqref{eq:linesearch_armijo}) or
    regularized (reg \eqref{eq:linesearch_reg}) steps.
    \label{fig:stepsize_linquad}}
    }
	\end{center}
\end{figure}

In Fig.~\ref{fig:stepsize_linquad}, we plot the stepsizes taken by algorithms
using linear-quadratic approximations for the pendulum and the simple
model of a car.
\revised{
\begin{enumerate}
	\item On the pendulum example, the stepsizes used by directional steps quickly
	tend to $1$ which means that the algorithms (GN or DDP-LQ) are then taking the
	largest possible stepsize for this strategy and may exhibit quadratic
	convergence. On the other hand, for the regularized steps, on the pendulum
	example, the regularization (i.e. the reciprocal of the stepsizes) quickly
	converges to $0$, which means that, as the number of iterations increases, the
	regularized and directional steps coincide. 
	\item For the car example, the step sizes for the directional steps slowly
	increase to one. We note yet that while stepsizes taken by DDP-LQ and GN are
	similar, DDP-LQ displays a faster convergence in terms of gradient norm. For
	the regularized steps, the regularizations (i.e. the reciprocal of the
	stepsizes) tend to remain low and stable. As for the directional steps, the
	regularizations taken by regularized steps are similar between DDP-LQ and GN,
  yet DDP-LQ displays faster convergence in terms of gradient norm.
\end{enumerate}
}

\begin{figure}[t]
	\begin{center}
    \includegraphics[width=0.9\linewidth]{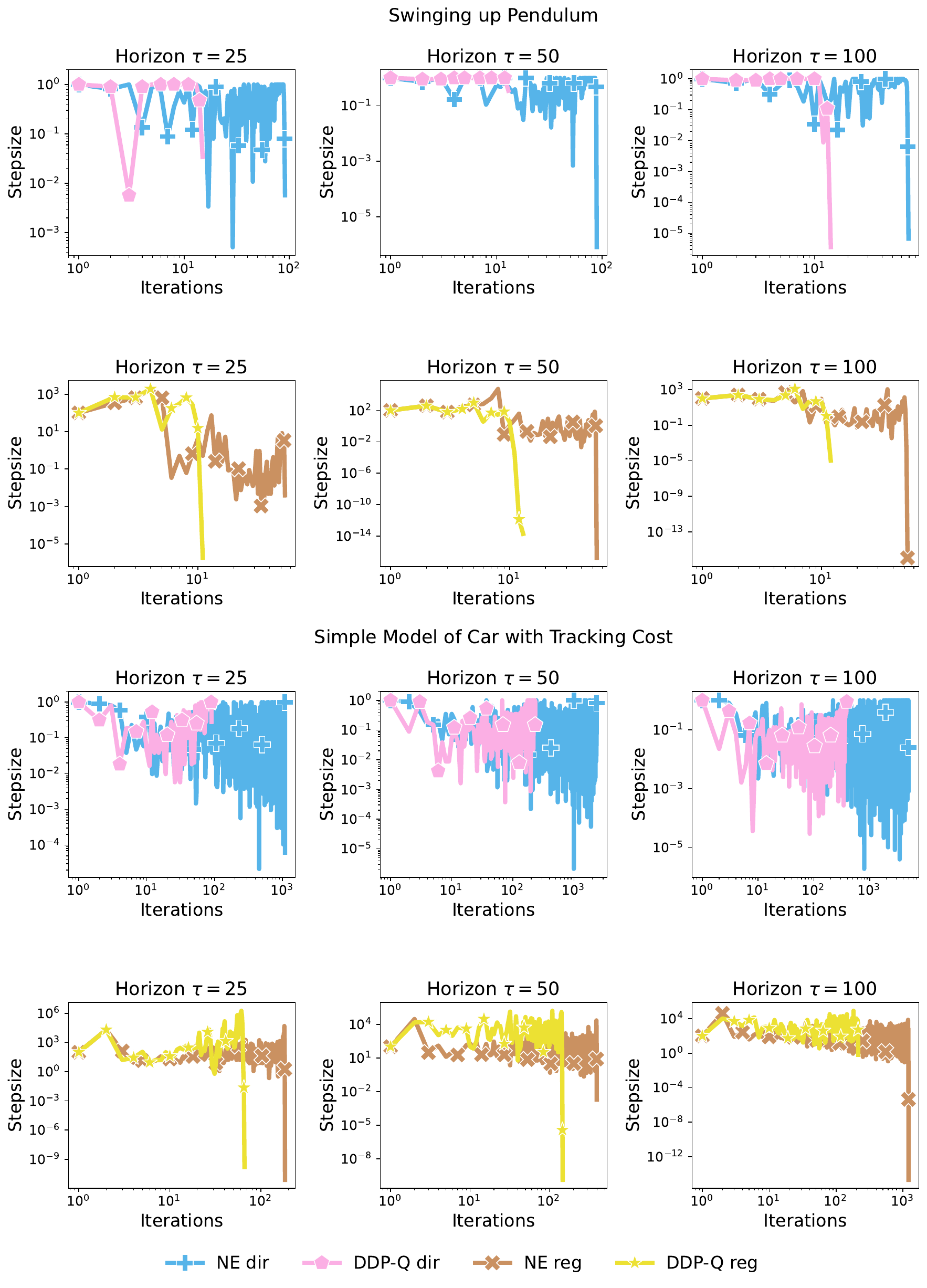}
    \caption{
    \revised{
		Stepsizes taken along the iterations on various control problems
    (see Appendix~\ref{app:exp_details}) with algorithms using quadratic
    approximations (NE, DDP-Q, see Fig.\ref{fig:taxonomy} for taxonomy details)
    and directional (dir \eqref{eq:linesearch_armijo}) or
    regularized (reg \eqref{eq:linesearch_reg}) steps.
    \label{fig:stepsize_quad}}
    }
	\end{center}
\end{figure}

In Fig.~\ref{fig:stepsize_quad}, we compare the stepsizes taken by the methods
using quadratic approximations.
\revised{
\begin{enumerate}
	\item In terms of directional steps, DDP-Q appears to take relatively large
	steps while its NE counterpart displays more variations on, e.g., the pendulum
	example. For both algorithms, the stepsizes tend to oscillate for the car
	example, never steadily taking full steps (stepsize of 1) closer to
	convergence.
	\item In terms of regularized steps, DDP-Q tends to take larger steps
  (smaller regularizations) than its Newton counterpart, in particular on the car
  example.
\end{enumerate}
}

\subsection{Comparison of Inner Solver}
\begin{figure}[t]
    \centering
    \includegraphics[width=0.8\linewidth]{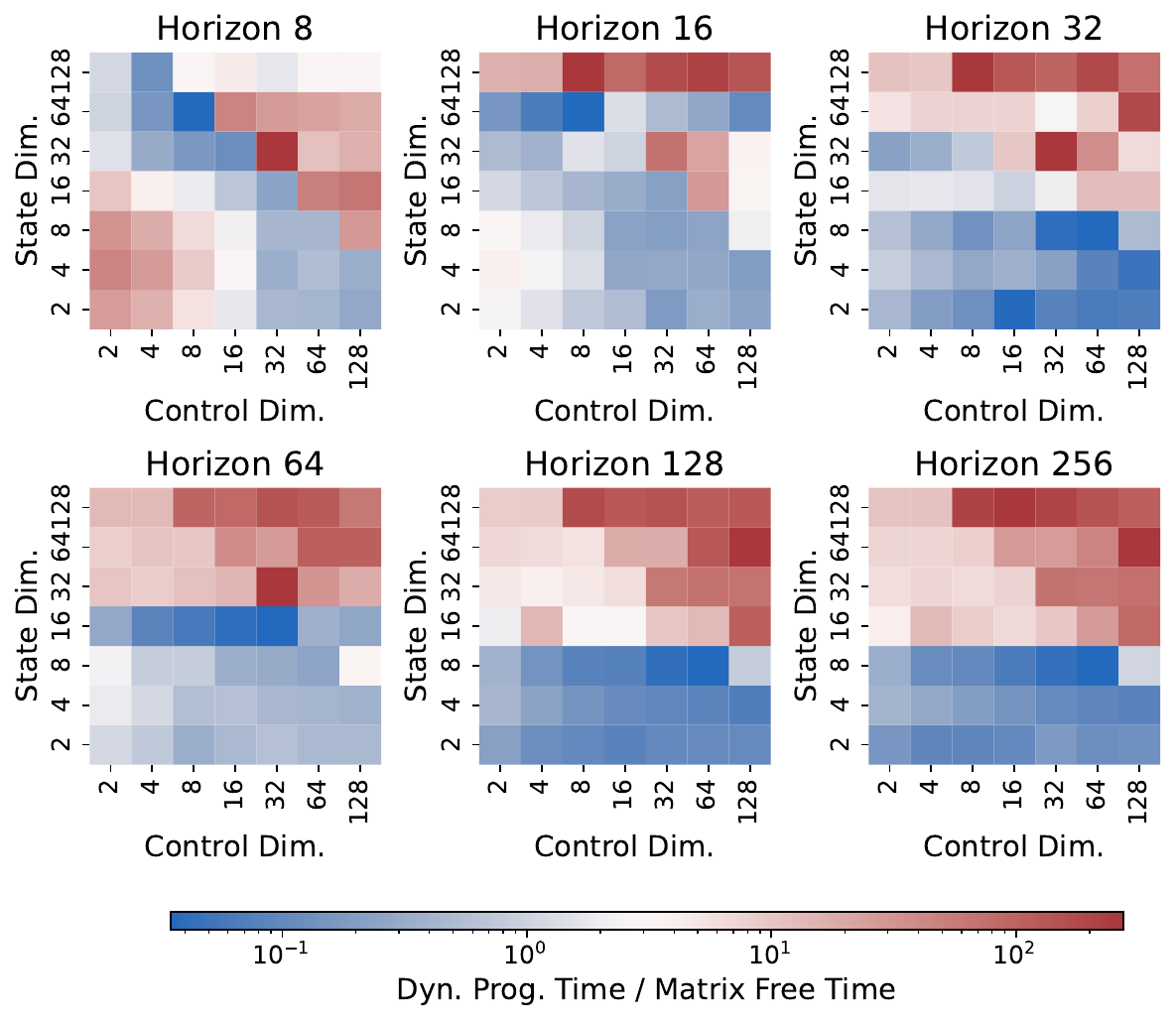}
    \caption{Comparison of time needed to solve synthetic linear quadratic
    control problems with either a matrix-free implementation or a dynamic
    programming implementation as presented in this manuscript. Blue cells
    indicate that dynamic programming is faster than matrix-free procedures,
    while red cells indicate the opposite. \label{fig:compa_sparse}}    
\end{figure}

As presented in Appendix~\ref{app:sparse_solvers}, we may consider using
directly Hessian-vector products to solve the linear quadratic controls arising
from the computation of Gauss-Newton and Newton steps. In
Fig.~\ref{fig:compa_sparse}, we plotted the ratio of time between an
implementation using dynamic programming and an implementation using matrix-free
solvers for varying dimensions of the state, the control and various horizons on
synthetic linear quadratic control problems. Namely, for each triplet
$(\dimstate, \dimctrl, \horizon)$, we generated five linear quadratic control
problems, solved each problem as if those were nonlinear dynamics for which we
are computing a Gauss-Newton step, by using each of the aforementioned methods.
We then averaged the time needed for each method over the five instances and
computed the ratio of time between an implementation by dynamic programming and
an implementation by matrix-free solvers. These values are recorded in a heatmap
in~\ref{fig:compa_sparse}, where blue cells correspond to instances where
dynamic programming is faster than the matrix-free program and red cells
correspond to instances where dynamic programming is slower than its
counterpart. 

The matrix-free solver approach can readily be implemented in any differentiable
programming framework such as CasADI~\citep{andersson2018cassadi} which takes
advantage of the differentiable dynamic programming framework to cast nonlinear
control problems as numerical optimization problems fed into off-the-shelf
solvers like IPOPT~\citep{wachter2006implementation}. 

The results presented in~\ref{fig:compa_sparse} show that for small state
dimensions, the dynamic programming approach is generally faster. As soon as the
state dimension exceeds a few dozen dimensions, the matrix-free approach is
generally faster. This observation matches the computation complexities
delineated in Appendix~\ref{app:sparse_solvers} and
Sec.~\ref{sec:comput_cplxity} as the matrix-free approach a priori scales
quadratically in terms of the state dimension while the dynamic programming
approach scales cubically. 

Beyond the time comparisons, each approach has different advantages. The
matrix-free approach enables simple introduction of constraints in control
variables by casting the whole problem as an optimization problem solved by
interior-point methods. On the other-hand, the dynamic programming approach can
be adapted to differential dynamic programming procedures as explained in this
manuscript.

\begin{figure}[t]
	\begin{center}
    \includegraphics[width=0.9\linewidth]{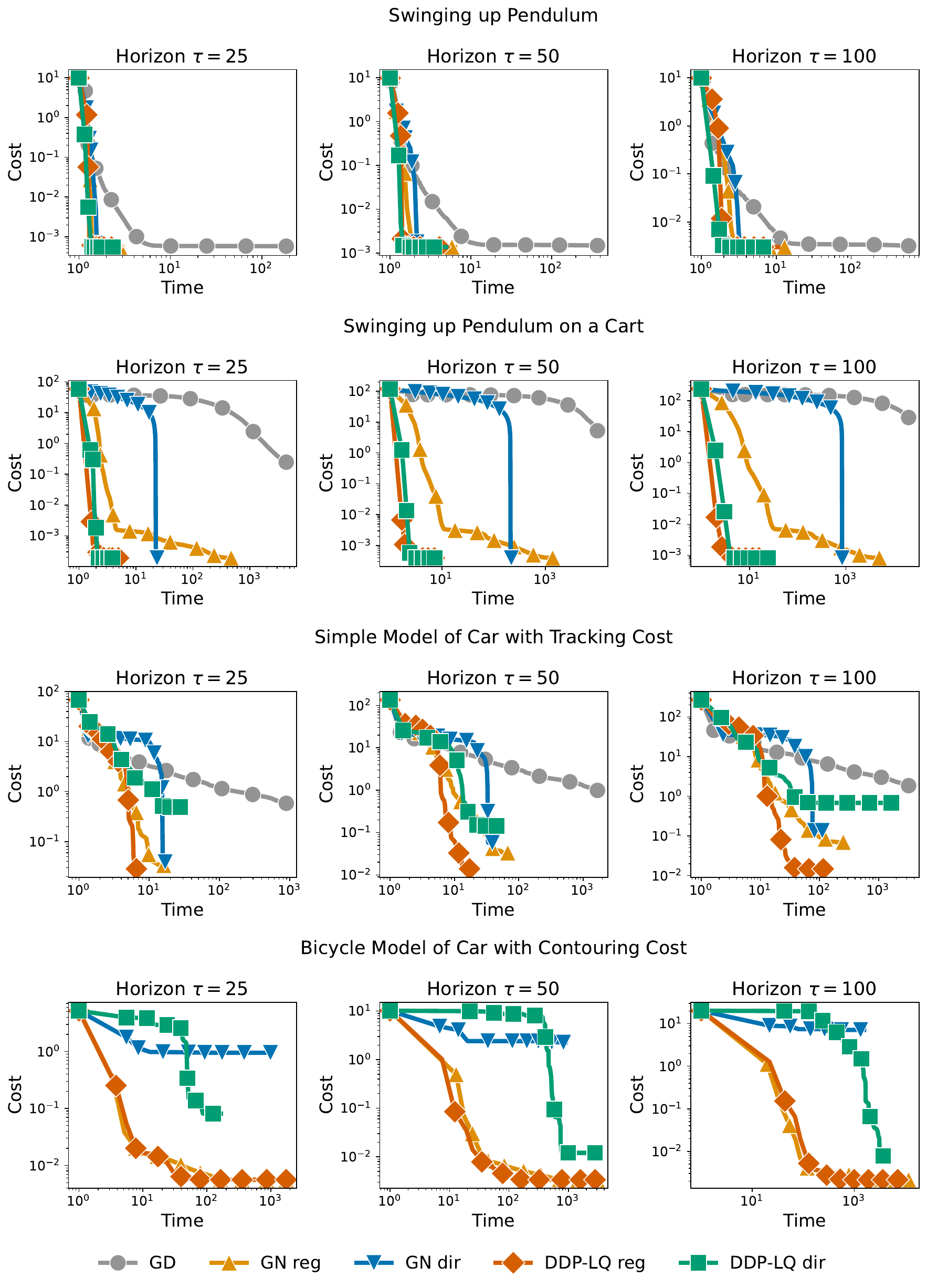}
    \caption{
    \revised{
		Cost along computational time on various control problems (see
    Appendix~\ref{app:exp_details}) with algorithms using linear (GD) or
    linear-quadratic approximations (GN, DDP-LQ, see Fig.~\ref{fig:taxonomy} for
    taxonomy details) and directional (dir \eqref{eq:linesearch_armijo}) or
    regularized (reg \eqref{eq:linesearch_reg}) steps.
    \label{fig:conv_linquad_time}}
    }
	\end{center}
\end{figure}

\begin{figure}[t]
	\begin{center}
    \includegraphics[width=0.9\linewidth]{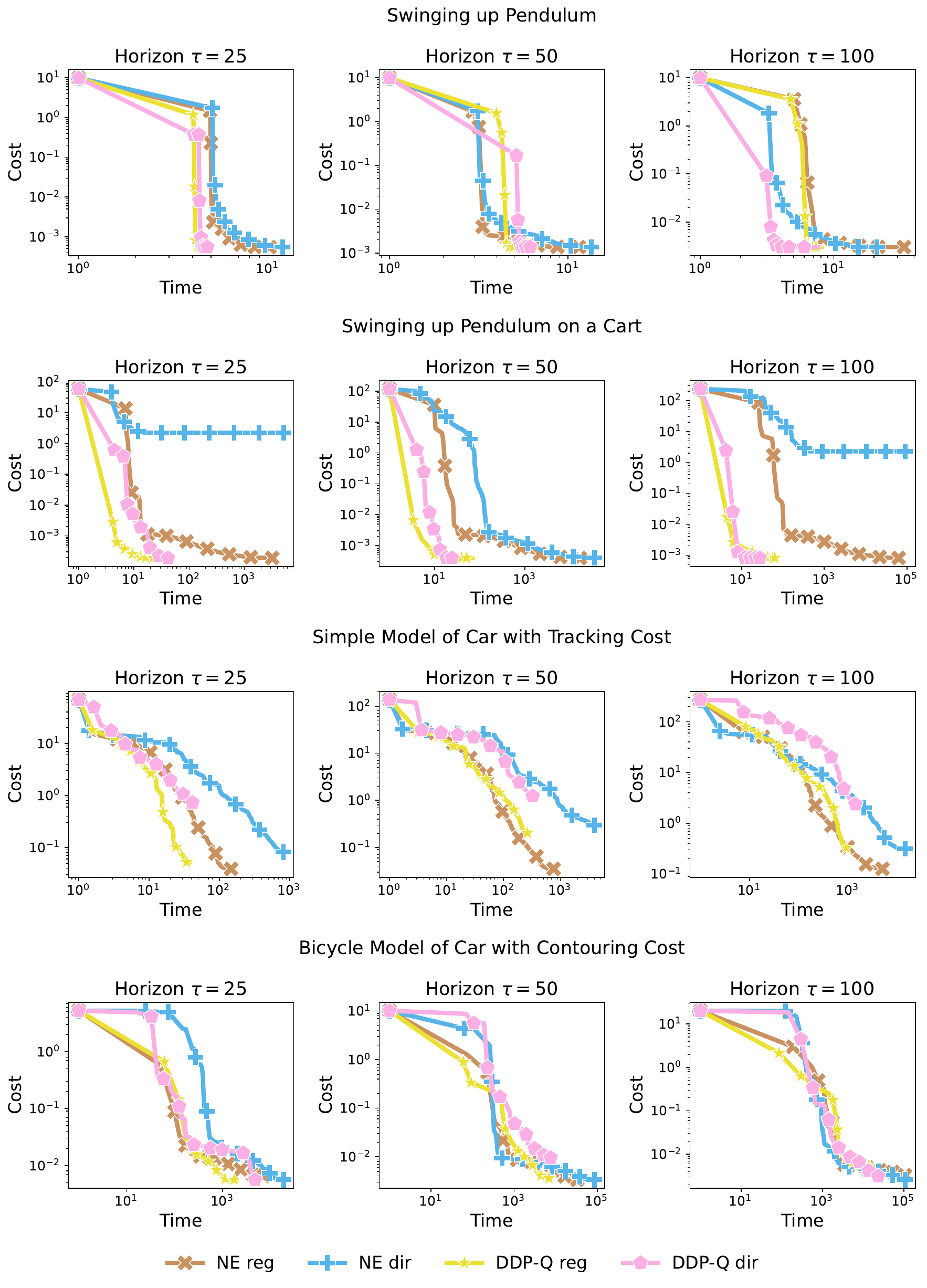}
    \caption{
    \revised{
		Cost along computational time on various control problems (see
		Appendix~\ref{app:exp_details}) with algorithms using quadratic
		approximations (NE, DDP-Q, see Fig.\ref{fig:taxonomy} for taxonomy details)
		and directional (dir \eqref{eq:linesearch_armijo}) or
    regularized (reg \eqref{eq:linesearch_reg}) steps.
		\label{fig:conv_quad_time}}
    }
	\end{center}
\end{figure}

We already presented in Sec.~\ref{sec:exp} numerical comparisons in
\emph{iterations} of classical optimization methods (Gauss-Newton or Newton)
against their differential dynamic programming counterparts (iLQR or DDP). By
using matrix-free solvers instead of dynamic programming procedures to implement
Gauss-Newton or Newton steps, these behaviors in iterations would not change.
The comparisons in time presented in Fig.~\ref{fig:conv_linquad_time} and
Fig.~\ref{fig:conv_quad_time} can change by using matrix-free solvers as
suggested by the heatmaps presented in Fig.~\ref{fig:compa_sparse}. However, the
qualitative conclusions presented in this manuscript remain the same and suggest
that differential dynamic programming methods may offer overall gains over
classical optimization algorithms.
 
\clearpage

\begin{figure}[t]
	\begin{center}
    \includegraphics[width=0.9\linewidth]{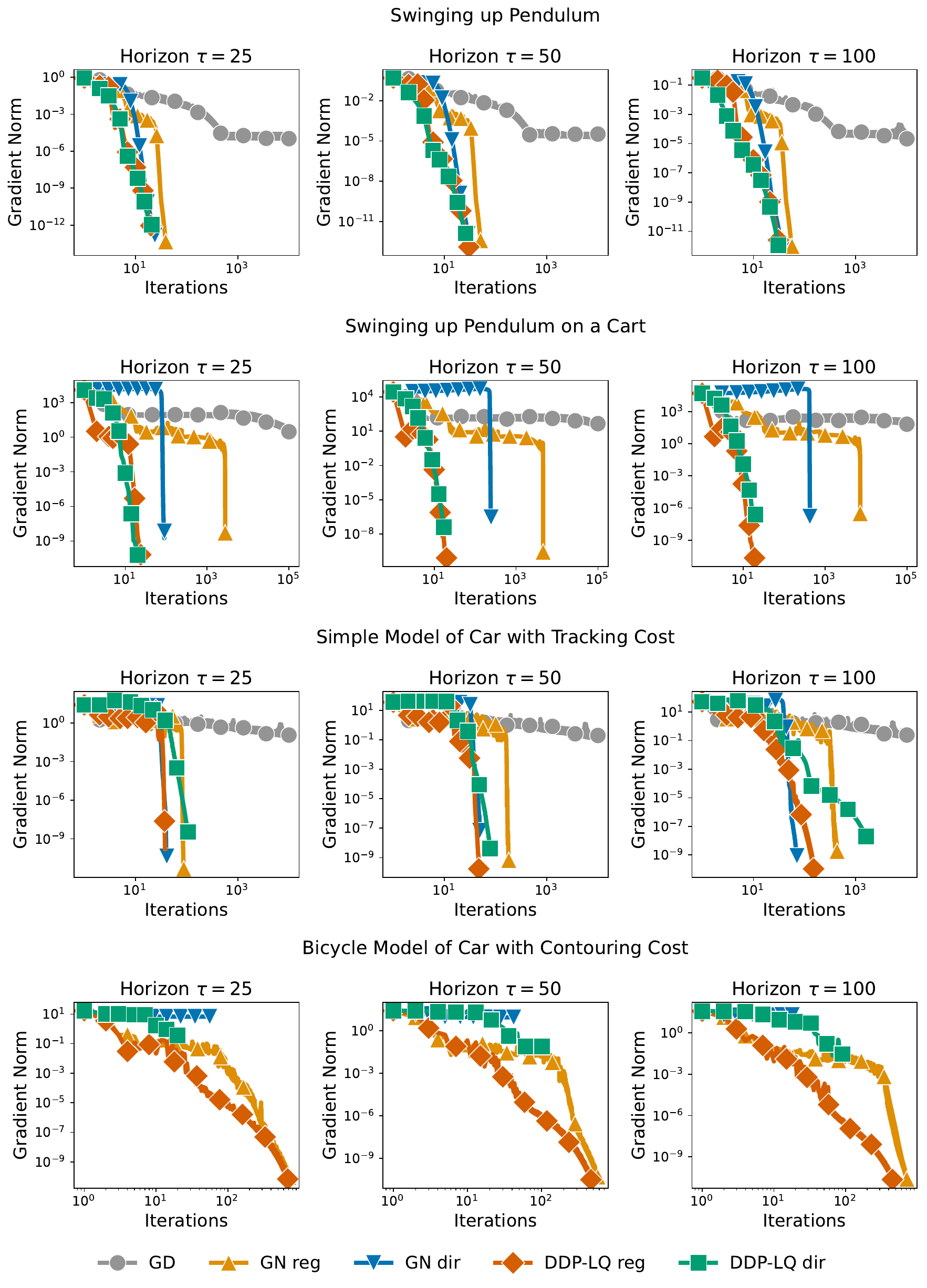}
    \caption{
    \revised{
		Gradient norm along iterations on various control problems (see
    Appendix~\ref{app:exp_details}) with algorithms using linear (GD) or
    linear-quadratic approximations (GN, DDP-LQ, see Fig.~\ref{fig:taxonomy} for
    taxonomy details) and directional (dir \eqref{eq:linesearch_armijo}) or
    regularized (reg \eqref{eq:linesearch_reg}) steps.
    \label{fig:norm_grad_linquad}}
    }
	\end{center}
\end{figure}

\begin{figure}[t]
	\begin{center}
    \includegraphics[width=0.9\linewidth]{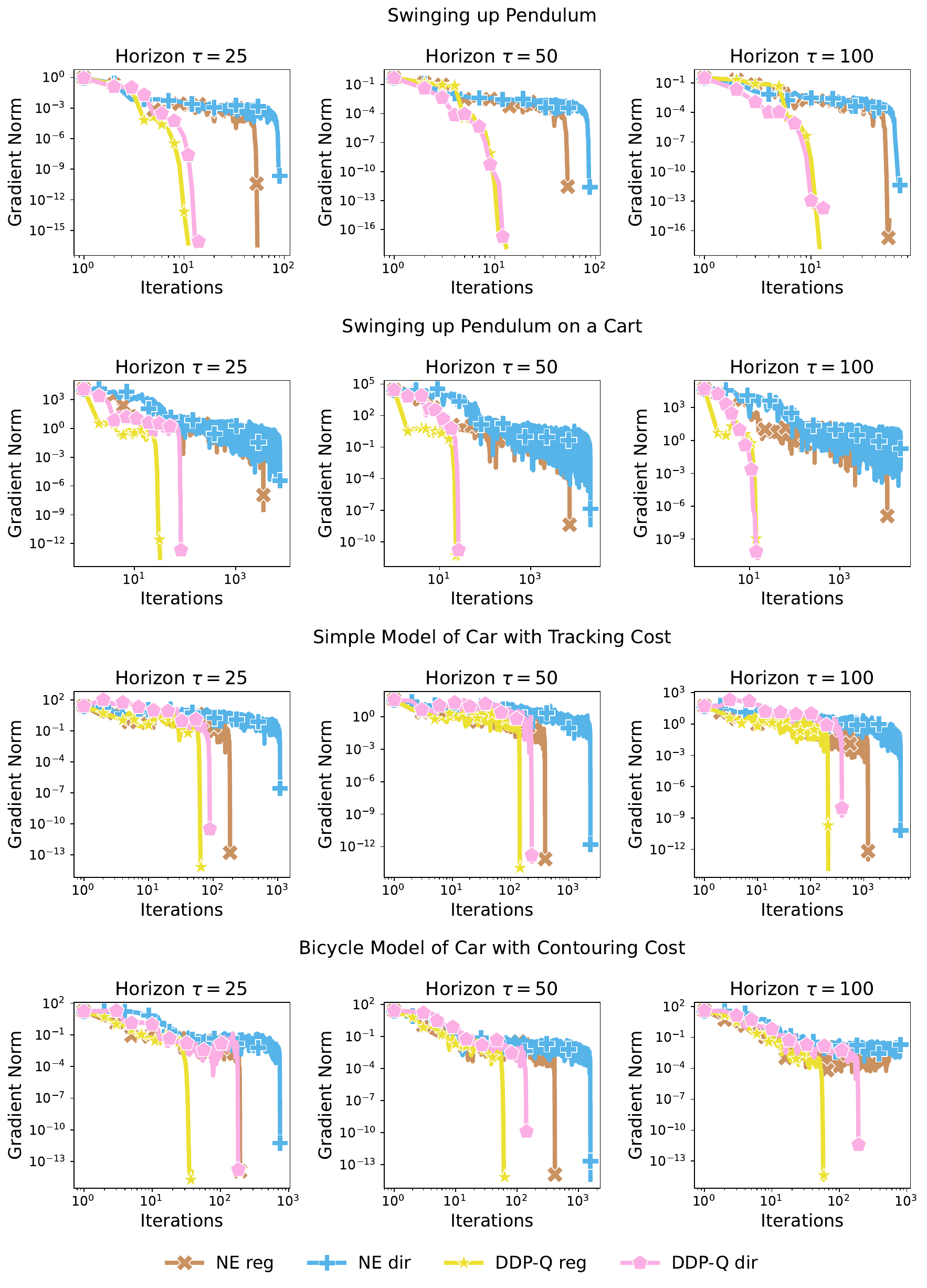}
    \caption{
    \revised{
		Gradient norm along iterations on various control problems (see
		Appendix~\ref{app:exp_details}) with algorithms using quadratic
		approximations (NE, DDP-Q, see Fig.\ref{fig:taxonomy} for taxonomy details)
		and directional (dir \eqref{eq:linesearch_armijo}) or
    regularized (reg \eqref{eq:linesearch_reg}) steps.
		\label{fig:norm_grad_quad}}
    }
	\end{center}
\end{figure}

\begin{figure}[t]
  \begin{center}
    \includegraphics[width=0.9\linewidth]{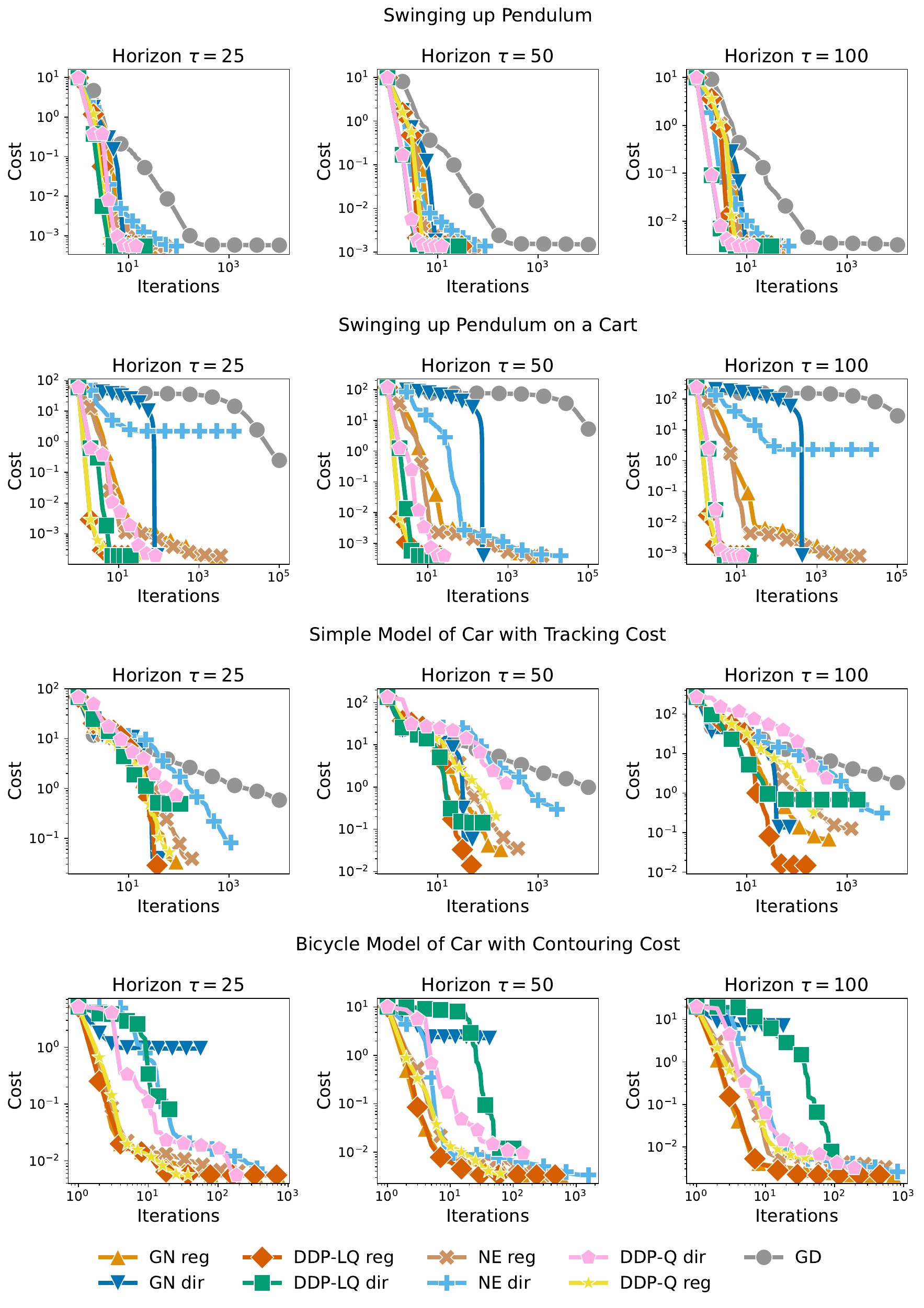}
    \caption{
    \revised{
    Cost along iterations on various control problems (see
    Appendix~\ref{app:exp_details}) with the algorithms presented in
    Fig.~\ref{fig:taxonomy} and directional (dir \eqref{eq:linesearch_armijo}) or
    regularized (reg \eqref{eq:linesearch_reg}) steps.
    \label{fig:conv_all_iteration}}
    }
  \end{center}
\end{figure}

\begin{figure}[t]
  \begin{center}
    \includegraphics[width=0.9\linewidth]{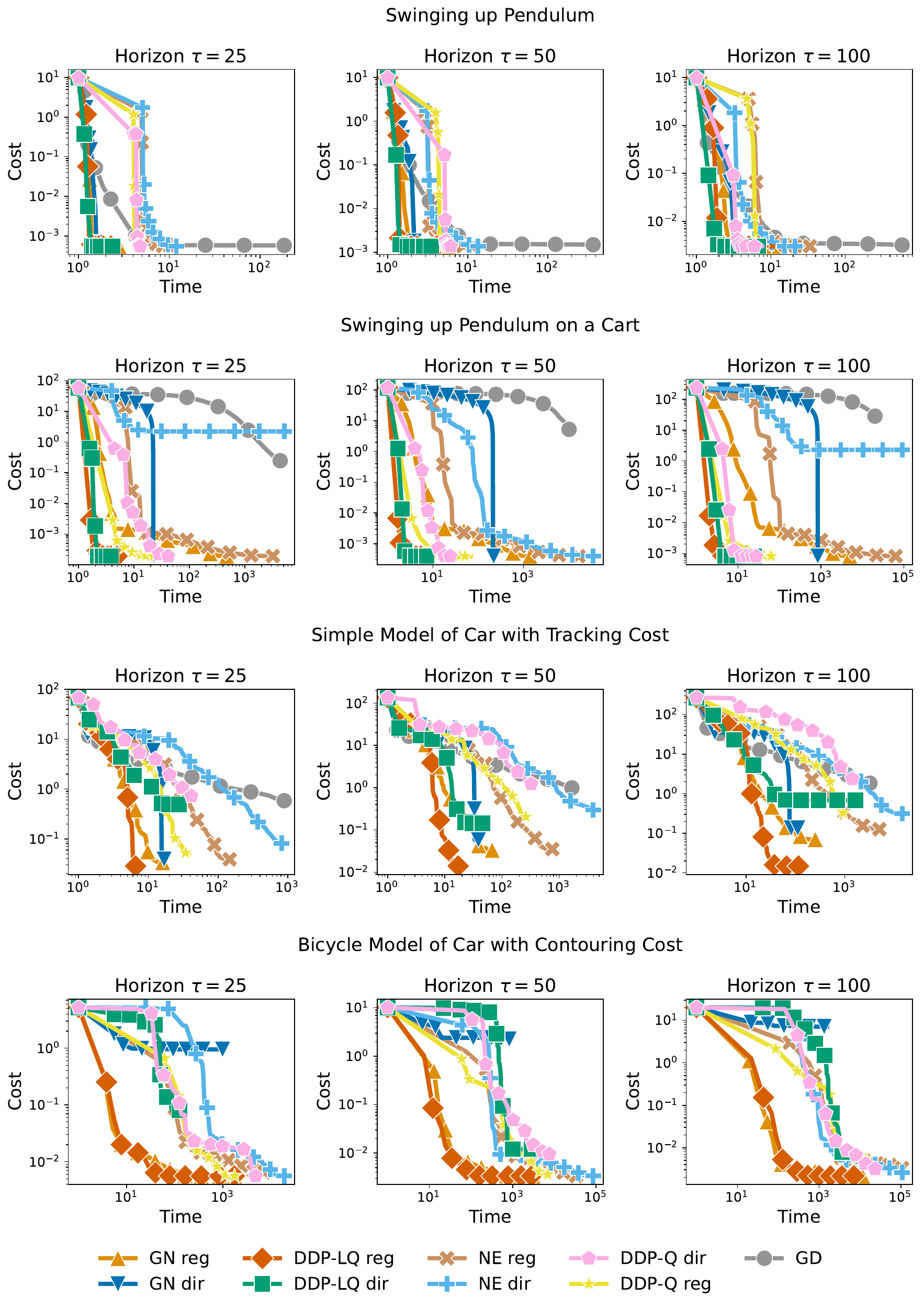}
    \caption{
    \revised{
    Cost along computational time on various control problems (see
    Appendix~\ref{app:exp_details}) with the algorithms presented in
    Fig.~\ref{fig:taxonomy} and directional (dir \eqref{eq:linesearch_armijo}) or
    regularized (reg \eqref{eq:linesearch_reg}) steps.
    \label{fig:conv_all_time}}
    }
  \end{center}
\end{figure}

  \clearpage
  \bibliography{ilqc_refs}

\begin{thebibliography}{}

\bibitem[Abadi et~al., 2015]{tensorflow2015-whitepaper}
Abadi, M., Agarwal, A., Barham, P., Brevdo, E., Chen, Z., Citro, C., Corrado,
  G.~S., Davis, A., Dean, J., Devin, M., Ghemawat, S., Goodfellow, I., Harp,
  A., Irving, G., Isard, M., Jia, Y., Jozefowicz, R., Kaiser, L., Kudlur, M.,
  Levenberg, J., Man\'{e}, D., Monga, R., Moore, S., Murray, D., Olah, C.,
  Schuster, M., Shlens, J., Steiner, B., Sutskever, I., Talwar, K., Tucker, P.,
  Vanhoucke, V., Vasudevan, V., Vi\'{e}gas, F., Vinyals, O., Warden, P.,
  Wattenberg, M., Wicke, M., Yu, Y., and Zheng, X. (2015).
\newblock {TensorFlow}: Large-scale machine learning on heterogeneous systems.

\bibitem[Andersson et~al., 2018]{andersson2018cassadi}
Andersson, J. A.~E., Gillis, J., Horn, G., Rawlings, J.~B., and Diehl, M.
  (2018).
\newblock {CasADi} -- {A} software framework for nonlinear optimization and
  optimal control.
\newblock {\em Mathematical Programming Computation}.

\bibitem[Bambade et~al., 2022]{bambade2022prox}
Bambade, A., El-Kazdadi, S., Taylor, A., and Carpentier, J. (2022).
\newblock Prox-qp: Yet another quadratic programming solver for robotics and
  beyond.
\newblock In {\em RSS 2022-Robotics: Science and Systems}.

\bibitem[Baur and Strassen, 1983]{baur1983complexity}
Baur, W. and Strassen, V. (1983).
\newblock The complexity of partial derivatives.
\newblock {\em Theoretical computer science}, 22(3):317--330.

\bibitem[Baydin et~al., 2018]{baydin2018automatic}
Baydin, A.~G., Pearlmutter, B.~A., Radul, A.~A., and Siskind, J.~M. (2018).
\newblock Automatic differentiation in machine learning: a survey.
\newblock {\em Journal of machine learning research}, 18(153):1--43.

\bibitem[Bellman, 1971]{bellman1971introduction}
Bellman, R. (1971).
\newblock {\em Introduction to the mathematical theory of control processes},
  volume~2.
\newblock Academic press.

\bibitem[Betts, 2010]{betts2010practical}
Betts, J. (2010).
\newblock {\em Practical methods for optimal control and estimation using
  nonlinear programming}.
\newblock SIAM.

\bibitem[Bock and Plitt, 1984]{bock1984multiple}
Bock, H.~G. and Plitt, K.-J. (1984).
\newblock A multiple shooting algorithm for direct solution of optimal control
  problems.
\newblock {\em IFAC Proceedings Volumes}, 17(2):1603--1608.

\bibitem[Bolte and Pauwels, 2020]{bolte2020mathematical}
Bolte, J. and Pauwels, E. (2020).
\newblock A mathematical model for automatic differentiation in machine
  learning.
\newblock In {\em Advances in Neural Information Processing Systems},
  volume~33.

\bibitem[Boyd and Vandenberghe, 1997]{boyd1997semidefinite}
Boyd, S. and Vandenberghe, L. (1997).
\newblock Semidefinite programming relaxations of non-convex problems in
  control and combinatorial optimization.
\newblock In {\em Communications, Computation, Control, and Signal Processing},
  pages 279--287. Springer.

\bibitem[Boyd and Vandenberghe, 2004]{boyd2004convex}
Boyd, S. and Vandenberghe, L. (2004).
\newblock {\em Convex optimization}.
\newblock Cambridge university press.

\bibitem[Bradbury et~al., 2018]{jax2018github}
Bradbury, J., Frostig, R., Hawkins, P., Johnson, M.~J., Leary, C., Maclaurin,
  D., Necula, G., Paszke, A., Vander{P}las, J., Wanderman-{M}ilne, S., and
  Zhang, Q. (2018).
\newblock {JAX}: composable transformations of {P}ython+{N}um{P}y programs.

\bibitem[Bynum et~al., 2021]{bynum2021pyomo}
Bynum, M.~L., Hackebeil, G.~A., Hart, W.~E., Laird, C.~D., Nicholson, B.~L.,
  Siirola, J.~D., Watson, J.-P., and Woodruff, D.~L. (2021).
\newblock {\em Pyomo--optimization modeling in python}, volume~67.
\newblock Springer Science \& Business Media, third edition.

\bibitem[Diehl et~al., 2006]{diehl2006fast}
Diehl, M., Bock, H.~G., Diedam, H., and Wieber, P.-B. (2006).
\newblock Fast direct multiple shooting algorithms for optimal robot control.
\newblock {\em Fast motions in biomechanics and robotics: optimization and
  feedback control}, pages 65--93.

\bibitem[Diehl et~al., 2009]{diehl2009efficient}
Diehl, M., Ferreau, H.~J., and Haverbeke, N. (2009).
\newblock Efficient numerical methods for nonlinear mpc and moving horizon
  estimation.
\newblock {\em Nonlinear model predictive control: towards new challenging
  applications}, pages 391--417.

\bibitem[Dunn and Bertsekas, 1989]{dunn1989efficient}
Dunn, J. and Bertsekas, D. (1989).
\newblock Efficient dynamic programming implementations of {Newton}'s method
  for unconstrained optimal control problems.
\newblock {\em Journal of Optimization Theory and Applications}, 63(1):23--38.

\bibitem[Dunning et~al., 2017]{dunning2017jump}
Dunning, I., Huchette, J., and Lubin, M. (2017).
\newblock Jump: A modeling language for mathematical optimization.
\newblock {\em SIAM review}, 59(2):295--320.

\bibitem[Farshidian et~al., 2017]{farshidian2017efficient}
Farshidian, F., Neunert, M., Winkler, A.~W., Rey, G., and Buchli, J. (2017).
\newblock An efficient optimal planning and control framework for quadrupedal
  locomotion.
\newblock In {\em 2017 IEEE International Conference on Robotics and Automation
  (ICRA)}, pages 93--100. IEEE.

\bibitem[Frasch et~al., 2015]{frasch2015parallel}
Frasch, J.~V., Sager, S., and Diehl, M. (2015).
\newblock A parallel quadratic programming method for dynamic optimization
  problems.
\newblock {\em Mathematical programming computation}, 7(3):289--329.

\bibitem[Giftthaler et~al., 2018]{giftthaler2018family}
Giftthaler, M., Neunert, M., St{\"a}uble, M., Buchli, J., and Diehl, M. (2018).
\newblock A family of iterative {Gauss-Newton} shooting methods for nonlinear
  optimal control.
\newblock In {\em 2018 IEEE/RSJ International Conference on Intelligent Robots
  and Systems (IROS)}, pages 1--9.

\bibitem[Gilbert, 1992]{gilbert1992automatic}
Gilbert, J.~C. (1992).
\newblock Automatic differentiation and iterative processes.
\newblock {\em Optimization methods and software}, 1(1):13--21.

\bibitem[Gill et~al., 2005]{gill2005snopt}
Gill, P.~E., Murray, W., and Saunders, M.~A. (2005).
\newblock Snopt: An sqp algorithm for large-scale constrained optimization.
\newblock {\em SIAM review}, 47(1):99--131.

\bibitem[Goodfellow et~al., 2016]{goodfellow2016deep}
Goodfellow, I., Bengio, Y., and Courville, A. (2016).
\newblock {\em Deep learning}.
\newblock MIT press.

\bibitem[Griewank and Walther, 2008]{griewank2008evaluating}
Griewank, A. and Walther, A. (2008).
\newblock {\em Evaluating derivatives: principles and techniques of algorithmic
  differentiation}.
\newblock SIAM.

\bibitem[Houska and Diehl, 2013]{houska2013quadratically}
Houska, B. and Diehl, M. (2013).
\newblock A quadratically convergent inexact sqp method for optimal control of
  differential algebraic equations.
\newblock {\em Optimal Control Applications and Methods}, 34(4):396--414.

\bibitem[Jacobson and Mayne, 1970]{jacobson1970differential}
Jacobson, D. and Mayne, D. (1970).
\newblock {\em Differential Dynamic Programming}.
\newblock Elsevier.

\bibitem[Jallet et~al., 2023]{jallet2023proxddp}
Jallet, W., Bambade, A., Arlaud, E., El-Kazdadi, S., Mansard, N., and
  Carpentier, J. (2023).
\newblock Proxddp: Proximal constrained trajectory optimization.

\bibitem[Kakade et~al., 2020]{kakade2020information}
Kakade, S., Krishnamurthy, A., Lowrey, K., Ohnishi, M., and Sun, W. (2020).
\newblock Information theoretic regret bounds for online nonlinear control.
\newblock {\em Advances in Neural Information Processing Systems},
  33:15312--15325.

\bibitem[LeCun, 1988]{lecun1988theoretical}
LeCun, Y. (1988).
\newblock A theoretical framework for back-propagation.
\newblock In {\em 1988 Connectionist Models Summer School, CMU, Pittsburg, PA}.

\bibitem[Li and Todorov, 2007]{li2007iterative}
Li, W. and Todorov, E. (2007).
\newblock Iterative linearization methods for approximately optimal control and
  estimation of non-linear stochastic system.
\newblock {\em International Journal of Control}, 80(9):1439--1453.

\bibitem[Liao and Shoemaker, 1991]{liao1991convergence}
Liao, L.-Z. and Shoemaker, C. (1991).
\newblock Convergence in unconstrained discrete-time differential dynamic
  programming.
\newblock {\em IEEE Transactions on Automatic Control}, 36(6):692--706.

\bibitem[Liao and Shoemaker, 1992]{liao1992advantages}
Liao, L.-Z. and Shoemaker, C.~A. (1992).
\newblock Advantages of differential dynamic programming over {N}ewton's method
  for discrete-time optimal control problems.
\newblock Technical report, Cornell University.

\bibitem[Liniger et~al., 2015]{liniger2015optimization}
Liniger, A., Domahidi, A., and Morari, M. (2015).
\newblock Optimization-based autonomous racing of 1: 43 scale {RC} cars.
\newblock {\em Optimal Control Applications and Methods}, 36(5):628--647.

\bibitem[Lions, 1982]{lions1982generalized}
Lions, P.-L. (1982).
\newblock {\em Generalized Solutions of Hamilton-Jacobi Equations}.
\newblock Pitman.

\bibitem[Magdy et~al., 2019]{magdy2019modeling}
Magdy, M., El~Marhomy, A., and Attia, M.~A. (2019).
\newblock Modeling of inverted pendulum system with gravitational search
  algorithm optimized controller.
\newblock {\em Ain Shams Engineering Journal}, 10(1):129--149.

\bibitem[Mayne and Polak, 1975]{mayne1975first}
Mayne, D. and Polak, E. (1975).
\newblock First-order strong variation algorithms for optimal control.
\newblock {\em Journal of Optimization Theory and Applications},
  16(3):277--301.

\bibitem[Messerer et~al., 2021]{messerer2021survey}
Messerer, F., Baumg{\"a}rtner, K., and Diehl, M. (2021).
\newblock Survey of sequential convex programming and generalized
  {Gauss-Newton} methods.
\newblock {\em ESAIM. Proceedings and Surveys}, 71:64.

\bibitem[Murray and Yakowitz, 1984]{murray1984differential}
Murray, D. and Yakowitz, S. (1984).
\newblock Differential dynamic programming and {Newton}'s method for discrete
  optimal control problems.
\newblock {\em Journal of Optimization Theory and Applications},
  43(3):395--414.

\bibitem[Nesterov, 2018]{nesterov2018lectures}
Nesterov, Y. (2018).
\newblock {\em Lectures on convex optimization}.
\newblock Springer.

\bibitem[Nganga and Wensing, 2021]{nganga2021accelerating}
Nganga, J. and Wensing, P. (2021).
\newblock Accelerating second-order differential dynamic programming for
  rigid-body systems.
\newblock {\em IEEE Robotics and Automation Letters}, 6(4):7659--7666.

\bibitem[Nocedal and Wright, 2006]{nocedal2006numerical}
Nocedal, J. and Wright, S. (2006).
\newblock {\em Numerical optimization}.
\newblock Springer Science \& Business Media.

\bibitem[Pantoja, 1988]{de1988differential}
Pantoja, J. (1988).
\newblock Differential dynamic programming and {Newton}'s method.
\newblock {\em International Journal of Control}, 47(5):1539--1553.

\bibitem[Paszke et~al., 2019]{paszke2017automatic}
Paszke, A., Gross, S., Massa, F., Lerer, A., Bradbury, J., Chanan, G., Killeen,
  T., Lin, Z., Gimelshein, N., Antiga, L., Desmaison, A., Kopf, A., Yang, E.,
  DeVito, Z., Raison, M., Tejani, A., Chilamkurthy, S., Steiner, B., Fang, L.,
  Bai, J., and Chintala, S. (2019).
\newblock Pytorch: An imperative style, high-performance deep learning library.
\newblock In {\em Advances in Neural Information Processing Systems},
  volume~32.

\bibitem[Polak, 1971]{polak1971computational}
Polak, E. (1971).
\newblock {\em Computational methods in optimization: a unified approach},
  volume~77.
\newblock Academic press.

\bibitem[Rao et~al., 1998]{rao1998application}
Rao, C., Wright, S., and Rawlings, J. (1998).
\newblock Application of interior-point methods to model predictive control.
\newblock {\em Journal of optimization theory and applications},
  99(3):723--757.

\bibitem[Recht, 2019]{recht2019tour}
Recht, B. (2019).
\newblock A tour of reinforcement learning: The view from continuous control.
\newblock {\em Annual Review of Control, Robotics, and Autonomous Systems},
  2:253--279.

\bibitem[Roulet et~al., 2019]{roulet2019iterative}
Roulet, V., Srinivasa, S., Drusvyatskiy, D., and Harchaoui, Z. (2019).
\newblock Iterative linearized control: stable algorithms and complexity
  guarantees.
\newblock In {\em Proceedings of the 36th International Conference on Machine
  Learning}, pages 5518--5527.

\bibitem[Rumelhart et~al., 1986]{rumelhart1985learning}
Rumelhart, D.~E., Hinton, G.~E., and Williams, R.~J. (1986).
\newblock Learning representations by back-propagating errors.
\newblock {\em Nature}, 323(6088):533--536.

\bibitem[Schmidhuber, 1990]{schmidhuber1990making}
Schmidhuber, J. (1990).
\newblock {\em Making the world differentiable: on using self supervised fully
  recurrent neural networks for dynamic reinforcement learning and planning in
  non-stationary environments}, volume 126.
\newblock Inst. f{\"u}r Informatik.

\bibitem[Sideris and Bobrow, 2005]{sideris2005efficient}
Sideris, A. and Bobrow, J. (2005).
\newblock An efficient sequential linear quadratic algorithm for solving
  nonlinear optimal control problems.
\newblock In {\em Proceedings of the 2005 American Control Conference}, pages
  2275--2280.

\bibitem[Srinivasan and Todorov, 2015]{srinivasan2015graphical}
Srinivasan, A. and Todorov, E. (2015).
\newblock Graphical newton.
\newblock {\em arXiv preprint arXiv:1508.00952}.

\bibitem[Tassa et~al., 2007]{tassa2007receding}
Tassa, Y., Erez, T., and Smart, W. (2007).
\newblock Receding horizon differential dynamic programming.
\newblock {\em Advances in neural information processing systems}, 20.

\bibitem[Tassa et~al., 2012]{tassa2012synthesis}
Tassa, Y., Erez, T., and Todorov, E. (2012).
\newblock Synthesis and stabilization of complex behaviors through online
  trajectory optimization.
\newblock In {\em 2012 IEEE/RSJ International Conference on Intelligent Robots
  and Systems}, pages 4906--4913.

\bibitem[Tassa et~al., 2014]{tassa2014control}
Tassa, Y., Mansard, N., and Todorov, E. (2014).
\newblock Control-limited differential dynamic programming.
\newblock In {\em 2014 IEEE International Conference on Robotics and Automation
  (ICRA)}, pages 1168--1175.

\bibitem[Todorov et~al., 2012]{todorov2012mujoco}
Todorov, E., Erez, T., and Tassa, Y. (2012).
\newblock Mujoco: A physics engine for model-based control.
\newblock In {\em International Conference on Intelligent Robots and Systems
  (IROS)}, pages 5026--5033. IEEE.

\bibitem[Verschueren et~al., 2021]{Verschueren2021}
Verschueren, R., Frison, G., Kouzoupis, D., Frey, J., van Duijkeren, N.,
  Zanelli, A., Novoselnik, B., Albin, T., Quirynen, R., and Diehl, M. (2021).
\newblock acados -- a modular open-source framework for fast embedded optimal
  control.
\newblock {\em Mathematical Programming Computation}.

\bibitem[Verschueren et~al., 2016]{verschueren2016exploiting}
Verschueren, R., van Duijkeren, N., Quirynen, R., and Diehl, M. (2016).
\newblock Exploiting convexity in direct optimal control: a sequential convex
  quadratic programming method.
\newblock In {\em 2016 IEEE 55th Conference on Decision and Control (CDC)},
  pages 1099--1104. IEEE.

\bibitem[Von~Stryk, 1993]{von1993numerical}
Von~Stryk, O. (1993).
\newblock {\em Numerical solution of optimal control problems by direct
  collocation}.
\newblock Springer.

\bibitem[W{\"a}chter and Biegler, 2006]{wachter2006implementation}
W{\"a}chter, A. and Biegler, L.~T. (2006).
\newblock On the implementation of an interior-point filter line-search
  algorithm for large-scale nonlinear programming.
\newblock {\em Mathematical programming}, 106:25--57.

\bibitem[Werbos, 1994]{werbos1994roots}
Werbos, P. (1994).
\newblock {\em The Roots of Backpropagation: From Ordered Derivatives to Neural
  Networks and Political Forecasting}.
\newblock Wiley-Interscience.

\bibitem[Wright, 1990]{wright1990solution}
Wright, S. (1990).
\newblock Solution of discrete-time optimal control problems on parallel
  computers.
\newblock {\em Parallel Computing}, 16(2-3):221--237.

\bibitem[Wright, 1991a]{wright1991partitioned}
Wright, S. (1991a).
\newblock Partitioned dynamic programming for optimal control.
\newblock {\em SIAM Journal on optimization}, 1(4):620--642.

\bibitem[Wright, 1991b]{wright1991structured}
Wright, S. (1991b).
\newblock Structured interior point methods for optimal control.
\newblock In {\em Proceedings of the 30th IEEE Conference on Decision and
  Control}, pages 1711--1716.

\bibitem[Wright, 1993]{wright1993interior}
Wright, S.~J. (1993).
\newblock Interior point methods for optimal control of discrete time systems.
\newblock {\em Journal of Optimization Theory and Applications},
  77(1):161--187.

\bibitem[Zhang et~al., 2023]{zhang2021dive}
Zhang, A., Lipton, Z.~C., Li, M., and Smola, A.~J. (2023).
\newblock {\em Dive into Deep Learning}.
\newblock Cambridge University Press.

\end{thebibliography}
	\bibliographystyle{apalike}
	
\end{document}